\pdfoutput=1
\RequirePackage{ifpdf}
\ifpdf 
\documentclass[pdftex]{sigma}
\else
\documentclass{sigma}
\fi

\newtheorem{Theorem}{Theorem}[section]
\newtheorem{Corollary}[Theorem]{Corollary}
\newtheorem{Lemma}[Theorem]{Lemma}
\newtheorem{Proposition}[Theorem]{Proposition}
 { \theoremstyle{definition}
\newtheorem{Definition}[Theorem]{Definition}
\newtheorem{Notation}[Theorem]{Notation}
\newtheorem{Example}[Theorem]{Example}
\newtheorem{Remark}[Theorem]{Remark} }

\numberwithin{equation}{section}

\newcommand{\CP}{\mathbb{C}P^1}
\newcommand{\MC}{\mathcal{X}_m}
\newcommand{\MCOV}{\mathcal{X}_m^{\text{ov}}}
\newcommand{\MCOVI}{\mathcal{X}_m^{\text{ov}-1}}

\begin{document}
\allowdisplaybreaks

\newcommand{\arXivNumber}{1912.00261}

\renewcommand{\PaperNumber}{005}

\FirstPageHeading

\ShortArticleName{The Ooguri--Vafa Space as a Moduli Space of Framed Wild Harmonic Bundles}

\ArticleName{The Ooguri--Vafa Space as a Moduli Space\\ of Framed Wild Harmonic Bundles}

\Author{Iv\'an TULLI}

\AuthorNameForHeading{I.~Tulli}

\Address{Department of Pure Mathematics, University of Sheffield, Sheffield, S3 7RH, UK}
\Email{\href{mailto:i.n.tulli@sheffield.ac.uk}{i.n.tulli@sheffield.ac.uk}}

\ArticleDates{Received April 08, 2024, in final form January 05, 2025; Published online January 14, 2025}

\Abstract{The Ooguri--Vafa space is a 4-dimensional incomplete hyperk\"ahler manifold, defined on the total space of a singular torus fibration with one singular nodal fiber. It has been proposed that the Ooguri--Vafa hyperk\"ahler metric should be part of the local model of the hyperk\"ahler metric of the Hitchin moduli spaces, near the most generic kind of singular locus of the Hitchin fibration. In order to relate the Ooguri--Vafa space with the Hitchin moduli spaces, we show that the Ooguri--Vafa space can be interpreted as a set of rank 2, framed wild harmonic bundles over $\mathbb{C}P^1$, with one irregular singularity. Along the way we show that a certain twistor family of holomorphic Darboux coordinates, which describes the hyperk\"ahler geometry of the Ooguri--Vafa space, has an interpretation in terms of Stokes data associated to our framed wild harmonic bundles.}

\Keywords{Higgs bundles; hyperk\"ahler geometry; twistor spaces; Stokes data}

\Classification{53C26; 53C28; 34M40}

\section{Introduction}

This work originated as an effort to understand, in an explicit way, the hyperk\"ahler metric of the Hitchin moduli spaces $\mathcal{M}_{\text{Hit}}$ associated to a (possibly punctured) compact Riemann surface $C$ \cite{BB04,Hitchin,Konno,Nakajima}. See also the foundational work \cite{Cor,Don,Sim,Sim92}. Roughly speaking, $\mathcal{M}_{\text{Hit}}$ parametrizes equivalence classes of harmonic bundles. The latter are tuples $\bigl(E,\overline{\partial}_E,\theta,h\bigr)$, where $\bigl(E,\overline{\partial}_E\bigr)\to C$ is a holomorphic bundle, $\theta$ is an endomorphism valued (1,0)-form $\theta\in \Omega^{(1,0)}(C,\operatorname{End}(E))$ known as the Higgs field, and $h$ is a hermitian metric on $E\to C$ such that the Hitchin equations are satisfied
\begin{equation}\label{hiteq}
 F\bigl(D\bigl(\overline{\partial}_E,h\bigr)\bigr)+\bigl[\theta,\theta^{\dagger_h}\bigr]=0 , \qquad \overline{\partial}_E(\theta)=0.
\end{equation}
In the above equation, $F\bigl(D\bigl(\overline{\partial}_E,h\bigr)\bigr)$ is the curvature of the Chern connection $D\bigl(\overline{\partial}_E,h\bigr)$ on \smash{$\bigl(E,\overline{\partial}_E,h\bigr)\to C$}, and $\dagger_{h}$ denotes the adjoint with respect to $h$. Furthermore, if $C$ is has punctures, then appropriate boundary conditions should be imposed at the punctures. In the original work of Hitchin \cite{Hitchin}, the equations \eqref{hiteq} were obtained by dimensional reduction to 2 dimensions of the self-dual solutions of the Yang--Mills equations over $\mathbb{R}^4$. The hyperk\"ahler structure was then obtained by interpreting $\mathcal{M}_{\text{Hit}}$ as a hyperk\"ahler quotient of an infinite-dimensional hyperk\"ahler affine space by the action of an infinite-dimensional gauge group, where the hyperk\"ahler moment map is determined by the left hand side of the two equations in \eqref{hiteq}.

The moduli space $\mathcal{M}_{\text{Hit}}$ also carries the structure of a complex integrable system $\pi\colon\mathcal{M}_{\text{Hit}}\to \mathcal{B}_{\text{Hit}}$, known as the Hitchin integrable system. Namely, with respect to one of the complex structures of the hyperk\"ahler structure of $\mathcal{M}_{\text{Hit}}$, $\pi\colon\mathcal{M}_{\text{Hit}}\to \mathcal{B}_{\text{Hit}}$ is a holomorphic fibration, $\mathcal{M}_{\text{Hit}}$ is holomorphic symplectic, and the generic fibers of $\pi$ are compact Lagrangian tori. There is a divisor $\mathcal{B}_{\text{sing}}\subset \mathcal{B}_{\text{Hit}}$, where the fibers become singular.

Besides their purely mathematical interest, the Hitchin moduli spaces arise in several contexts in the physics literature (see, for example, the references mentioned at the beginning of \cite{neitzke2014hitchin}). Of special interest to us is its appearance in a certain class of 4d $\mathcal{N}=2$ supersymmetric field theories, called theories of class $\mathcal{S}$ \cite{GMN1,GMN2,neitzke2014hitchin}. More precisely, they arise as moduli spaces associated to theories of class $\mathcal{S}$ compactified on a circle $S^1$. From this picture several conjectures about the behaviour of the hyperk\"ahler metric on $\mathcal{M}_{\text{Hit}}$ arise:
\begin{itemize}\itemsep=0pt
 \item One one hand, it is expected that away from $\pi^{-1}(\mathcal{B}_{\text{sing}})\subset \mathcal{M}_{\text{Hit}}$, the metric can be asymptotically approximated by a certain explicit and simpler semiflat\footnote{Here by semiflat we mean that the metric restricted to the fibers of $\pi$ becomes flat. It is also known as the rigid c-map metric \cite{SK}.} hyperk\"ahler metric, up to exponentially suppressed corrections. For a more precise statement of the conjecture and works where this has been proved in certain cases, see \cite{DN,LF,FMSW,MSWW,mochizuki2023asymptotic}. If $\mathcal{B}_{\text{reg}}:=\mathcal{B}-\mathcal{B}_{\text{sing}}$, then the semiflat metric is defined on $\pi^{-1}(\mathcal{B}_{\text{reg}})$ and it is purely determined by the affine special K\"ahler geometry on $\mathcal{B}_{\text{reg}}$ \cite{LF}.
 \item On the other hand, near the most generic singular locus of $\pi^{-1}(\mathcal{B}_{\text{sing}})\subset \mathcal{M}_{\text{Hit}}$ it is conjectured that the Ooguri--Vafa hyperk\"ahler metric should be part of the local model for its approximate description. The Ooguri--Vafa metric is a 4-dimensional incomplete hyperk\"aher manifold, defined on the total space of a singular torus fibration with a single nodal fiber. A more precise statement of the conjecture is given in \cite[Section 7]{A2}, while in \cite{GMN1,GW, OV} the Ooguri--Vafa metric is discussed in detail.
\end{itemize}

We remark that this proposed picture of the hyperk\"ahler metric is very similar to the one given by Gross--Wilson for the hyperk\"ahler metric of K$3$ surfaces (see \cite{GW}). In the picture of~\cite{GW}, we have a generic elliptic fibration of a K$3$ surface $f\colon X\to \CP$ with $24$ singular nodal fibers; the hyperk\"ahler metric of $X$ is then approximated by taking the semiflat metric away from the singular fibers, and gluing in the Ooguri--Vafa metric in a neighborhood of each singular fiber.

Motivated by these facts, our goal in this paper is to relate the Ooguri--Vafa space with the objects present in $\mathcal{M}_{\text{Hit}}$, namely harmonic bundles $\bigl(E,\overline{\partial}_E,\theta,h\bigr)$. We will find an interpretation of the Ooguri--Vafa space in terms of a set $\mathfrak{X}^{{\rm fr}}$ of equivalence classes of certain framed wild harmonic bundles $\bigl(E,\overline{\partial}_E,\theta,h,g\bigr)$. Roughly speaking, the connection between the two is done as follows:
\begin{itemize}\itemsep=0pt
 \item On one hand, the hyperk\"ahler structure of the Ooguri--Vafa space $\mathcal{M}^{\text{ov}}$ can be encoded in its associated twistor space $\mathcal{Z}^{\text{ov}}=\mathcal{M}^{\text{ov}}\times \mathbb{C}P^1$ \cite{SUSY}. In particular, $\mathcal{M}^{\text{ov}}$ has an associated twistor family of holomorphic symplectic forms $\Omega^{\text{ov}}(\xi)$, $\xi\in \mathbb{C}P^1$. In~\cite{GMN1}, the twistor family~$\Omega^{\text{ov}}(\xi)$ is encoded via a twistor family of holomorphic Darboux coordinates $\log(\mathcal{X}_{e}^{\text{ov}}(\xi))$ and $\log(\mathcal{X}_{m}^{\text{ov}}(\xi))$. The family of complex coordinates $\mathcal{X}_e^{\text{ov}}(\xi)$ and $\mathcal{X}_{m}^{\text{ov}}(\xi)$ are referred to as the ``electric'' and ``magnetic'' twistor coordinates, and encode the hyperk\"ahler structure of $\mathcal{M}^{\text{ov}}$.
 \item On the other hand, given a framed wild harmonic bundle $\bigl(E,\overline{\partial}_E,\theta,h,g\bigr)$ we will associate a twistor family of framed flat bundles, which in turn has an associated twistor family of Stokes data (or ``generalized monodromy data''). From the Stokes data, we will define analogous coordinates $\mathcal{X}_e(\xi)$ and~$\mathcal{X}_{m}(\xi)$ for our set of framed wild harmonic bundles~$\mathfrak{X}^{{\rm fr}}$, and show that $\mathcal{X}_{e}(\xi)=\mathcal{X}_{e}^{\text{ov}}(\xi)$ and $\mathcal{X}_{m}(\xi)=\mathcal{X}_{m}^{\text{ov}}(\xi)$ under an appropriate correspondence of parameters. This last fact was anticipated in \cite[Section 9.4.3]{GMN2}, and will let us interpret~$\mathcal{M}^{\text{ov}}$ as a moduli space of framed wild harmonic bundles.
\end{itemize}

Our main result is the following.

\begin{Theorem}\label{maintheorem} Let $\mathcal{M}^{\textnormal{ov}}$ be the Ooguri--Vafa space with cut-off $\Lambda=4{\rm i}$, and let $\mathcal{B}\subset \mathbb{C}$ be the base of the singular torus fibration $\mathcal{M}^{\textnormal{ov}}\to \mathcal{B}$. Fixing an affine coordinate $z\in \mathbb{C}\subset \CP$, let
\begin{equation*}
 \mathfrak{X}^{\textnormal{fr}}(\mathcal{B}):=\bigl\{\bigl[E,\overline{\partial}_E,\theta,h,g\bigr]\in \mathfrak{X}^{\textnormal{fr}} \big| \textnormal{Det}(\theta)=-\bigl(z^2+2m\bigr){\rm d}z^2 \ \text{with}\ -2{\rm i}m \in \mathcal{B} \bigr\},
\end{equation*}
and let $\mathfrak{X}^{\textnormal{fr}}_*(\mathcal{B})\subset \mathfrak{X}^{\textnormal{fr}}(\mathcal{B})$ be the subset of elements with $m\neq 0$.
Then there is a one-to-one correspondence between $\mathfrak{X}^{\textnormal{fr}}(\mathcal{B})$ and $\mathcal{M}^{\textnormal{ov}}$ such that $\mathcal{X}_e=\mathcal{X}_e^{\text{ov}}$ and $\mathcal{X}_m=\mathcal{X}_m^{\text{ov}}$. Under this correspondence $\mathfrak{X}^{\textnormal{fr}}(\mathcal{B})$ gets an induced hyperk\"ahler structure, whose twistor family of holomorphic symplectic forms $\Omega(\xi)$ restricted to $\mathfrak{X}^{\textnormal{fr}}_*(\mathcal{B})$ is described by
\begin{equation*}
 \Omega(\xi)=\frac{{\rm d}\mathcal{X}_e(\xi)}{\mathcal{X}_e(\xi)}\wedge \frac{{\rm d}\MC(\xi)}{\MC(\xi)} \qquad \text{for}\quad \xi \in \mathbb{C}^*.
\end{equation*}
\end{Theorem}

Since our identification will use framed wild harmonic bundles, we remark that this set does not match any of the usual wild moduli spaces $\mathcal{M}_{\text{Hit}}$. In the usual story of moduli spaces of wild harmonic bundles over a punctured compact Riemann surface, one fixes the singular part of the Higgs field and the parabolic structure at the punctures. Under certain stability conditions, one obtains moduli spaces of these objects, with the natural hyperk\"ahler metric \cite{BB04}. On the other hand, in our set of wild harmonic bundles we will allow the simple pole of the Higgs field and the parabolic structure to vary. Furthermore, we will have the additional data of a ``framing''. Hence our moduli space must a priori be different from the usual moduli spaces of wild harmonic bundles. We remark that the moduli space of logarithmic connections equipped with framings on an $n$-pointed Riemann surface has been constructed in \cite{biswas2021modulispacesframedlogarithmic} as a Deligne--Mumford stack, so an extension of their results to the irregular case should be useful for the present case.

We hope that our interpretation of the Ooguri--Vafa space in terms of wild harmonic bundles serves as a first step to establish part of the conjectural picture of \cite{GMN1,GMN2} and \cite[Section 7]{A2} mentioned above. For now, we leave the question of the specific relation between the Ooguri--Vafa metric and the hyperk\"ahler metric of the Hitchin moduli spaces for future work. Independently of this problem, we also hope that our construction and methods can be generalized to produce hyperk\"ahler structures for similar sets of wild harmonic bundles.

\subsection{Summary and strategy of the paper}

We start Section~\ref{OVsec} by defining the Ooguri--Vafa hyperk\"ahler space \cite{OV} and describe its hyperk\"ahler structure following \cite{GMN1}. This space is built using the so called ``Gibbons--Hawking ansatz'' \cite{GHansatz}. This ansatz takes a positive harmonic function on an open set $U\subset \mathbb{R}^3$ and, provided some integrality condition is satisfied, produces a~principal ${\rm U}(1)$-bundle $X\to U$ with connection, whose total space carries a hyperk\"ahler metric. In our particular case, this principal ${\rm U}(1)$-bundle will have an extra $\mathbb{Z}$-shift symmetry. Upon dividing by this symmetry, we~obtain a~principal ${\rm U}(1)$-bundle of the form \smash{$\widetilde{\pi}\colon \widetilde{X}\to \bigl(\mathcal{B}\times S^1\bigr)\setminus \bigl(\{0\}^2\times\{1\}\bigr)$}, where $\mathcal{B}$ is an open subset of $\mathbb{C}$ containing the origin. By adding a point to $\widetilde{X}$, its hyperk\"ahler structure extends over $\{0\}^2\times \{1\}$, and we call such a space the Ooguri--Vafa space $\mathcal{M}^{\text{ov}}$. Strictly speaking, the definition of $\mathcal{M}^{\text{ov}}$ also depends on a choice of ``cut-off parameter'' $\Lambda \in \mathbb{C}^*$, but we omit this point until Section \ref{OVsec}.

\begin{figure}[ht]\centering
\includegraphics[width=6.0cm]{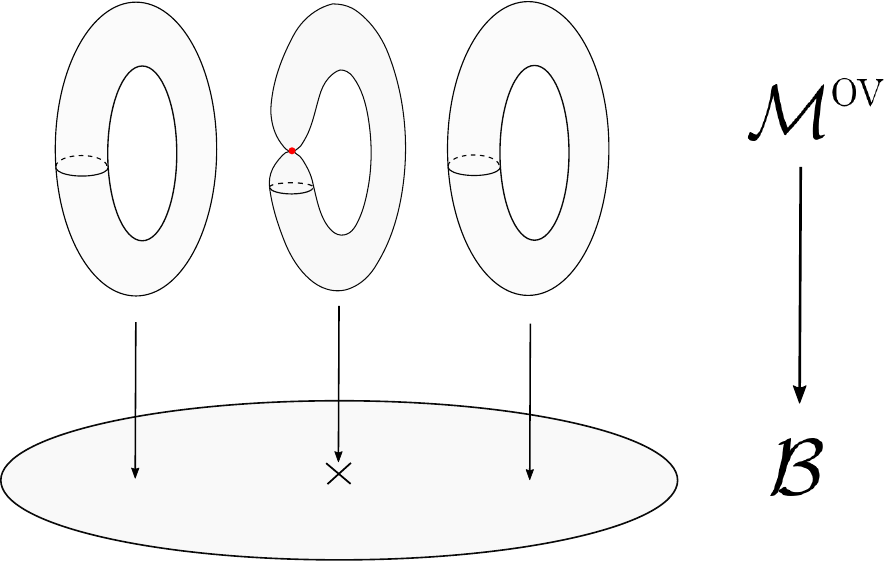}
\caption{$\mathcal{M}^{\textnormal{ov}}$ as a singular torus fibration over $\mathcal{B}\subset \mathbb{C}$. The central fiber at $0\in \mathcal{B}$ degenerates to a~torus with a node.}\label{fig1}
\end{figure}

From the above discussion, we can also think of $\mathcal{M}^{\text{ov}}$ as having a projection $p\colon \mathcal{M}^{\text{ov}} \to \mathcal{B} \subset \mathbb{C}$, making it a (singular) torus fibration. More precisely, for points $z\in \mathcal{B}\cap \mathbb{C}^*$ we have that $p^{-1}(z)$ is a torus, while $p^{-1}(0)$ is a torus with a node (see Figure~\ref{fig1}).

Since $\mathcal{M}^{\text{ov}}$ is hyperk\"ahler, it comes with a twistor family of holomorphic symplectic forms\footnote{In more global terms, let $\mathcal{Z}=\mathcal{M}^{\text{ov}}\times \CP$ be the associated twistor space, $\pi\colon\mathcal{Z}\to \CP$ the canonical projection into the second factor, and $T_{F}=\text{Ker}({\rm d}\pi)$ vertical bundle of $\pi$. Then the family of holomorphic symplectic forms $\Omega^{\text{ov}}$ gives a holomorphic section of the vector bundle $\wedge^2 T^*_{F}\otimes \pi^*\mathcal{O}(2)\to \mathcal{Z}$ (see \cite{SUSY}).} $\Omega^{\text{ov}}(\xi)$ for $\xi \in \mathbb{C}\subset\CP$. In \cite{GMN1}, this family is described away from the central fiber $p^{-1}(0)$ via the ``electric'' and ``magnetic'' twistor coordinates $\mathcal{X}_e^{\text{ov}}(\xi)$ and $\mathcal{X}_m^{\text{ov}}(\xi)$, satisfying
\begin{equation}\label{holsymform}
 \Omega^{\text{ov}}(\xi)=\frac{{\rm d}\mathcal{X}_e^{\text{ov}}(\xi)}{\mathcal{X}_e^{\text{ov}}(\xi)}\wedge \frac{{\rm d}\mathcal{X}_m^{\text{ov}}(\xi)}{\mathcal{X}_m^{\text{ov}}(\xi)} \qquad \text{for}\quad \xi \in \mathbb{C}^*,
\end{equation}
where ${\rm d}$ does not differentiate in the twistor parameter $\xi$. The coordinates $\mathcal{X}_e^{\text{ov}}(\xi)$ and $\mathcal{X}_m^{\text{ov}}(\xi)$ encode the hyperk\"ahler structure of $\mathcal{M}^{\text{ov}}$. In particular, $\log(\mathcal{X}_e^{\text{ov}}(\xi))$ and $\log(\mathcal{X}_m^{\text{ov}}(\xi))$ are a twistor family of holomorphic Darboux coordinates for the twistor family of holomorphic symplectic forms.

While $\mathcal{X}_e^{\text{ov}}(\xi)$ is holomorphic in $\xi \in \mathbb{C}^*$, $\mathcal{X}_m^{\text{ov}}(\xi)$ is only holomorphic in $\xi$ away from certain rays that depend on $z\in \mathcal{B}\setminus \{0\}$. More precisely, if we fix $z\in \mathcal{B}\setminus \{0\}$ and let
\begin{equation*}
 l_{\pm}(z):=\left\{\xi \in \mathbb{C}^* \mid \pm\frac{z}{\xi}<0\right\},
\end{equation*}
we then have that $\mathcal{X}_m^{\text{ov}}(\xi)$ is holomorphic in $\xi$ away from $l_{\pm}(z)$, and furthermore it has the following jumps
\begin{gather*}
 \mathcal{X}_m^{\text{ov}}(\xi)^{+}=\mathcal{X}_m^{\text{ov}}(\xi)^{-}(1-\mathcal{X}_e^{\text{ov}}(\xi))^{-1} \qquad\text{along} \quad \xi \in l_{+}(z),\\
 \mathcal{X}_m^{\text{ov}}(\xi)^{+}=\mathcal{X}_m^{\text{ov}}(\xi)^{-}\bigl(1-\mathcal{X}_e^{\text{ov}}(\xi)^{-1}\bigr) \qquad \text{along} \quad \xi \in l_{-}(z),
\end{gather*}
where the $\pm$ denotes the fact that we approach $l_{\pm}(z)$ clockwise or anticlockwise, respectively. Notice that even though $\mathcal{X}_m^{\text{ov}}(\xi)$ jumps along $l_{\pm}(z)$, the form of the jumps implies that $\Omega^{\text{ov}}(\xi)$ given by \eqref{holsymform} is continuous in $\xi \in \mathbb{C}^*$.

Looking forward, we remark that the jumps of $\mathcal{X}_m^{\text{ov}}(\xi)$ together with the asymptotics with respect to the twistor parameter as $\xi \to 0$ and as $\xi \to \infty$, will form our guiding principle for building the corresponding coordinate $\mathcal{X}_m(\xi)$ in the context of framed wild harmonic bundles. The plan is to set up a Riemann--Hilbert problem like the one in \cite{GMN1}, and use the uniqueness of solutions of such a problem to claim that $\mathcal{X}_m(\xi)=\mathcal{X}_m^{\text{ov}}(\xi)$, under an appropriate correspondence of parameters between the framed wild harmonic bundles and the Ooguri--Vafa space. The coordinate $\mathcal{X}_e(\xi)$ analogous to $\mathcal{X}_e^{\text{ov}}(\xi)$ will be easier to reproduce, and no Riemann--Hilbert problem will be necessary.

Once we finish with the necessary details of the Ooguri--Vafa space, we start Section \ref{sec3} by recalling the notions of unramified filtered Higgs bundles, unramified filtered flat bundles, unramified wild harmonic bundles, and the main results relating them. The notion of filtered bundle goes back to C. Simpson \cite{Sim}, and is equivalent to the notion of parabolic bundle due to Seshadri \cite{ParBun,ParBun2}. Many of the notations and conventions that we will use are the same as in~\cite{M}.

After recalling the basic definitions, we will define our set of ``framed'' wild harmonic bundles. Roughly speaking our set consists of tuples $\bigl(E,\overline{\partial}_E,\theta,h,g\bigr)$, where $(E,h)\to \mathbb{C}P^1$ is an~${\rm SU}(2)$ bundle, \smash{$\bigl(E|_{\mathbb{C}},\overline{\partial}_E,\theta,h\bigr)\to \mathbb{C}\subset \mathbb{C}P^1$} is a wild harmonic bundle, and $g$ is a unitary frame at~${\infty\in \mathbb{C}P^1}$.

The strategy will then be the following:
\begin{itemize}\itemsep=0pt
 \item To each $\bigl(E,\overline{\partial}_E,\theta,h,g\bigr)$ and $\xi \in \mathbb{C}^*$, we will associate a ``framed filtered flat bundle'' that we will denote by \smash{$\bigl(\mathcal{P}^h_*\mathcal{E}^{\xi},\nabla^{\xi}, \tau_{*}^{\xi}\bigr)\to \bigl(\CP,\infty\bigr)$}. We will give a more precise definition of this object later in Section \ref{sec3}. For now, it should be thought as a collection $\bigl(\mathcal{P}^h_a\mathcal{E}^{\xi},\nabla^{\xi}, \tau_{a}^{\xi}\bigr)\to \bigl(\CP,\infty\bigr)$ indexed by $a \in \mathbb{R}$, where $\mathcal{P}^h_a\mathcal{E}^{\xi} \to \CP$ is a holomorphic bundle, $\nabla^{\xi}$ is a~meromorphic (and hence flat) connection with a pole at $\infty$ given by
 \begin{equation*}
 \nabla^{\xi}=D\bigl(\overline{\partial}_E,h\bigr)+\xi^{-1}\theta + \xi \theta^{\dagger_h},
 \end{equation*}
 and $\tau_a^{\xi}$ is a certain frame of $\mathcal{P}^h_a\mathcal{E}^{\xi}|_{\infty}$ (see \eqref{singform}).

 \item To each such $\bigl(\mathcal{P}^h_*\mathcal{E}^{\xi},\nabla^{\xi}, \tau_{*}^{\xi}\bigr)\to \bigl(\CP,\infty\bigr)$, we will associate its Stokes data. Roughly speaking, the Stokes data will consist of transition functions between certain sectorial flat frames of $\nabla^{\xi}$ near the singularity (also known as the Stokes matrices), and the ``formal monodromy'' of $\nabla^{\xi}$ around the singularity (i.e., the monodromy of the formal diagonalization near $\infty$). We remark that Stokes data in this sense is usually associated to ``compatibly framed meromorphic flat bundles'' as in \cite{B}. We will show that in our case, the Stokes data can be actually associated to the framed filtered flat bundles that we consider. A~reference for the subject of Stokes data can be found in the classical book \cite{Wasow}. We will follow mainly~\cite{Bir,B,B2}.
 \item We will then construct the analog $\mathcal{X}_e(\xi)$ of $\mathcal{X}_e^{\text{ov}}(\xi)$ by taking the formal monodromy of the associated framed filtered flat bundle. Furthermore, we will construct the analog $\mathcal{X}_m(\xi)$ of~$\mathcal{X}_m^{\text{ov}}(\xi)$ in terms of the non-trivial entries of certain Stokes matrices. We will define~$\mathcal{X}_m(\xi)$ in such a way that we get the same jumping behaviour as $\mathcal{X}_m^{\text{ov}}(\xi)$. The quantities $\mathcal{X}_e(\xi)$ and $\mathcal{X}_m(\xi)$ will later be shown to be coordinates on $\mathfrak{X}^{{\rm fr}}$.
\end{itemize}

As we can see from the description so far, there is a heavy emphasis on the fact that all of our objects are framed. One of the reasons for taking framed objects, is to achieve that the non-trivial Stokes matrix entries are actual coordinates on the isomorphism classes of our objects. Without the framing, only the Stokes data up to conjugation by diagonal matrices is well defined on isomorphism classes. We should also remark that the idea of using these types of framed wild harmonic bundles, and using Stokes data to build the twistor coordinates, is heavily inspired by the observations made in \cite[Section 9.4.3]{GMN2}.

The rest of Section \ref{sec3} is concerned with the holomorphic dependence, the asymptotics, and jumping behavior of $\mathcal{X}_e$ and $\mathcal{X}_m$, which will be needed to match them with $\mathcal{X}_e^{\text{ov}}$ and $\mathcal{X}_m^{\text{ov}}$. More precisely,
\begin{itemize}\itemsep=0pt
 \item First, we must show that the coordinates $\mathcal{X}_e(\xi)$ and $\mathcal{X}_m(\xi)$ that we built out of Stokes data depend holomorphically on the twistor parameter $\xi$. This will be clear for $\mathcal{X}_e(\xi)$ due to the holomorphic dependence of the formal monodromy on $\xi$ see \eqref{formmon}. However, since the formal type of $\nabla^{\xi}$ with respect to the frame $\tau_{a}^{\xi}$ has anti-holomorphic dependence in~$\xi$ (see \eqref{formtype}), a more involved argument will be needed to show holomorphic dependence of~$\mathcal{X}_{m}(\xi)$. We will require to do a ``deformation of irregular values'' before being able to glue the family of framed filtered flat bundles into a meromorphic family in $\xi \in \mathbb{C}^*$. Roughly speaking, the procedure of deformation of irregular values varies the compatibly framed flat bundle, while keeping the Stokes data the same (and hence also goes by the name of isomonodromy). The idea of ``deformation of irregular values'' goes back to Jimbo--Miwa--Ueno \cite{JMU}. A reference where this is applied to the twistor family of meromorphic bundles associated to a wild harmonic bundle can be found in \cite[Chapter 9]{M}.

 \item Second, we compute the asymptotics of $\mathcal{X}_e(\xi)$ and $\mathcal{X}_m(\xi)$ as $\xi \to 0$ and as $\xi \to \infty$. While this computation does not require fancy machinery, it will require some results about asymptotic formulas of solutions to the parallel transport equation corresponding to $\nabla^{\xi}$. For this part, we mainly use techniques that can be found in \cite{Lev48,Vag10, Wasow}.
\end{itemize}

Finally, in Section \ref{sec4}, we explain the correspondence between the isomorphism classes of our set of framed wild harmonic bundles and the Ooguri--Vafa space. We will give here a brief description of how this correspondence works.

We will denote the set of isomorphism classes of our set of framed wild harmonic bundles by~$\mathfrak{X}^{{\rm fr}}$. For these classes, there is the parameter $m\in \mathbb{C}$ specifying the simple pole term of the Higgs field at $\infty\in \mathbb{C}P^1$, and a parameter $m^{(3)}\in \bigl(-\frac12,\frac12\bigr]\subset \mathbb{R}$ specifying the parabolic structure of the associated filtered Higgs bundles. If $\mathfrak{X}^{{\rm fr}}\bigl(m,m^{(3)}\bigr)\subset \mathfrak{X}^{{\rm fr}}$ denotes the set of isomorphism classes with associated parameters $m$ and \smash{$m^{(3)}$}, then in Proposition \ref{U1A} and Lemma~\ref{U1A2} we show that $\mathfrak{X}^{{\rm fr}}\bigl(m,m^{(3)}\bigr)$ is a ${\rm U}(1)$ torsor as long as $m$ and $m^{(3)}$ are not both $0$. When $m=m^{(3)}=0$, we show that \smash{$\mathfrak{X}^{{\rm fr}}\bigl(m,m^{(3)}\bigr)$} reduces to a point. The ${\rm U}(1)$ torsor structure basically comes from the framing of our objects.

Focusing on the case $m\neq 0$, the asymptotics of $\mathcal{X}_m(\xi)$ will then give us a natural way to locally trivialize the torsors. These local trivializations will have, after an appropriate correspondence of parameters, the same transition functions as the ${\rm U}(1)$ principal bundle appearing in the construction of the Ooguri--Vafa space. By ``appropriate correspondence of parameters'' we mean the following: if $\bigl(z=x^1+{\rm i}x^2,{\rm e}^{2\pi {\rm i} x^3}\bigr)$ denotes the canonical coordinates on $\mathcal{B}\times S^1$, then we have that $z$ corresponds to $-2{\rm i}m$ and $x^{3}$ corresponds to \smash{$m^{(3)}$}. The reason for this correspondence of parameters will arise naturally from the specific formulas of the twistor coordinates. Furthermore, we will see that in the case $m=0$ we have the same picture as the singular fiber of $\mathcal{M}^{\text{ov}}$ from Figure~\ref{fig1}: $\mathfrak{X}^{{\rm fr}}\bigl(0,m^{(3)}\bigr)$ is a ${\rm U}(1)$ torsor for $m^{(3)}\neq 0$, and degenerates into a point for $m^{(3)}=0$. Hence, we are able to identify $\mathcal{M}^\text{ov}$ with the subset of the elements of $\mathfrak{X}^{{\rm fr}}$ having associated parameter $m$ satisfying the condition $-2{\rm i}m \in \mathcal{B}$. Our main Theorem \ref{maintheorem} will then follow.

\section{The Ooguri--Vafa space}\label{OVsec}
In this section, we define the Ooguri--Vafa space \cite{OV} and give a description of its twistor coordinates. Most of what we say in this section can be found in \cite{GMN1} or \cite{GW}, so we will try to be concise in explaining what we need about the Ooguri--Vafa space.

\subsection{The Gibbons--Hawking ansatz}

The Ooguri--Vafa space can be constructed using the Gibbons--Hawking ansatz \cite{GHansatz}. We start with an open set $U\subset \mathbb{R}^{3}$ and a positive harmonic function $V\colon U\to \mathbb{R}$. We let $F= 2\pi {\rm i}\star {\rm d}V \in \Omega^{2}(U,{\rm i}\mathbb{R})$, where $\star$ denotes the Hodge star in $\mathbb{R}^3$ with the canonical Euclidean metric. We further assume that the cohomology class $\bigl[\frac{{\rm i}}{2\pi}F\bigr]\in \operatorname{Im}\bigl(H^2(U,\mathbb{Z})\to H^2(U,\mathbb{R})\bigr)$. Hence we can find a~principal ${\rm U}(1)$-bundle $\pi\colon X\to U$ with connection $\Theta\in \Omega^1(X,{\rm i}\mathbb{R})$, such that $\pi^{*}F={\rm d}\Theta$.

We now define for $j=1,2,3$, the following 2-forms on $X$:
\begin{equation}\label{sympGH}
 \omega_j=\left(\frac{{\rm i}}{2\pi}\Theta\right)\wedge \pi^*{\rm d}x^j +\pi^*\bigl(V \star {\rm d}x^{j}\bigr),
\end{equation}
where $x^j$ denotes the canonical coordinates of $\mathbb{R}^3$. The non-degeneracy is easy to check, and the fact that they are closed follows from the fact that $\frac{{\rm i}}{2\pi}{\rm d}\Theta=\frac{{\rm i}}{2\pi}\pi^*F=-\pi^*(\star {\rm d}V)$, so they define symplectic forms on $X$. The three symplectic forms $\omega_j$ correspond to the K\"ahler forms of a hyperk\"ahler structure $(X,g,I_1,I_2,I_3)$, where $g$ is the metric, $I_i$ are complex structures satisfying the imaginary quaternion relations, and $\omega_i(-,-)=g(I_i-,-)$. The metric can be written explicitly as
\begin{equation}\label{HKmetric}
 g=V^{-1}\left(\frac{{\rm i}}{2\pi}\Theta\right)\otimes \left(\frac{{\rm i}}{2\pi}\Theta\right)+V\pi^*\bigl({\rm d}x^1\otimes {\rm d}x^1+{\rm d}x^2\otimes {\rm d}x^2+{\rm d}x^3\otimes {\rm d}x^3\bigr).
\end{equation}
Hence, from the data of a positive harmonic function $V\colon U\subset \mathbb{R}^{3}\to \mathbb{R}$ such that $[\star {\rm d}V]\in \operatorname{Im}\bigl(H^2(U,\mathbb{Z})\to H^2(U,\mathbb{R})\bigr)$, the Gibbons--Hawking ansatz produces a hyperk\"ahler manifold $(X,g,I_1,I_2,I_3)$, where $\pi\colon X\to U$ is a principal ${\rm U}(1)$-bundle admitting a connection $\Theta$ with curvature ${\rm d}\Theta=\pi^*(2\pi {\rm i}\star {\rm d}V)$.

\subsection{Construction of the Ooguri--Vafa space}\label{constructionOV}

We now construct the Ooguri--Vafa space by using the Gibbons--Hawking ansatz. This is a~hyperk\"ahler space that we denote by $\mathcal{M}^{\text{ov}}(\Lambda)$, depending on a parameter $\Lambda \in \mathbb{C}^*$. Furthermore, for some open subset $\mathcal{B}\subset \mathbb{C}$ containing the origin, we will get a torus fibration $\mathcal{M}^{\text{ov}}(\Lambda)\to \mathcal{B}$ with a singular fiber at $0 \in \mathcal{B}$.

We start by taking the harmonic function on $\mathbb{R}^3\setminus \{0\}^{2}\times \mathbb{Z}$ defined by\footnote{$V$ can be though as the electro-magnetic potential of a $\mathbb{Z}$-worth of point unit charges, evenly distributed along the $x^3$-axis (see \cite{OV}).}
\begin{equation*}
 V\bigl(x^1,x^2,x^3\bigr):=\frac{1}{4\pi}\sum_{n=-\infty}^{\infty}\left(\frac{1}{\sqrt{\bigl(x^1\bigr)^2+\bigl(x^2\bigr)^2+\bigl(x^3+n\bigr)^2}} -c_n\right),
\end{equation*}
where $c_n \in \mathbb{R}_{\geq 0}$ are certain regularization constants making the sum converge \cite{OV}. After doing Poisson resummation, one obtains the following expression for $V$:
\begin{equation*}
 V\bigl(x^1,x^2,x^3\bigr):=-\frac{1}{2\pi}\operatorname{Log}\left(\frac{|z|}{|\Lambda|}\right)+\frac{1}{2\pi}\sum_{n\neq 0, n\in \mathbb{Z}} {\rm e}^{2\pi {\rm i}n x^3}K_0(2\pi|nz|),
\end{equation*}
where $z=x^1+{\rm i}x^2$, $|\Lambda| \in \mathbb{R}_{>0}$ is a constant related to the choice of the $c_n$, and $K_0$ is the $0$-th modified Bessel function of the second kind. Following \cite{GMN1}, we will denote the logarithm term of $V$ by $V^{\text{sf}}$ (the ``semiflat'' part), and the term with the series by $V^{\text{inst}}$ (the ``instanton'' part).

Now we would like to apply the Gibbons--Hawking ansatz to $V$. We have that $V$ is positive in an open subset of the form $U:=\mathcal{B}\times \mathbb{R}\setminus \{0\}^2\times \mathbb{Z}$, where $\mathcal{B}\subset \mathbb{R}^2$ is a neighborhood of the origin. Furthermore, it is also easy to check that the integrality condition for $[\star {\rm d}V]$ is satisfied. Hence, the Gibbons Hawking ansatz produces a principal ${\rm U}(1)$-bundle $\pi\colon X\to U$ with connection $\Theta$, such that $X$ carries a hyperk\"ahler metric of the form given by \eqref{HKmetric}. We remark that, because of the topology of $U$, all possible pairs $(\pi\colon X\to U,\Theta)$ with ${\rm d}\Theta=\pi^*(2\pi {\rm i} \star {\rm d}V)$ are gauge equivalent (see, for example, \cite[Theorem 2.5.1]{Kostant}), so they give isometric hyperk\"ahler spaces.

The pair $(\pi\colon X\to U, \Theta)$ actually admits an extra piece of structure, coming from the fact that $V$ (and hence $F=2\pi {\rm i} \star {\rm d}V$) is invariant under shifts $x^3\to x^3+n$ with $n\in \mathbb{Z}$. This extra piece of structure is a lift of the $\mathbb{Z}$-action to the total space $X$, preserving $\Theta$. There is a~${\rm U}(1)$-worth of ways of lifting the $\mathbb{Z}$-action, and we record this choice in the phase of $\Lambda \in \mathbb{C}^*$. By taking the quotient by the lift of the $\mathbb{Z}$-action we obtain a ${\rm U}(1)$ principal bundle $\widetilde{\pi}\colon\widetilde{X}\to \widetilde{U}$ with connection $\widetilde{\Theta}$, where $\widetilde{U}:=\bigl(\mathcal{B}\times S^1\bigr)\setminus \bigl(\{0\}^2\times \{1\}\bigr)$. Furthermore, the hyperk\"ahler structure of $X$ clearly descends to \smash{$\widetilde{X}$}.

Finally, one can show that by adding a point $\widetilde{X}$, the map $\widetilde{\pi}\colon\widetilde{X}\to \widetilde{U}$ smoothly extends to a~map $\widetilde{\pi}\colon \mathcal{M}^{\text{ov}}(\Lambda)\to \mathcal{B}\times S^1$ \cite[Proposition 3.2]{GW}, where $\mathcal{M}^{\text{ov}}(\Lambda)$ denotes $\widetilde{X}$ with the extra point. The hyperk\"ahler structure of $\widetilde{X}$ also smoothly extends to $\mathcal{M}^{\text{ov}}(\Lambda)$, which is what we call the Ooguri--Vafa space \cite[Section 4.1]{GMN1}. Furthermore, by composing $\widetilde{\pi}$ with the projection into the first factor, we can think of $\mathcal{M}^{\text{ov}}(\Lambda)$ as a singular torus fibration $p\colon\mathcal{M}^{\text{ov}}(\Lambda)\to \mathcal{B}$, with $p^{-1}(z)$ a torus for $z\in \mathcal{B}\cap \mathbb{C}^*$, and $p^{-1}(0)$ a torus with a node (see Figure~\ref{fig1} from the introduction).

We will now give a more explicit coordinate description of $\mathcal{M}^{\text{ov}}(\Lambda)$, following \cite[Section~4.1]{GMN1}. We start by writing an explicit solution to the equation ${\rm d}A=\star {\rm d}V$. Let $\mathcal{B}^*=\mathcal{B}\setminus\{0\}$. In the coordinates $\bigl(z=x^1+{\rm i}x^2, x^3\bigr)$ for $\mathcal{B}^*\times \mathbb{R}$, a solution to ${\rm d}A=\star {\rm d}V$ is given by
\begin{gather}
 A=A^{\text{sf}}+A^{\text{inst}},\qquad
 A^{\text{sf}}:=\frac{{\rm i}}{4\pi}\left(\operatorname{Log}\left(\frac{z}{\Lambda}\right)-\overline{\operatorname{Log}\left(\frac{z}{\Lambda}\right)}\right){\rm d}x^3,\nonumber\\
 A^{\text{inst}}:=-\frac{1}{4\pi}\Bigg(\sum_{n\neq 0} \operatorname{sgn}(n){\rm e}^{2\pi {\rm i}nx^3}|z|K_1(2\pi|nz|)\Bigg)\left(\frac{{\rm d}z}{z}-\frac{{\rm d}\overline{z}}{\overline{z}}\right),\label{localconnform}
\end{gather}
for some fixed $\Lambda \in \mathbb{C}^*$. Here $K_1$ denotes the first modified Bessel function of the second kind, $\operatorname{sgn}(n)=\frac{n}{|n|}$, and $\operatorname{Log}(z)$ uses the principal branch, so that $\operatorname{Im}(\operatorname{Log}(z))=\operatorname{Arg}(z)\in(-\pi,\pi)$. Consider the open cover of $\mathcal{B}^*\times \mathbb{R}$ given by the open sets
\begin{equation*}
D:=\bigl\{\bigl(z,x^3\bigr)\in \mathcal{B}^*\times \mathbb{R} \mid {\rm e}^{z/\Lambda}\not\in \mathbb{R}_{<0}\bigr\}, \qquad \widetilde{D}:=\bigl\{\bigl(z,x^3\bigr)\in \mathcal{B}^*\times \mathbb{R} \mid {\rm e}^{z/\Lambda}\not\in \mathbb{R}_{>0}\bigr\}.
\end{equation*}
We further have $D\cap \widetilde{D}=U_{+}\cup U_{-}$, where
\begin{equation*}
 U_{+}:=\bigl\{\bigl(z,x^3\bigr)\in \mathcal{B}^*\times \mathbb{R} \mid \operatorname{Im}\bigl({\rm e}^{z/\Lambda}\bigr)>0\bigr\}, \qquad U_{-}:=\bigl\{\bigl(z,x^3\bigr)\in \mathcal{B}^*\times \mathbb{R} \mid \operatorname{Im}\bigl({\rm e}^{z/\Lambda}\bigr)<0\bigr\}.
\end{equation*}
Over $D$ consider the trivial principal ${\rm U}(1)$-bundle with connection $\Theta:={\rm i}{\rm d}\theta_m +2\pi {\rm i}A$, where $\theta_m$ is an angle coordinate for the ${\rm U}(1)$-fiber (called the magnetic angle in \cite{GMN1}). Similarly, consider on $\widetilde{D}$ the trivial ${\rm U}(1)$-bundle with connection $\widetilde{\Theta}:={\rm i}{\rm d}\widetilde{\theta}_m+ 2\pi {\rm i}\widetilde{A}$, where $\widetilde{A}$ is given by the same~formula from \eqref{localconnform} but with the branch of $\operatorname{Log}(z)$ such that $\operatorname{Im}(\operatorname{Log}(z))\in (0,2\pi)$. By the Gibbons--Hawking ansatz, $(D\times {\rm U}(1),\Theta)$ and $\bigl(\widetilde{D}\times {\rm U}(1), \widetilde{\Theta}\bigr)$ define HK metrics on their total space, with metric and K\"ahler forms given by \eqref{HKmetric} and \eqref{sympGH}. On the other hand, note that
\begin{equation}\label{translocalconn}
 \widetilde{A}=A \quad \text{on}\quad U_{+}, \qquad \widetilde{A}=A -{\rm d}x^3 \quad \text{on}\quad U_{-}.
\end{equation}
To obtain $\mathcal{M}^{\text{ov}}(\Lambda)|_{\mathcal{B}^*}$, one identifies $D\times {\rm U}(1)$ \big(with coordinates $\bigl(z,x^3,{\rm e}^{{\rm i}\theta_m}\bigr)$\big) and $\widetilde{D}\times {\rm U}(1)$ \big(with coordinates \smash{$\bigl(z,x^3,{\rm e}^{{\rm i}\widetilde{\theta}_m}\bigr)$}\big) by $\phi\colon (U_{+}\cup U_{-})\times {\rm U}(1)\to (U_{+}\cup U_{-})\times {\rm U}(1)$ given by
\begin{equation*}
 \phi\bigl(z,x^3,{\rm e}^{{\rm i}\theta_m}\bigr)=\begin{cases}
 \bigl(z,x^3,{\rm e}^{{\rm i}\theta_m}\bigr) & \text{if} \ \bigl(z,x^3\bigr) \in U_{+},\\
 \bigl(z,x^3,{\rm e}^{{\rm i}\theta_m+2\pi {\rm i}x^3 -\mathrm{i\pi}}\bigr) & \text{if} \ \bigl(z,x^3\bigr)\in U_{-},
 \end{cases}
\end{equation*}
and then one performs a quotient by translations $x^3\to x^3+n$, $n\in \mathbb{Z}$. In particular, note that on $U_{-}$ we have
\begin{equation}\label{magangletrans}
 \widetilde{\theta}_m=\theta_m+2\pi x^3 -\pi \quad \text{(mod $2\pi$)}.
\end{equation}
Both \eqref{translocalconn} and \eqref{magangletrans} imply that the corresponding expressions \eqref{HKmetric} and \eqref{sympGH} on $D\times {\rm U}(1)$ and~${\widetilde{D}\times {\rm U}(1)}$ give a well defined global expression on $\mathcal{M}^{\text{ov}}(\Lambda)|_{\mathcal{B}^*}$, since
\begin{equation*}
\frac{{\rm d}\widetilde{\theta}_m}{2\pi}+\widetilde{A}=\frac{{\rm d}\theta_e}{2\pi}+A \qquad \text{on}\quad U_{\pm}\times {\rm U}(1),
\end{equation*}
and both $A$ and $\widetilde{A}$ are periodic in $x^3$. The shift by $-\pi$ in \eqref{magangletrans} is justified in \cite[Section 4.1]{GMN1} as the choice that allows the above HK structure to extend over $z=0$ and match the Ooguri--Vafa hyperk\"ahler structure constructed at the beginning of the section.

\subsection{Twistor coordinates}

We now recall the description of the hyperk\"ahler structure of $\mathcal{M}^{\text{ov}}(\Lambda)$ using the twistor coordinates from \cite{GMN1}. These coordinates are used to describe the twistor family of holomorphic symplectic forms $\Omega^{\text{ov}}(\xi)$, $\xi\in \CP$, encoding the hyperk\"ahler structure of $\mathcal{M}^{\text{ov}}(\Lambda)$. We use the conventions of \cite{GMN1}, so in particular we have the following formulas for the $\CP$-worth of complex structures and holomorphic symplectic forms associated to $\mathcal{M}^{\text{ov}}(\Lambda)$
\begin{gather*}
 I(\xi)=\frac{{\rm i}\bigl(-\xi +\overline{\xi}\bigr)}{1+|\xi|^2}I_1 - \frac{\xi + \overline{\xi}}{1+|\xi|^2}I_2 + \frac{1-|\xi|^2}{1+|\xi|^2}I_3 \qquad \text{for} \quad \xi \in \mathbb{C}\subset \CP,\\
 \Omega^{\text{ov}}(\xi)=-\frac{{\rm i}}{2}\xi^{-1}(\omega_1+{\rm i}\omega_2)+\omega_3-\frac{{\rm i}}{2}\xi(\omega_1-{\rm i}\omega_2) \qquad\text{for} \quad \xi \in \mathbb{C}^*.
\end{gather*}
To obtain the holomorphic symplectic forms corresponding to $\xi=0$ and $\xi=\infty$ from the above formula, we consider $\xi\Omega^{\text{ov}}(\xi)|_{\xi=0}$ for $\xi=0$, and $\xi^{-1}\Omega^{\text{ov}}(\xi)|_{\xi=\infty}$ for $\xi=\infty$.

Using \eqref{sympGH} and \eqref{localconnform}, we find that the K\"ahler forms $\omega_i$ of $\mathcal{M}^{\text{ov}}(\Lambda)$ are given by
\begin{equation}\label{sympOV}
 \omega_i = {\rm d}x^i\wedge \left(\frac{{\rm d}\theta_m}{2\pi} + A\right) + V\star {\rm d}x^i.
\end{equation}
Plugging \eqref{sympOV} into the expression of $\Omega^{\text{ov}}(\xi)$, we have that, away from $z=0\in \mathcal{B}$,
 $
 \Omega^{\text{ov}}(\xi)=\frac{1}{4\pi^2}\xi_m\wedge\xi_e$,
 where
 \begin{equation*}
 \xi_e=2\pi {\rm i}{\rm d}x^3 + \frac{\pi}{\xi}{\rm d}z+\pi\xi {\rm d}\overline{z},\qquad
 \xi_m=\pi {\rm i}V\left(\frac{1}{\xi}{\rm d}z-\xi {\rm d}\overline{z}\right)+ {\rm i}{\rm d}\theta_m +2\pi {\rm i}A.
 \end{equation*}

From the fact that $(\mathcal{M}^{\text{ov}}(\Lambda),I(\xi),\Omega^{\text{ov}}(\xi))$ is a holomorphic symplectic manifold, we conclude that $\xi_e(\xi)$ and $\xi_m(\xi)$ must be of type $(1,0)$ in holomorphic structure $I(\xi)$.

Now we define the twistor coordinates from \cite{GMN1}. Denoting $\theta_e:=2\pi x^3$ (the ``electric angle''), we have that $\xi_e={\rm d}\mathcal{X}_e^\text{ov} / \mathcal{X}_e^\text{ov}$, where
\begin{equation}\label{OV electric}
 \mathcal{X}_e^\text{ov}(\xi):=\exp \left(\frac{\pi}{\xi}z+{\rm i}\theta_e +\pi\xi \overline{z}\right).
\end{equation}
We will call $\mathcal{X}_e^\text{ov}$ the electric twistor coordinate. Since $\xi_e(\xi)$ is a $(1,0)$ form is complex structure~$I(\xi)$, we conclude that $\mathcal{X}_e^\text{ov}(\xi)$ defines a holomorphic function on $\mathcal{M}^{\text{ov}}(\Lambda)$ in complex structure~$I(\xi)$.

Now we define the magnetic twistor coordinate $\mathcal{X}_m^\text{ov}$. This coordinate satisfies
\begin{equation}\label{familyholsym}
 \Omega^{\text{ov}}(\xi)=-\frac{1}{4\pi^2}\frac{{\rm d}\mathcal{X}_e^\text{ov}(\xi)}{\mathcal{X}_e^\text{ov}(\xi)}\wedge \frac{{\rm d}\mathcal{X}_m^\text{ov}(\xi)}{\mathcal{X}_m^\text{ov}(\xi)},
\end{equation}
so $\mathcal{X}_m^\text{ov}(\xi)$ also gives a holomorphic function on $\mathcal{M}^{\text{ov}}(\Lambda)$ in holomorphic structure $I(\xi)$. To define this coordinate we first write $\mathcal{X}_m^\text{ov}=\mathcal{X}_m^{\text{sf}}\mathcal{X}_m^{\text{inst}}$. The first factor $\mathcal{X}_m^{\text{sf}}$ is defined by
\begin{equation}\label{Xmsf}
 \mathcal{X}_m^{\text{sf}}(\xi):=\exp \left(\frac{1}{\xi}\frac{(z\operatorname{Log}(z/\Lambda)-z)}{2 {\rm i}}+{\rm i}\theta_m -\xi\frac{\bigl(\overline{z} \operatorname{Log}\bigl(\overline{z}/\overline{\Lambda}\bigr)-\overline{z}\bigr)}{2 {\rm i}}\right).
\end{equation}
On the other hand, $\mathcal{X}_m^{\text{inst}}$ is given by
\begin{align}
 \mathcal{X}_m^{\text{inst}}(\xi)={}&\exp \biggl(\frac{{\rm i}}{4\pi}\int_{l_{+}(z)}\frac{{\rm d}\xi'}{\xi'}\frac{\xi+\xi'}{\xi'-\xi}\operatorname{Log}(1-\mathcal{X}_e^\text{ov}(\xi'))\nonumber\\
 &- \frac{{\rm i}}{4\pi}\int_{l_{-}(z)}\frac{{\rm d}\xi'}{\xi'}\frac{\xi+\xi'}{\xi'-\xi}\operatorname{Log}\bigl(1-(\mathcal{X}_e^\text{ov}(\xi'))^{-1}\bigr)\biggr),\label{Xminst}
\end{align}
where the integration contours $l_{\pm}(z)$ are the rays
\begin{equation}\label{rays}
 l_{\pm}(z)=\biggl\{\xi \in \mathbb{C}^* \mid \pm\frac{z}{\xi}<0\biggr\}
\end{equation}
oriented from $0$ to $\infty$.

In \cite{GMN1}, it is verified that $\mathcal{X}_e^\text{ov}(\xi)$ and $\mathcal{X}_m^\text{ov}(\xi)$ indeed satisfy \eqref{familyholsym}, so that $\operatorname{Log}(\mathcal{X}_e^\text{ov}(\xi))$ and~$\operatorname{Log}(\mathcal{X}_m^\text{ov}(\xi))$ give twistor holomorphic Darboux coordinates.

\subsection{Some properties of the magnetic twistor coordinate}

Here we state some of the properties satisfied by the magnetic twistor coordinate $\MC^\text{ov}(\xi)$. These will serve as guiding principles to construct an analog in the wild harmonic bundle setting. As before, the main reference and proof of the statements can be found in \cite{GMN1}.

\begin{Proposition}[jumps of the twistor coordinate]\label{jumpsOV} For $z\in \mathcal{B}\setminus\{0\}$, consider the rays~$l_{\pm}(z)$ defined in \eqref{rays}. We then have that $\mathcal{X}_m^\textnormal{ov}(\xi)$ is holomorphic in $\xi$ away from $l_{\pm}(z)$, and furthermore it has the following jumps along the rays
\begin{gather*}
 \mathcal{X}_m^\textnormal{ov}(\xi)^{+}=\mathcal{X}_m^\textnormal{ov}(\xi)^{-}(1-\mathcal{X}_e^\textnormal{ov}(\xi))^{-1} \qquad \text{along} \quad \xi \in l_{+}(z),\\
 \mathcal{X}_m^\textnormal{ov}(\xi)^{+}=\mathcal{X}_m^\textnormal{ov}(\xi)^{-}\bigl(1-\mathcal{X}_e^\textnormal{ov}(\xi)^{-1}\bigr)\qquad \text{along} \quad \xi \in l_{-}(z),
\end{gather*}
where the $\pm$ denotes the fact that we approach $l_{\pm}(z)$ clockwise or anticlockwise, respectively.
\end{Proposition}

\begin{Proposition}[asymptotics of the twistor coordinate]\label{asympov} $\MC^\textnormal{ov}(\xi)$ has the following asymptotics:
\begin{gather}\label{aft}
 \MC^\textnormal{ov}(\xi) \sim
 \begin{cases}
\displaystyle \exp \left(-\frac{{\rm i}}{2\xi}(z \operatorname{Log}(z/\Lambda)-z) +{\rm i}\theta_m +\frac{1}{2\pi {\rm i}}\sum_{s\neq 0}\frac{1}{s}{\rm e}^{{\rm i}s\theta_e}K_0(2\pi |sz|) \right) \\ \qquad \text{as}\quad \xi \to 0,\\
\displaystyle \exp \left(\frac{{\rm i}\xi}{2}\bigl(\overline{z} \operatorname{Log}\bigl(\overline{z}/\overline{\Lambda}\bigr)-\overline{z}\bigr) +{\rm i}\theta_m -\frac{1}{2\pi {\rm i}}\sum_{s\neq 0}\frac{1}{s}{\rm e}^{{\rm i}s\theta_e}K_0(2\pi |sz|) \right)\\ \qquad\text{as} \quad \xi \to \infty.
 \end{cases}
\end{gather}
\end{Proposition}

\begin{Proposition}[reality condition] $\MC^\textnormal{ov}(\xi)$ satisfies the following reality condition:
\begin{equation}\label{rc}
 \MC^\textnormal{ov}(\xi)=\overline{\MC^{\textnormal{ov}}\bigl(-1/\overline{\xi}\bigr)}^{-1}.
\end{equation}
\end{Proposition}

\section{Framed wild harmonic bundles}\label{sec3}

We now go to the subject of wild harmonic bundles. This section is roughly divided into two big parts:
\begin{itemize}\itemsep=0pt
 \item In the first part, consisting of Sections \ref{filtbunsec}--\ref{secdeffwhb}, we start by recalling the notions of filtered Higgs bundles, filtered flat bundles, and wild harmonic bundles. We then define what we mean by ``framed wild harmonic bundles'', and specify which type of framed wild harmonic bundles we will consider for our moduli space.
 \item In the second part, consisting of Sections \ref{sec3.4}--\ref{nonvanishing}, we start by recalling some facts about Stokes data, and then define the analogs of $\mathcal{X}_e^{\text{ov}}(\xi)$ and $\mathcal{X}_m^{\text{ov}}(\xi)$ for our set of framed wild harmonic bundles, which we denote by $\mathcal{X}_e(\xi)$ and $\mathcal{X}_m(\xi)$. The rest of the section is devoted to showing holomorphic dependence of $\mathcal{X}_e(\xi)$ and $\mathcal{X}_m(\xi)$ with respect to $\xi\in \mathbb{C}^{*}$, and computing the asymptotics of $\mathcal{X}_e(\xi)$ and $\mathcal{X}_m(\xi)$ as $\xi\to 0$ and $\xi\to \infty$. These properties will be crucial for the identification with the Ooguri--Vafa space in Section \ref{sec4}.
\end{itemize}

\subsection{Filtered bundles on curves}\label{filtbunsec}

We start by recalling the notions of filtered \cite{Sim2,Sim} and parabolic bundles \cite{ParBun,ParBun2} of finite rank, and parabolic degree. Most of what we say in this section and the notations that we use are based on \cite{M,M2,M5}.

Let $X$ be a Riemann surface and $D\subset X$ a discrete set of points. We will denote by $\mathbb{R}^D$ the set of maps $D\to \mathbb{R}$, and its elements by $\boldsymbol{a}$.
\begin{Definition} \label{deffilt} A filtered bundle $\mathcal{P}_*\mathcal{E}:=\bigl(\mathcal{E},\bigl\{\mathcal{P}_{\boldsymbol{a}}\mathcal{E} \mid \boldsymbol{a}\in \mathbb{R}^{D} \bigr\}\bigr)$ over $(X,D)$ consists of the following data:
\begin{itemize}\itemsep=0pt
 \item $\mathcal{E}\to (X,D)$ is a locally free $\mathcal{O}_{X}(*D)$-module over $X$ with finite rank, i.e., $\mathcal{E}$ is a meromorphic bundle over $X$ with poles along $D$.
 \item Each $\mathcal{P}_{\boldsymbol{a}}\mathcal{E}$ is a locally free $\mathcal{O}_{X}$-submodule of $\mathcal{E}$, i.e., a holomorphic subbundle of $\mathcal{E}$ over $X$.
\end{itemize}
This data must satisfy the following conditions:
\begin{itemize}\itemsep=0pt
 \item $\mathcal{P}_{\boldsymbol{a}}\mathcal{E}\otimes_{\mathcal{O}_X}\mathcal{O}_{X}(*D)=\mathcal{E}$ for $\boldsymbol{a}\in \mathbb{R}^D$. In particular, $\mathcal{P}_{\boldsymbol{a}}\mathcal{E}|_{X\setminus D}=\mathcal{E}|_{X\setminus D}$.
 \item If $p\in D$, the stalk $\mathcal{P}_{\boldsymbol{a}}\mathcal{E}|_p$ depends only on $a_p:=\boldsymbol{a}(p)\in \mathbb{R}$. Hence, we will sometimes write~${\mathcal{P}_{a_p}\mathcal{E}|_p:=\mathcal{P}_{\boldsymbol{a}}\mathcal{E}|_p}$.
 \item For $p \in D$, we have $\mathcal{P}_{\boldsymbol{a}}\mathcal{E}|_p\subset \mathcal{P}_{\boldsymbol{b}}\mathcal{E}|_p$ if $a_p\leq b_p$. Furthermore, for any $\boldsymbol{a} \in \mathbb{R}^D$, there is $\epsilon>0$ such that $\mathcal{P}_{a_p}\mathcal{E}|_p=\mathcal{P}_{a_p+\epsilon}\mathcal{E}|_p$.
 \item For $n \in \mathbb{Z}$ and $p \in D$, we have that $\mathcal{P}_{a_p+n}\mathcal{E}|_p=\mathcal{P}_{a_p}\mathcal{E}|_p\otimes_{\mathcal{O}_{X,p}} \mathcal{O}_{X,p}(np)$.
\end{itemize}
\end{Definition}

Let $\boldsymbol{c}\in \mathbb{R}^{D}$, and consider the filtration $\mathcal{F}$ of $\mathcal{P}_{\boldsymbol{c}}\mathcal{E}$ indexed by $\bigl\{\boldsymbol{d}\in \mathbb{R}^{D} \mid \boldsymbol{d}(p)\in (c_p-1,c_p] \bigr\}$ and defined by $\mathcal{F}_{\boldsymbol{d}}(\mathcal{P}_{\boldsymbol{c}}\mathcal{E}):=\bigcup_{\boldsymbol{a}\leq \boldsymbol{d}}\mathcal{P}_{\boldsymbol{a}}\mathcal{E}$, where $\boldsymbol{a}\leq \boldsymbol{d}$ if and only if $\boldsymbol{a}(p)\leq \boldsymbol{d}(p)$ for all $p\in D$. Then the filtration is parabolic in the sense that for each $p \in D$, the set $\bigl\{d\in (c_p-1,c_p] \mid \text{Gr}^{\mathcal{F}}_d(\mathcal{P}_{\boldsymbol{c}}\mathcal{E}|_p)\neq 0\bigr\}$ is finite, where \smash{$\text{Gr}^{\mathcal{F}}_{d}(\mathcal{P}_{\boldsymbol{c}}\mathcal{E}|_p)
:=\mathcal{F}_d(\mathcal{P}_{\boldsymbol{c}}\mathcal{E}|_p)/\mathcal{F}_{<d}(\mathcal{P}_{\boldsymbol{c}}\mathcal{E}|_p)$}. The data of this filtration is called the $\boldsymbol{c}$-truncation of the filtered bundle $\mathcal{P}_{*}\mathcal{E}\to (X,D)$.

\begin{Definition} Given $\boldsymbol{c}\in \mathbb{R}^{D}$, a $\boldsymbol{c}$-parabolic bundle $_{\boldsymbol{c}}\mathcal{E}$ over $(X,D)$ consists of the following data:
\begin{itemize}\itemsep=0pt
 \item $_{\boldsymbol{c}}\mathcal{E} \to X$ is a holomorphic bundle.
 \item For each $p\in D$, we have an increasing filtration $\mathcal{F}_d(_{\boldsymbol{c}}\mathcal{E}|_p)$ of the fiber $_{\boldsymbol{c}}\mathcal{E}|_p$ indexed by~${d\in (c_p-1,c_p]}$.
\end{itemize}
This data must satisfy the following conditions:
\begin{itemize}\itemsep=0pt
 \item For $p\in D$, $\mathcal{F}_a(_{\boldsymbol{c}}\mathcal{E}|_p)=\bigcap_{a<d}\mathcal{F}_d(_{\boldsymbol{c}}\mathcal{E}|_p)$.
 \item For $p\in D$, $_{\boldsymbol{c}}\mathcal{E}|_p=\bigcup_{d}\mathcal{F}_d(_{\boldsymbol{c}}\mathcal{E}|_p)$.
\end{itemize}
The set of $d\in (c_p-1,c_p]$ such that $\text{Gr}^{\mathcal{F}}_{d}(_{\boldsymbol{c}}\mathcal{E}|_p):=\mathcal{F}_d(_{\boldsymbol{c}}\mathcal{E}|_p)/\mathcal{F}_{<d}(_{\boldsymbol{c}}\mathcal{E}|_p)\neq 0$ is called the parabolic weights at $p$ of the $\boldsymbol{c}$-parabolic bundle.
\end{Definition}
It is easy to see that given a filtered bundle $\mathcal{P}_*\mathcal{E}\to (X,D)$ and $\boldsymbol{c}\in \mathbb{R}^{D}$, its $\boldsymbol{c}$-truncation gives rise to a $\boldsymbol{c}$-parabolic bundle $(_{\boldsymbol{c}}\mathcal{E}\to (X,D),\mathcal{F})$. Conversely, given a $\boldsymbol{c}$-parabolic bundle~$(_{\boldsymbol{c}}\mathcal{E},\mathcal{F})$, one can obtain a filtered bundle by taking $\mathcal{E}:= {}_{\boldsymbol{c}}\mathcal{E}\otimes_{\mathcal{O}_X} \mathcal{O}_X(*D)$, and $\bigl\{\mathcal{P}_{\boldsymbol{a}}\mathcal{E} \mid \boldsymbol{a}\in \mathbb{R}^{D} \bigr\}$ induced from $\mathcal{F}$. Since one can reconstruct a filtered bundle from any of its $\boldsymbol{c}$-truncations, the data of~$\mathcal{P}_*\mathcal{E}$ is equivalent to the data of $_{\boldsymbol{c}}\mathcal{E}$.

With this terminology, if $\boldsymbol{c}=\boldsymbol{0}$, where $\boldsymbol{0}(p)=0$ for all $p\in D$, we have that a $\boldsymbol{0}$-parabolic bundle is what is usually called a parabolic bundle.

\begin{Definition} Let $X$ be a compact Riemann surface and $D\subset X$ a finite subset of points. Given a $\boldsymbol{c}$-parabolic bundle $_{\boldsymbol{c}}E\to (X,D)$, we define its parabolic degree $\operatorname{pdeg}(_{\boldsymbol{c}}\mathcal{E})$ as follows:
\begin{equation*}
 \textnormal{pdeg}(_{\boldsymbol{c}}\mathcal{E}):= \textnormal{deg}(_{\boldsymbol{c}}\mathcal{E}) - \sum_{p\in D}\sum_{d\in (c_p-1,c_p]} d\cdot \textnormal{dim}_{\mathbb{C}}\left(\frac{\mathcal{F}_d(_{\boldsymbol{c}}\mathcal{E}|_p)}{{\mathcal{F}_{<d}(_{\boldsymbol{c}}\mathcal{E}|_p)}}\right),
\end{equation*}
where $\textnormal{deg}(-)$ is the usual degree of a vector bundle.

If $\mathcal{P}_*\mathcal{E} \to (X,D)$ is a filtered bundle, we define its parabolic degree by
$
 \textnormal{pdeg}(\mathcal{P}_*\mathcal{E}):=\textnormal{pdeg}(_{\boldsymbol{c}}\mathcal{E})$,
where $_{\boldsymbol{c}}\mathcal{E}$ is any $\boldsymbol{c}$-truncation of $\mathcal{P}_*\mathcal{E}$. It is easy to check that this is well defined.

\end{Definition}

\subsection{Wild harmonic bundles and their associated filtered objects}\label{WHB}

Let $X$ be a compact Riemann surface and $D\subset X$ a finite subset.
\begin{Definition} A harmonic bundle over $X\setminus D$ is a tuple $\bigl(E,\overline{\partial}_E,\theta,h\bigr)$ such that
\begin{itemize}\itemsep=0pt
 \item $\bigl(E,\overline{\partial}_E\bigr)\to (X\setminus D)$ is a holomorphic bundle with hermitian metric $h$.
 \item \smash{$\theta\in \Omega^{(1,0)}_{X\setminus D}(\textnormal{End}(E))$} and $\overline{\partial}_E(\theta)=0$. This endomorphism valued 1-form is known as the Higgs field.
 \item The connection $\nabla = D\bigl(\overline{\partial}_E,h\bigr)+\theta + \theta^{\dagger_h}$ is flat,\footnote{Given $E\to (X\setminus D)$ as above and a flat smooth connection $\nabla$ on $E$, we will frequently equip $E$ with the holomorphic structure given by the $(0,1)$ part of $\nabla$. The connection $\nabla$ then induces a holomorphic connection on the holomorphic bundle $\bigl(E,\nabla^{(0,1)}\bigr)$, which we will also denote by $\nabla$. We hope that this abuse of notation does not cause confusion.} where $D\bigl(\overline{\partial}_E,h\bigr)$ denotes the Chern connection, and $\theta^{\dagger_h}$ is the adjoint of $\theta$ with respect to $h$. Equivalently, the Hitchin equation is satisfied
\smash{$ F\bigl(D\bigl(\overline{\partial}_E,h\bigr)\bigr)+\bigl[\theta,\theta^{\dagger_h}\bigr]=0$},
 where $F(D\bigl(\overline{\partial}_E,h\bigr))$ denotes the curvature of the Chern connection.
\end{itemize}

We say that the harmonic bundle $\bigl(E,\overline{\partial}_E,\theta,h\bigr) \to X\setminus D$ is unramified and wild over $(X,D)$,
 if for every $p\in D$ there is a holomorphic coordinate neighborhood $(U_p,z)$ with $z(p)=0$, and a~finite set of (non-zero) irregular values $\textnormal{Irr}(\theta)_p\subset z^{-1}\mathbb{C}\bigl[z^{-1}\bigr]$ such that
 \begin{equation*}
 \bigl(E,\overline{\partial}_E,\theta\bigr)|_{U_p\setminus \{p\}}=\bigoplus_{\mathfrak{a}\in \textnormal{Irr}(\theta)_p}\bigl(E_{\mathfrak{a}},\overline{\partial}_{E_{\mathfrak{a}}},\theta_{\mathfrak{a}}\bigr),
\end{equation*}
and $\theta_{\mathfrak{a}}-{\rm d}\mathfrak{a}\cdot \textnormal{Id}_{E_{\mathfrak{a}}}$ has at most a simple pole at $p$, where $\textnormal{Id}_{E_{\mathfrak{a}}}$ denotes the identity map of $E_{\mathfrak{a}}$.
\end{Definition}

\begin{Definition}\label{filteredhiggsdef} An unramified filtered Higgs bundle over $(X,D)$ is a pair $(\mathcal{P}_*\mathcal{E},\theta)$, where
\begin{itemize}\itemsep=0pt
 \item $\mathcal{P}_*\mathcal{E}=\bigl(\mathcal{E},\bigl\{\mathcal{P}_{\boldsymbol{a}}\mathcal{E} \mid \boldsymbol{a} \in \mathbb{R}^{D}\bigr\}\bigr)$ is a filtered bundle over $(X,D)$.
 \item \smash{$\theta \in \Omega^{(1,0)}_X(\textnormal{End}(\mathcal{E}))$} and $\overline{\partial}_{\mathcal{E}}\theta=0$.
\end{itemize}
The pair $(\mathcal{P}_*\mathcal{E},\theta)$ must satisfy the following condition:
\begin{itemize}\itemsep=0pt
 \item For every $p\in D$ and every $\boldsymbol{a}\in \mathbb{R}^{D}$, there is a holomorphic coordinate neighborhood $(U_p,z)$ with $z(p)=0$ and a finite set of irregular values $\textnormal{Irr}(\theta)_p\subset z^{-1}\mathbb{C}\bigl[z^{-1}\bigr]$ such that
 \begin{equation}\label{irregdecomp}
 (\mathcal{P}_{\boldsymbol{a}}\mathcal{E},\theta)|_{U_p}=\bigoplus_{\mathfrak{a}\in \textnormal{Irr}(\theta)_p}(\mathcal{P}_{\boldsymbol{a}}\mathcal{E}_{\mathfrak{a}},\theta_{\mathfrak{a}}),
\end{equation}
where $\theta_{\mathfrak{a}}-{\rm d}\mathfrak{a}\cdot\textnormal{Id}_{\mathcal{P}_{\boldsymbol{a}}\mathcal{E}_{\mathfrak{a}}}$ has at most a simple pole at $p$ as a meromorphic endomorphism of~$\mathcal{P}_{\boldsymbol{a}}\mathcal{E}_{\mathfrak{a}}$.
\item Given any $p\in D$ and $\boldsymbol{a}\in \mathbb{R}^D$, we assume that the irregular decomposition \eqref{irregdecomp} is compatible with the parabolic filtration on $\mathcal{P}_{\boldsymbol{a}}\mathcal{E}|_p$.
\end{itemize}

By taking the $\boldsymbol{c}$-truncation of the filtered bundle, we also have the corresponding notion of unramified $\boldsymbol{c}$-parabolic Higgs bundle.
\end{Definition}

The notion of unramified filtered flat bundle is similar to the notion of unramified filtered Higgs bundle, but it differs in the condition that we put on the splitting near the points $p\in D$.

\begin{Definition} An unramified filtered flat bundle over $(X,D)$ is a pair $(\mathcal{P}_*\mathcal{E},\nabla)$, where
\begin{itemize}\itemsep=0pt
 \item $\mathcal{P}_*\mathcal{E}=\bigl(\mathcal{E},\bigl\{\mathcal{P}_{\boldsymbol{a}}\mathcal{E} \mid \boldsymbol{a} \in \mathbb{R}^{D}\bigr\}\bigr)$ is a filtered bundle over $(X,D)$.
 \item $\nabla\colon\mathcal{E}\to \Omega^1_X\otimes \mathcal{E}$ is a flat meromorphic connection.
\end{itemize}
The pair $(\mathcal{P}_*\mathcal{E},\nabla)$ satisfies the following condition:
\begin{itemize}\itemsep=0pt
 \item For every $p\in D$ and every $\boldsymbol{a}\in \mathbb{R}^{D}$, there is a holomorphic coordinate neighborhood $(U_p,z)$ with $z(p)=0$ and a finite set of irregular values $\textnormal{Irr}(\nabla)_p\subset z^{-1}\mathbb{C}\bigl[z^{-1}\bigr]$ such that
 \begin{equation}\label{irrdec}
 (\mathcal{P}_{\boldsymbol{a}}\mathcal{E},\nabla)|_{U_p}\otimes_{\mathcal{O}(U_p)}\mathbb{C}[[z]]=\bigoplus_{\mathfrak{a}\in \textnormal{Irr}(\nabla)_p}\bigl(\widehat{\mathcal{P}_{\boldsymbol{a}}\mathcal{E}}_{\mathfrak{a}},\widehat{\nabla}_{\mathfrak{a}}\bigr),
\end{equation}
where \smash{$\widehat{\mathcal{P}_{\boldsymbol{a}}\mathcal{E}}_{\mathfrak{a}}$} are free $\mathbb{C}[[z]]$-modules with formal meromorphic connection $\smash{\widehat{\nabla}_{\mathfrak{a}}\colon \widehat{\mathcal{P}_{\boldsymbol{a}}\mathcal{E}}_{\mathfrak{a}}}\to \smash{\widehat{\mathcal{P}_{\boldsymbol{a}}\mathcal{E}}_{\mathfrak{a}} \otimes _{\mathcal{O}(U_p)} \Omega^1_{U_p}(* p)}$, and \smash{$\widehat{\nabla}_{\mathfrak{a}}-d\mathfrak{a}\cdot\textnormal{Id}_{\widehat{\mathcal{P}_{\boldsymbol{a}}\mathcal{E}_{\mathfrak{a}}}}$} has at most a simple pole at $p$.\footnote{If in a holomorphic trivialization near $p$ we have that $\nabla={\rm d} +A_k\frac{{\rm d}z}{z^k} + \text{lower order terms}$, for some $k>1$ and with \smash{$A_k\in \operatorname{End}\bigl(\mathbb{C}^{\text{rank}(\mathcal{P}_{\boldsymbol{a}}\mathcal{E})}\bigr)$}; then a sufficient condition for \eqref{irrdec} to hold is to have $A_k$ diagonalizable with distinct eigenvalues (see, for example, \cite[Lemma 1]{B2}).} In other words, the holomorphic bundle with flat meromorphic connection $(\mathcal{P}_{\boldsymbol{a}}\mathcal{E},\nabla)$ can be ``block diagonalized'' by a formal gauge transformation near $p$.
\item Given any $p\in D$ and $\boldsymbol{a}\in \mathbb{R}^D$, we assume that the irregular decomposition~\eqref{irrdec} is compatible with the filtration of the parabolic structure on $\mathcal{P}_{\boldsymbol{a}}\mathcal{E}|_p$.
\end{itemize}

By taking the $\boldsymbol{c}$-truncation of the filtered bundle, we also have the corresponding notion of unramified $\boldsymbol{c}$-parabolic flat bundle.
\end{Definition}
\begin{Remark} Since all the objects we will consider are unramified, we will drop the adjective from now on.
\end{Remark}

We now explain how to obtain a filtered Higgs bundle and a filtered flat bundle from a wild harmonic bundle. The main idea is that the filtered structure comes from the growth of certain sections with respect to the harmonic metric.

\begin{Definition}\label{phdef} Let $\bigl(E,\overline{\partial}_E,\theta,h\bigr)$ be a wild harmonic bundle over $X\setminus D$, and consider the holomorphic bundle \smash{$\mathcal{E}^{\xi}:=\bigl(E,\overline{\partial}_E+\xi\theta^{\dagger_h}\bigr)\to X\setminus D$}, where $\xi\in \mathbb{C}$. Given $\boldsymbol{a}\in \mathbb{R}^D$, we define the holomorphic bundle $\mathcal{P}_{\boldsymbol{a}}^h\mathcal{E}^{\xi}\to X$ as follows:
\begin{itemize}\itemsep=0pt
 \item If $U\subset X$ is an open subset with $U\cap D=\varnothing$, we let the holomorphic sections over $U$ be~$\mathcal{P}_{\boldsymbol{a}}^h\mathcal{E}^{\xi}(U)=\mathcal{E}^{\xi}(U)$.
 \item If $U\cap D\neq \varnothing$ and $s\in \mathcal{E}^{\xi}(U-D)$, we say that $s\in \mathcal{P}^h_{\boldsymbol{a}}\mathcal{E}^{\xi}(U)$ if for every $p\in D\cap U$ we have~${|s|_h=\mathcal{O}(|z|^{-a_p-\epsilon})}$ for every $\epsilon>0$, where $z$ is a holomorphic coordinate vanishing at~$p$.
\end{itemize}
We denote $\mathcal{P}^h_*\mathcal{E}^{\xi}:=\bigl(\mathcal{P}_{\boldsymbol{0}}^h\mathcal{E}^{\xi}\otimes_{\mathcal{O}_X}\mathcal{O}_X(*D), \bigl\{\mathcal{P}_{\boldsymbol{a}}^h\mathcal{E}^{\xi} \mid \boldsymbol{a}\in \mathbb{R}^D\bigr\}\bigr)$.
\end{Definition}
\begin{Theorem}[{\cite[Theorems 7.4.3 and 7.4.5, Sections 8.1.2 and 8.1.3]{M}}]\label{compirrpar} $\bigl(\mathcal{P}^h_*\mathcal{E}^{0},\theta\bigr)\to (X,D)$ is a filtered Higgs bundle. Furthermore, if $\xi \in \mathbb{C}^*$ and we denote \smash{$\nabla^{\xi}=D\bigl(\overline{\partial}_E,h\bigr)+\xi^{-1}\theta + \xi \theta^{\dagger_h}$}, then $\bigl(\mathcal{P}^h_*\mathcal{E}^{\xi},\nabla^{\xi}\bigr)\to (X,D)$ is a filtered flat bundle with \[
\textnormal{Irr}_p\bigl(\nabla^{\xi}\bigr)=\left\{\frac{1+|\xi|^2}{\xi}\mathfrak{a} \mid \mathfrak{a}\in \textnormal{Irr}_p(\theta) \right\}.
\]
\end{Theorem}

We now want to state the main results that give the required conditions to go from a filtered Higgs bundle or filtered meromorphic flat bundle to a wild harmonic bundle. To do this, we need to introduce the appropriate stability notions.

\begin{Definition} Let $(_{\boldsymbol{c}}\mathcal{E},\theta)\to (X,D)$ be a $\boldsymbol{c}$-parabolic Higgs bundle. Any subbundle~${\mathcal{H}\subset {_{\boldsymbol{c}}\mathcal{E}}}$ gets an induced $\boldsymbol{c}$-parabolic structure $_{\boldsymbol{c}}\mathcal{H}$, given by $\mathcal{F}_d(_{\boldsymbol{c}}\mathcal{H}|_p):=\mathcal{H}|_p\cap \mathcal{F}_d(_{\boldsymbol{c}}\mathcal{E}|_p)$. We say that~$(_{\boldsymbol{c}}\mathcal{E},\theta)$ is stable (resp.\ semistable), if for every proper non-trivial subbundle $\mathcal{H}$ with $\theta(\mathcal{H})\subset \mathcal{H}\otimes \Omega^1_X(*D)$ we have that
\begin{equation*}
 \frac{\textnormal{pdeg}(_{\boldsymbol{c}}\mathcal{H})}{\textnormal{rank}(_{\boldsymbol{c}}\mathcal{H})}<\frac{\textnormal{pdeg}(_{\boldsymbol{c}}\mathcal{E})}{\textnormal{rank}(_{\boldsymbol{c}}\mathcal{E})} \qquad \left(\textnormal{resp.} \quad \frac{\textnormal{pdeg}(_{\boldsymbol{c}}\mathcal{H})}{\textnormal{rank}(_{\boldsymbol{c}}\mathcal{H})}\leq\frac{\textnormal{pdeg}
 (_{\boldsymbol{c}}\mathcal{E})}{\textnormal{rank}(_{\boldsymbol{c}}\mathcal{E})}\right) .
\end{equation*}
Furthermore, we say that $(_{\boldsymbol{c}}\mathcal{E},\theta)$ is polystable if it is semistable, and $(_{\boldsymbol{c}}\mathcal{E},\theta)=\bigoplus_{i}(_{\boldsymbol{c}}\mathcal{E}_i,\theta_i)$ with each $(_{\boldsymbol{c}}\mathcal{E}_i,\theta_i)$ stable and satisfying
\begin{equation*}
 \frac{\textnormal{pdeg}(_{\boldsymbol{c}}\mathcal{E}_i)}{\textnormal{rank}(_{\boldsymbol{c}}\mathcal{E}_i)}=\frac{\textnormal{pdeg}(_{\boldsymbol{c}}\mathcal{E})}{\textnormal{rank}(_{\boldsymbol{c}}\mathcal{E})}.
\end{equation*}
Similarly, we have the definitions of stable, semistable, and polystable for a $\boldsymbol{c}$-parabolic flat bundle $(_{\boldsymbol{c}}\mathcal{E},\nabla)$. We say that a filtered Higgs bundle or filtered flat bundle is stable/semistable /polystable if any of its $\boldsymbol{c}$-truncations is stable/semistable/polystable. This is well defined (i.e., it does not depend on $\boldsymbol{c}$).
\end{Definition}

Now we can state the following known results due to Biquard and Boalch \cite{BB04}.

\begin{Theorem} \label{HMFB} Let $(\mathcal{P}_*\mathcal{E},\nabla)\to (X,D)$ be a filtered flat bundle, and let $(\mathcal{E}|_{X\setminus D},\nabla) \to X\setminus D$ be its restriction to $X\setminus D$. Then there is a harmonic metric $h$ for $(\mathcal{E}|_{X\setminus D},\nabla)$ adapted to the filtration $\bigl($i.e., $\mathcal{P}^h_*(\mathcal{E}|_{X\setminus D})=\mathcal{P}_*\mathcal{E}\bigr)$ if and only if $(\mathcal{P}_*\mathcal{E},\nabla)\to (X,D)$ is polystable with~$\textnormal{pdeg}(\mathcal{P}_*\mathcal{E})=0$. The harmonic metric is unique up to multiplication by positive constants.
\end{Theorem}

\begin{Remark} A harmonic metric $h$ for a flat bundle $(E,\nabla)$ is a hermitian metric for $E$ such that the decomposition $\nabla=D+\Phi$ into a unitary connection $D$ and a self adjoint endomorphism~$\Phi$ satisfies that $\bigl(E,D^{(0,1)},\Phi^{(1,0)}\bigr)$ is a Higgs bundle, where $D^{(0,1)}$ \big(resp.\ $\Phi^{(1,0)}$\big) denotes the $(0,1)$ (resp.\ $(1,0)$ part) of $D$ (resp.\ $\Phi$).
\end{Remark}

We also have the Higgs bundle version of the previous theorem.

\begin{Theorem} \label{HMHB} Let $(\mathcal{P}_*\mathcal{E},\theta)\to (X,D)$ be a filtered Higgs bundle and let $(\mathcal{E}|_{X\setminus D},\theta) \to X\setminus D$ be its restriction to $X\setminus D$. Then there is a harmonic metric $h$ for \smash{$(\mathcal{E}|_{X\setminus D},\theta)$} adapted to the filtration $\bigl($i.e., \smash{$\mathcal{P}^h_*(\mathcal{E}|_{X\setminus D})=\mathcal{P}_*\mathcal{E}\bigr)$} if and only if $(\mathcal{P}_*\mathcal{E},\theta)\to (X,D)$ is polystable with~$\textnormal{pdeg}(\mathcal{P}_*\mathcal{E})=0$. The harmonic metric is unique up to multiplication by positive constants.
\end{Theorem}

\subsection{Definition of our set of framed wild harmonic bundles}\label{secdeffwhb}

Now we focus on the case of interest to us. For the rest of the paper, we take $X=\CP$, we fix a holomorphic coordinate $z$ on $\mathbb{C}\subset \CP$, and we take $D=\{\infty\}$.

\begin{Definition}
We denote by $\mathcal{H}$ the set of rank $2$, wild harmonic bundles $\bigl(E,\overline{\partial}_E,\theta,h\bigr)\to \bigl(\CP\setminus \{\infty\}\bigr)$, such that $\textnormal{Tr}(\theta)=0$ and $\textnormal{Det}(\theta)=-\bigl(z^2+2m\bigr){\rm d}z^2$ for some $m \in \mathbb{C}$.

\end{Definition}

To illustrate some examples, we show the following lemma.

\begin{Lemma}\label{C1} For every $m\in \mathbb{C}^*$, there is an element of $\mathcal{H}$ such that $\operatorname{Det}(\theta)=-\bigl(z^2+2m\bigr){\rm d}z^2$. Furthermore, if \smash{$m^{(3)}\in \bigl(-\frac12,\frac12\bigr]$}, we can find an element of $\mathcal{H}$ whose associated $\frac{1}{2}$-parabolic Higgs bundle has parabolic weights \smash{$\pm m^{(3)}$} if \smash{$m^{(3)}\neq \frac{1}{2}$}, and with parabolic weights equal to $\frac{1}{2}$ if~\smash{$m^{(3)}=\frac{1}{2}$}.
\end{Lemma}

\begin{proof}
Let $m\in \mathbb{C}^*$, and consider the trivial bundle $V:=\bigl(\CP \setminus \{\infty\}\bigr)\times \mathbb{C}^2 \to \bigl(\CP \setminus \{\infty\}\bigr)$ with canonical global frame $(e_1,e_2)$ and holomorphic structure $\overline{\partial}_V:=\overline{\partial}$. Let $\theta$ be given in this frame by
\begin{equation*}
\theta= \begin{bmatrix}
 0& 1\\
 z^2+2m & 0\\
 \end{bmatrix}{\rm d}z.
\end{equation*}
The plan is to extend $V$ to a $\frac{1}{2}$-parabolic Higgs bundle on $\bigl(\CP,\infty\bigr)$ in such a way that we can apply Theorem \ref{HMHB}.

Since the eigenvalues of $\theta$ near $\infty$ are different and unramified, we can find a punctured neighborhood $U^*_{\infty}:=U_{\infty}\setminus \{\infty\}$ such that there is a holomorphic eigenframe $(\eta_1,\eta_2)$ of $\theta$. We can furthermore assume that we pick the holomorphic eigenframe $(\eta_1,\eta_2)$ such that $e_1\wedge e_2 = \eta_1\wedge \eta_2$.

Let $E\to \CP$ be the holomorphic vector bundle defined by extending $V$ using the frame $(\eta_1,\eta_2)$ near $\infty$. Explicitly, this means that holomorphic sections of $E$ in a neighborhood $U_{\infty}$ of $\infty$ are of the form $f^1\eta_1+f^2\eta_2$, where $f^i$ are holomorphic functions on $U_{\infty}$. Because of the construction, $e_1\wedge e_2$ on $\CP \setminus \{\infty\}$ and $\eta_1\wedge \eta_2$ on $U_{\infty}$ glue together to give a global frame of~$\operatorname{Det}(E)$, so~$\operatorname{deg}(E)=0$.

Now we explain how to put several possible parabolic structures on $E$. On a punctured neighborhood $U_{\infty}^*$ of $\infty$, we know that $E|_{U_{\infty}^*}=L_1\oplus L_2$, where $L_i$ are the eigenlines of $\theta$ near~$\infty$. These bundles extend to line bundles over $U_{\infty}$ and we denote them by $L_i$ as before. First let us consider the filtration of $E_{\infty}$ given by attaching the parabolic weight $m^{(3)} \in \bigl(-\frac12,\frac12\bigr)$ to~$L_1|_{\infty}$ and the parabolic weight $-m^{(3)}$ to $L_2|_{\infty}$. With these choices, $E$ acquires a $\frac{1}{2}$-parabolic structure with
\begin{equation*}
 \operatorname{pdeg}(E)=\operatorname{deg}(E)-m^{(3)}+m^{(3)}=0.
\end{equation*}

If we let $w$ be the holomorphic coordinate related to $z$ by $w=\frac{1}{z}$ and let $\text{Irr}(\theta)_{\infty}:=\bigl\{\pm \frac{1}{w^2}\bigr\}$, then we clearly have a splitting
\begin{equation*}
 \bigl(E,\overline{\partial}_E,\theta\bigr)|_{U_{\infty}}=\bigoplus_{\mathfrak{a}\in \text{Irr}(\theta)_{\infty}} \bigl(E_{\mathfrak{a}},\overline{\partial}_{\mathfrak{a}},\theta_{\mathfrak{a}}\bigr)
\end{equation*}
with $\theta_{\mathfrak{a}}-d\mathfrak{a}\cdot \text{Id}_{E_{\mathfrak{a}}}$ having at most a simple pole at $z=\infty$ (and where the $\bigl(E_{\mathfrak{a}},\overline{\partial}_{\mathfrak{a}}\bigr)$ correspond to the $L_i$, in some order). We then conclude that \smash{$\bigl(E,\overline{\partial}_E,\theta\bigr) \to \bigl(\CP,\infty\bigr)$} is a $\frac{1}{2}$-parabolic Higgs bundle, and that the parabolic structure is compatible with the irregular decomposition.

The parabolic Higgs bundle $\bigl(E,\overline{\partial}_E,\theta\bigr)$ is clearly stable, since there are no global eigenlines preserved by $\theta$: if there were, this would imply that there is a global branch over $\mathbb{C}\subset \CP$ of~\smash{$\sqrt{z^2+2m}$} when $m\neq 0$. Hence, by Theorem $\ref{HMHB}$ we get a harmonic metric $h$.

Now let us consider another possible extension of $V$ with a parabolic structure that allows for a harmonic metric. We now extend $V$ by eigenframes $(\eta_1,\eta_2)$ satisfying $ze_1\wedge e_2=\eta_1\wedge \eta_2$. With this extension, we get $\deg(E)=1$. We choose the trivial filtration of $E_{\infty}$ with weight~$m^{(3)}=\frac{1}{2}$. With these choices, we get
\begin{equation*}
 \operatorname{pdeg}(E)=\operatorname{deg}(E)-\frac{1}{2}2=0.
\end{equation*}

By the same argument from before, we get a harmonic metric for the $\frac{1}{2}$-parabolic Higgs bundle $\bigl(E,\overline{\partial}_E,\theta\bigr)$.
\end{proof}

It is also possible to explicitly build wild harmonic bundles in $\mathcal{H}$ in the case $m=0$. This is explained in Appendix \ref{ae}.

We now define the main set of wild harmonic bundles that we will consider.

\begin{Definition}\label{maindef}
We will denote by $\mathcal{H}^{{\rm fr}}$ the set of tuples $\bigl(E,\overline{\partial}_E,\theta,h,g\bigr)$, where
\begin{itemize}\itemsep=0pt
 \item $(E,h)\to \CP$ is an ${\rm SU}(2)$-vector bundle, so in particular comes with a natural volume form $\omega$ trivializing $\operatorname{Det}(E)$.
 \item $\bigl(E,\overline{\partial}_E,\theta,h\bigr)\to \CP \setminus \{\infty\}$ is an element of $\mathcal{H}$, and it is compatible with the ${\rm SU}(2)$ structure of $(E,h,\omega)$ in the sense that $D\bigl(\overline{\partial}_E,h\bigr)(\omega)=0$, where $D\bigl(\overline{\partial}_E,h\bigr)$ denotes the Chern connection.
 \item $g$ is an ${\rm SU}(2)$-frame of $E_\infty$, such that it extends to an ${\rm SU}(2)$-frame in a neighborhood of~$\infty$ where $\theta$ and $\overline{\partial}_E$ have the following form:
\begin{gather}
 \theta= -H\frac{{\rm d}w}{w^3}-mH\frac{{\rm d}w}{w}+\text{regular terms},\label{singHiggs}
\\
 \overline{\partial}_E=\overline{\partial}-\frac{m^{(3)}}{2}H\frac{{\rm d}\overline{w}}{\overline{w}}+\text{regular terms} \qquad \text{for some}\quad m^{(3)}\in \left(-\frac{1}{2},\frac{1}{2}\right],\label{singparabolic}
 \end{gather}
 where $w=\frac{1}{z}$ and $H=\bigl[\begin{smallmatrix}
 1&0\\
 0&-1\end{smallmatrix}\bigr]$.
 We will call such a frame $g$ at $\infty$ a compatible frame for the wild harmonic bundle.
\end{itemize}
\end{Definition}

Notice that the parameter $m$ appearing in equation \eqref{singHiggs} is the same as the $m$ parameter appearing in the condition $\textnormal{Det}(\theta)=-\bigl(z^2+2m\bigr){\rm d}z^2$. On the other hand, the parameter $m^{(3)}$ in \eqref{singparabolic} should be though as parametrizing the parabolic structures constructed in the proof of Lemma \ref{C1}. Finally, note that the bundles $E$ in the elements of $\mathcal{H}^{{\rm fr}}$ are actually bundles over~$\CP$, while the bundles $E$ in elements of $\mathcal{H}$ are only bundles over $\CP \setminus \{\infty\}$.

\begin{Definition} An isomorphism between two elements $\bigl(E_i,\overline{\partial}_{E_i},\theta_i,h_i,g_i\bigr)\in \mathcal{H}^{{\rm fr}}$, $i=1,2$, is an isomorphism $f\colon E_1\to E_2$ of vector bundles over $\mathbb{C}P^1$ such that $\overline{\partial}_{E_2}\circ f=f\circ \overline{\partial}_{E_1}$, $\theta_2\circ f=f\circ \theta_1$, $f^*h_2=h_1$, $f$ takes $g_1$ to $g_2$, and $f^*\omega_{2}=\omega_{1}$. The set of isomorphism classes of $\mathcal{H}^{{\rm fr}}$ will be denoted by $\mathfrak{X}^{{\rm fr}}$.
\end{Definition}
\begin{Remark}
Later we will consider the set $\mathfrak{X}^{{\rm fr}}$ of isomorphism classes of $\mathcal{H}^{{\rm fr}}$. It is worth mentioning at this point how the points of $\mathcal{M}^{\text{ov}}$ are going to correspond to point in $\mathfrak{X}^{{\rm fr}}$. First, it is not hard to see that $m$ and $m^{(3)}$ are isomorphism invariants of elements of $\mathcal{H}^{{\rm fr}}$ (see, for example, Lemma \ref{mm3inv} below). Furthermore, we will later show that for equivalence classes with fixed $m\neq 0$ and $m^{(3)}$ are a ${\rm U}(1)$-torsor for a certain ${\rm U}(1)$-action acting on the framing (see Proposition \ref{U1A}). The reader should then think of $m$ as being the analog of the base parameter~$z$ of $\mathcal{M}^{\text{ov}}$, \smash{${\rm e}^{2\pi {\rm i}m^{(3)}}$} as being the analog of ${\rm e}^{2\pi {\rm i}x^3}$, and the ${\rm U}(1)$-torsor as being the analog of the ${\rm U}(1)$-fiber of the principal ${\rm U}(1)$-bundle used in the construction of $\mathcal{M}^{\text{ov}}$. This interpretation will come naturally after we construct certain complex coordinates $\mathcal{X}_e$ and $\mathcal{X}_m$ on $\mathfrak{X}^{{\rm fr}}$ in terms of Stokes data, and compare them with the Ooguri--Vafa twistor coordinates $\mathcal{X}_e^{\text{ov}}$ and $\mathcal{X}_m^{\text{ov}}$.
\end{Remark}

\begin{Example}\label{trivWHB} Let us give an explicit (although trivial) example of such a compatibly framed wild harmonic bundle in the case, where $m=m^{(3)}=0$. We take $E\to \CP$ to be the trivial bundle~${E=\CP \times \mathbb{C}^2}$. If $(e_1,e_2)$ denotes the canonical global frame of $E$, we give an~${\rm SU}(2)$ structure to $E$ by considering the hermitian metric $h(e_i,e_j)=\delta_{ij}$ and the volume form ${\omega=e_1\wedge e_2}$. Furthermore, in the canonical frame $(e_1,e_2)$ consider the trivial holomorphic structure $\overline{\partial}_E=\overline{\partial}$, the Higgs field $\theta=zH{\rm d}z$, and the framing at infinity $g=(e_1,e_2)|_{\infty}$. We then have that~$D\bigl(\overline{\partial}_E,h\bigr)={\rm d}$, so Hitchin equations are clearly satisfied, and $\omega$ is parallel with respect to the Chern connection. On the other hand,
the frame $g$ extends to the global frame $(e_1,e_2)$, where $\overline{\partial}_E=\overline{\partial}$, and $\theta=-\frac{1}{w^3}H{\rm d}w$. Hence, we get an element of $\mathcal{H}^{{\rm fr}}$.
\end{Example}

\begin{Example} \label{exmorph} We now give examples of isomorphisms: let $\bigl(E,\overline{\partial}_E,\theta,h,(f_1,f_2)\bigr)\in \mathcal{H}^{{\rm fr}}$. Then it is easy to check that for $c>0$ we have that \smash{$\bigl(E,\overline{\partial}_E,\theta, c\cdot h, \bigl(\sqrt{c^{-1}}f_1,\sqrt{c^{-1}}f_2\bigr)\bigr) \in \mathcal{H}^{{\rm fr}}$}. These two elements are clearly isomorphic by the bundle map \smash{$\sqrt{c^{-1}}\cdot \text{Id}_E$}, where $\text{Id}_E$ is the identity map on $E$.
\end{Example}

In the following, we will explain how to get an element of $\mathcal{H}^{{\rm fr}}$ from an element in $\mathcal{H}$. To show this, we will require the following two results on wild harmonic bundles. The two result are stated with respect to our particular case, the more general statements are in the given references.

\begin{Theorem}[{wild version of Simpson's main estimate \cite[Theorem 7.2.1]{M}}]\label{FT1} Let $\bigl(E,\overline{\partial}_E,\theta,h\bigr)\allowbreak\in \mathcal{H}$, and consider the decomposition into eigenlines near $\infty$
\begin{equation*}
\bigl(E,\overline{\partial}_E,\theta\bigr)|_{U_{\infty}\setminus \{\infty\}}=\bigoplus_{\mathfrak{a}\in \textnormal{Irr}(\theta)_{\infty}} \bigl(E_{\mathfrak{a}},\overline{\partial}_{\mathfrak{a}},\theta_{\mathfrak{a}}\bigr),
\end{equation*}
where $\textnormal{Irr}(\theta)|_{\infty}=\bigl\{\pm \frac{1}{w^2}\bigr\}$. Let $v_{\mathfrak{a}}$ be a section of $\bigl(E_{\mathfrak{a}},\overline{\partial}_{E_{\mathfrak{a}}}\bigr)$ near $\infty$. Then for $\mathfrak{a}\neq \mathfrak{b}$, we have that~${
 |h(v_{\mathfrak{a}},v_{\mathfrak{b}})|\leq C|v_{\mathfrak{a}}|_h|v_{\mathfrak{b}}|_h\exp \bigl(-\epsilon|w|^{-2}\bigr)}$,
where $\epsilon>0$ and $w=\frac{1}{z}$ is a holomorphic coordinate vanishing at $w=0$. In particular, if $v_{\mathfrak{a}}$ and $v_{\mathfrak{b}}$ are sections of $\mathcal{P}^h_c\mathcal{E}^0$ for some $c \in \mathbb{R}$, then $v_{\mathfrak{a}}$ and $v_{\mathfrak{b}}$ are asymptotically exponentially orthogonal near $w=0$.
\end{Theorem}

\begin{Theorem}[{\cite[Proposition 8.1.1]{M}}]\label{FT2} Let $\bigl(E,\overline{\partial}_E,\theta,h\bigr)\in \mathcal{H}$ and let $\bigl(\mathcal{P}^h_{c}\mathcal{E}^{0},\theta\bigr)$ be the associated $\boldsymbol{c}$-parabolic Higgs bundle. Furthermore, let $v=(v_1,v_2)$ be a frame of $\mathcal{P}^h_{c}\mathcal{E}^{0}$ in a neighborhood~$U_{\infty}$ of $\infty$, compatible with the c-parabolic structure. Let $a(v_i)$ denote the parabolic weight corresponding to $v_{i}$, and define the following hermitian metric on $\mathcal{P}^h_{c}E|_{U_{\infty}}$
\begin{equation*}
 h_0(v_{i},v_{j})=\delta_{ij}|w|^{-2a(v_i)}.
\end{equation*}
Furthermore, let \smash{$\widetilde{v_{i}}=v_{i}|w|^{a(v_i)}$}. Then $h$ is mutually bounded with respect to $h_0$ in the sense that there are positive constants $C_1$, $C_2$ such that
$
 C_1<|H(\widetilde{v})|<C_2$,
where $H(\widetilde{v})$ is the matrix with entries $h(\widetilde{v}_{i},\widetilde{v}_{j})$.
\end{Theorem}

Now we can prove the following proposition.

\begin{Proposition} \label{constructionframe} Given $\bigl(E,\overline{\partial}_E,\theta,h\bigr)\in \mathcal{H}$, there is an extension of $(E,h)$ to an ${\rm SU}(2)$-vector bundle over $\CP$, and a compatible framing $g$ of $E_{\infty}$ such that \smash{$\bigl(E,\overline{\partial}_E,\theta,h, g\bigr)\in \mathcal{H}^{{\rm fr}}$}. Further\-mo\-re, if $g=(e_1,e_2)$, \smash{${\rm e}^{{\rm i}\theta}\in {\rm U}(1)$} and we let \smash{${\rm e}^{{\rm i}\theta}\cdot g:=\bigl({\rm e}^{{\rm i}\theta}e_1,{\rm e}^{-{\rm i}\theta}e_2\bigr)$}, then \smash{$\bigl(E,\overline{\partial}_E,\theta,h,\allowbreak {{\rm e}^{{\rm i}\theta}}\cdot g\bigr)\in \mathcal{H}^{{\rm fr}}$}.
\end{Proposition}

\begin{proof}Consider the associated filtered Higgs bundle $\bigl(\mathcal{P}_*^h\mathcal{E}^{0},\theta\bigr)$ (recall Theorem \ref{compirrpar}). We work with the associated $1/2$-parabolic bundle $\mathcal{P}_{1/2}^h\mathcal{E}^{0}$ and pick a holomorphic eigenframe $(v_1,v_2)$ of $\theta$ near $\infty$, compatible with its parabolic structure. The fact that we can do this follows from the compatibility of the irregular decomposition of \smash{$\bigl(\mathcal{P}^h_{1/2}\mathcal{E}^{0},\theta\bigr)$} at $z=\infty$ and the parabolic structure (recall the last point of Definition \ref{filteredhiggsdef}). We order this frame in such a way that
\begin{equation*}
\theta= -H\frac{{\rm d}w}{w^3}-mH\frac{{\rm d}w}{w}+\text{diagonal holomorphic terms.}
\end{equation*}

On the other hand, we claim the following.

\begin{Lemma}\label{connlem}
In the frame $(v_1,v_2)$ the Chern connection $D\bigl(\overline{\partial}_E,h\bigr)$ has the following form:
\begin{equation*}
 D\bigl(\overline{\partial}_E,h\bigr)=D_0 + \text{regular terms},
\end{equation*}
where $D_0=D\bigl(\overline{\partial}_E,h_0\bigr)$ is the Chern connection of the metric $h_0$ from Theorem {\rm\ref{FT2}}
\begin{equation*}
 D_0(v_i)=-a(v_i)v_i\frac{{\rm d}w}{w}.
\end{equation*}

\end{Lemma}
\begin{proof}
To see this, we follow a similar argument to \cite[Proposition 10.3.3]{M}, but modified for our special, simpler situation.

If we denote by $\hat{h}$ and $\hat{h}_0$ the matrices corresponding to $h$ and $h_0$ in the frame $(v_1,v_2)$, we then have $\hat{h}=\hat{h}_0\cdot g$ for some matrix valued $g$. Since $(v_1,v_2)$ is a holomorphic frame and $D$ and~$D_0$ are Chern connections for the same holomorphic structure $\overline{\partial}_E$, the connection matrix~$A$ of~$D$ in the frame $(v_1,v_2)$ (resp.\ $A_0$ of $D_0$) is given by
$A=\hat{h}^{-1}\partial \hat{h}$ (resp. $A_0=\hat{h}_0^{-1}\partial \hat{h}_0$).
A~quick computation then shows that $\hat{h}=\hat{h}_0g$ implies that
\begin{equation*}
 A=A_0+g^{-1}(\partial g +[A_0,g]).
\end{equation*}
Hence, $D$ and $D_0$ are related by the equation $D=D_0 +g^{-1}\partial_{h_0}(g)$, where ${\partial_{h_0}(g)=\partial g+[A_0,g]}$. The off-diagonal terms of $g^{-1}\partial_{h_0}(g)$ are exponentially decreasing near $w=0$ by \cite[Lem\-ma~10.1.3]{M}, so we only need to show that the diagonal terms of $g^{-1}\partial_{h_0}(g)$ are regular at $w=0$.

Because the frame $(v_1,v_2)$ is asymptotically exponentially orthogonal, it is easy to conclude that $\bigl[\theta,\theta^{\dagger_h}\bigr]$ is regular at $z=\infty$. Hence, by the Hitchin equation we get that $F\bigl(D\bigl(\overline{\partial}_E,h\bigr)\bigr)$ is regular at $\infty$. On the other hand, since $D_0$ is flat, we have that $F(D)=\overline{\partial}\bigl(g^{-1}\partial_{h_0}(g)\bigr)$. This lets us conclude that we can write $g^{-1}\partial_{h_0}(g)=\rho + \chi$, where $\rho$ is regular near $\infty$ and $\chi$ is holomorphic and defined in a punctured neighborhood of $\infty$.

By the argument in \cite[Proposition 10.3.3]{M}, we get that $g^{-1}\partial_{h_0}(g)$ is square integrable relative to $h_0$. Using this, and the fact that the off-diagonal entries of $g^{-1}\partial_{h_0}(g)$ are exponentially decreasing near $w=0$, it follows that the diagonal elements of $g^{-1}\partial_{h_0}(g)$ are square integrable. Since $\rho$ is regular near $\infty$, we conclude that the diagonal elements of $\chi$ are also square integrable. Furthermore, since the diagonal elements of $\chi$ are holomorphic functions on $U_{\infty}\setminus \{\infty\}$, we then conclude that they must extend to holomorphic functions on $U_{\infty}$. This shows the regularity at~${w=0}$ of the diagonal terms of $g^{-1}\partial_{h_0}(g)$.
\end{proof}

Now we go back to the proof of Proposition \ref{constructionframe}. By the wild version of Simpson's main estimate (see Theorem \ref{FT1}), $(v_1,v_2)$ are asymptotically exponentially orthogonal near $w=0$, so it is easy to check that if we orthonormalize the frame $(v_1,v_2)$ and obtain the (non-holomorphic) orthonormal frame $(e_1,e_2)$, then in this frame we get
\begin{equation}\label{higgscan}
 \theta=-H\frac{{\rm d}w}{w^3}-mH\frac{{\rm d}w}{w}+\text{regular terms},
\end{equation}
where the regular terms are no longer holomorphic. Similarly, using Theorem \ref{FT1} and Lem\-ma~\ref{connlem}, one can show that in the frame $(e_1,e_2)$, the Chern connection acquires the following form:
\begin{equation}\label{Dlemexp}
 D={\rm d} - \frac{1}{2}\begin{bmatrix} a(v_1) & 0 \\
 0 & a(v_2) \\
 \end{bmatrix}\left(\frac{{\rm d}w}{w}- \frac{{\rm d}\overline{w}}{\overline{w}}\right) + \text{regular}.
\end{equation}
We prove this statement in Appendix \ref{AA}.

Since we are picking the weights in $\bigl(-\frac12,\frac12\bigr]$ and the parabolic degree must be $0$ (and hence an integer), we have the following possibilities for $a(v_i)$: we either have $a(v_1)=-a(v_2)$, or $a(v_1)=a(v_2)=\frac12$.

If $a(v_1)=- a(v_2)$, let $m^{(3)}:=a(v_1)$, so that we can rewrite \eqref{Dlemexp} as
\begin{equation}\label{Dcan}
 D={\rm d}-m^{(3)}H\left(\frac{{\rm d}w}{2w}-\frac{{\rm d}\overline{w}}{2\overline{w}}\right) + \text{regular terms}.
\end{equation}

In the case, where $a(v_1)=a(v_2)=\frac12$, let $m^{(3)}:=\frac12$, consider the eigenframe $(v_1,wv_2)$, and then consider the associated orthonormal frame. We abuse notation and also denote it by~$(e_1,e_2)$. In this frame, we now get again the expressions \eqref{higgscan} and \eqref{Dcan}.

We now use the framing $(e_1,e_2)$ to do a unitary extension of the hermitian bundle $(E,h)\to \CP \setminus \{\infty\}$ to a hermitian bundle $(E,h)\to \CP$, where $(e_1,e_2)$ gives a unitary trivialization on a neighborhood $U_{\infty}$ of $\infty$. We denote the extended frame over $U_{\infty}$ induced by $(e_1,e_2)$ by the same notation. By construction, $g=(e_1,e_2)|_{\infty}$ extends to a unitary frame where $\theta$ and $\overline{\partial}_E$ have the appropriate forms \eqref{singHiggs} and \eqref{singparabolic}.

To finish the proof of Proposition \ref{constructionframe}, we need to show that we can reduce the structure group of $(E,h)\to \CP$ to ${\rm SU}(2)$ and that the induced volume form $\omega$ satisfies $D\bigl(\overline{\partial}_E,h\bigr)(\omega)=0$. We start with the following lemma.

\begin{Lemma} $D\bigl(\overline{\partial}_E,h\bigr)$ induces a flat connection on $\textnormal{Det}(E)\to \CP$.
\end{Lemma}

\begin{proof}Notice that by taking the trace of the Hitchin equation, we get that $\operatorname{Tr}(F(D))=0$, so that the connection $\operatorname{Det}(D)$ on $\operatorname{Det}(E)\to \CP \setminus \{\infty\}$ induced by $D\bigl(\overline{\partial}_E,h\bigr)$ is flat. On the other hand, since the singularity of $D$ in the frame $(e_1,e_2)$ is traceless (recall \eqref{Dcan}), we get that the connection form of $\operatorname{Det}(D)$ in the frame $e_1\wedge e_2$ is actually smooth at $w=0$, and hence defines a flat connection on $\operatorname{Det}(E)\to \CP$.
\end{proof}

If $D={\rm d}+A$ in the frame $(e_1,e_2)$, then $\operatorname{Det}(D)={\rm d}+\operatorname{Tr}(A)$ in the frame $e_1\wedge e_2$, with $\operatorname{Tr}(A)$ a well defined and closed form (since $\operatorname{Det}(D)$ is flat) on a neighborhood of $\infty$. Furthermore, since the frame $(e_1,e_2)$ is unitary and $D$ is the Chern connection, we get that $A$ is valued in~$\mathfrak{u}(2)$, and hence $\operatorname{Tr}(A)$ is valued in $\mathfrak{u}(1)={\rm i}\mathbb{R}$. We now perform the unitary (and non-singular) diagonal gauge transformation \smash{$e_i\to {\rm e}^{-\frac{1}{2}\int^{p}_{\infty}\operatorname{Tr}(A)}e_i:=\widetilde{e}_i$}, where the integral is performed along any path from $z=\infty$ to $p$ in $U_{\infty}$. The local frame $\widetilde{e}_1\wedge \widetilde{e}_2$ of $\operatorname{Det}(E)\to U_{\infty}$ is then flat with respect to $\operatorname{Det}(D)$, and by performing parallel transport, we get a global flat frame $\omega$ of~${\operatorname{Det}(E)\to \CP}$. Hence we can reduce the structure group of $(E,h)\to \CP$ to ${\rm SU}(2)$. Notice that in our ${\rm SU}(2)$-frame $(\widetilde{e}_1,\widetilde{e}_2)$, we have that $\theta$ and $\overline{\partial}_E$ still have the required form. Hence, $g=(e_1,e_2)|_{\infty}=(\widetilde{e}_1,\widetilde{e}_2)|_{\infty}$ satisfies our condition for a compatible frame. We then conclude that given $\bigl(E,\overline{\partial}_E,\theta,h\bigr) \in \mathcal{H}$, we can produce $\bigl(E,\overline{\partial}_E,\theta,h,g\bigr)\in \mathcal{H}^{{\rm fr}}$.

The remaining last statement that $\bigl(E,\overline{\partial}_E,\theta,h, {\rm e}^{{\rm i}\theta}\cdot g\bigr)\in \mathcal{H}^{{\rm fr}}$, where ${\rm e}^{{\rm i}\theta}\cdot g:=\bigl({\rm e}^{{\rm i}\theta}e_1,{\rm e}^{-{\rm i}\theta}e_2\bigr)$ and ${\rm e}^{{\rm i}\theta}\in {\rm U}(1)$ follows easily.
\end{proof}

\subsection{Twistor coordinates Part 1: Preliminaries and definition}\label{sec3.4}

We now begin with the construction of the twistor coordinates for $\mathfrak{X}^{{\rm fr}}$. Roughly speaking, to each element of $\mathcal{H}^{{\rm fr}}$ we can associate a ``compatibly framed'' filtered flat bundle, and to this object we can associate ``generalized monodromy data'' (more commonly known as Stokes data). The twistor coordinates for $\mathfrak{X}^{{\rm fr}}$ will be built from Stokes data.

\begin{subsubsection}{Stokes data}\label{GMD}
In this section, we recall the definition of Stokes data, following mostly \cite{Bir, B, B2, Wasow,Wit}. We assume that $X=\CP$ and $D\subset \CP$ is a point, although the definitions and results hold in more general settings (see, for example, the aforementioned references). We also fix a holomorphic bundle $E\to \CP$ of $\text{rank}(E)=n$, with a meromorphic connection $\nabla$ with poles along $D$. Since~$X$ is a~curve and $\nabla$ is meromorphic, it is automatically flat.

If we choose a local holomorphic coordinate $w$ such that $D$ is given by $w=0$, and a local holomorphic trivialization of the bundle $E$ near $D$, then the connection $\nabla$ has the form $\nabla={\rm d}+A$, where
\begin{equation} \label{normalform1}
 A=A_k\frac{{\rm d}w}{w^k}+A_{k-1}\frac{{\rm d}w}{w^{k-1}}+ \dots +A_1\frac{{\rm d}w}{w} + \text{holomorphic $(1,0)$ terms},
\end{equation}
and $A_i\in \operatorname{End}(\mathbb{C}^{n})$ for $i=1,2,\dots ,k$. We will assume that $k>1$ and that $A_k$ is generic in the sense that it is diagonalizable with distinct eigenvalues.

We will say that $(E,\nabla)$ is generic if the leading coefficient $A_k$ of the connection is generic in some holomorphic coordinate vanishing at $D$ and some holomorphic trivialization. It is easy to check that the order of the pole and the fact that the leading coefficient is generic does not depend on the holomorphic coordinate vanishing at $D$ and the holomorphic trivialization.

\begin{Definition}
A compatible framing at $D$ for a generic $(E,\nabla)\to (X,D)$ is a frame $g$ for~$E_{D}$, such that in some (and hence any) extension of $g$ to a local holomorphic trivialization, we have that the leading coefficient $A_k$ of the singularity of $\nabla$ is diagonal.
\end{Definition}

\begin{Lemma}[{\cite[Lemma 1]{B2}}] \label{FGT}
Let $(E,\nabla,g)\to (X,D)$ be a compatibly framed connection, and consider a local holomorphic trivialization $\tau$ extending $g$. In the local frame $\tau$, let $A_k$ be the leading coefficient of $\nabla$, as in \eqref{normalform1}. Then there is a unique formal gauge transformation~${\hat{F}\in {\rm GL}(n,\mathbb{C})[[w]]}$ and unique diagonal elements $A^0_j\in \operatorname{End}(\mathbb{C}^n)$, such that $\hat{F}(0)=1$, and such that in the formal frame \smash{$\tau \cdot \widehat{F}$}, the connection looks like
\begin{equation*}
 {\rm d}+A^0:={\rm d}+A_{k}^0\frac{{\rm d}w}{w^k}+A_{k-1}^0\frac{{\rm d}w}{w^{k-1}}+ \dots +A_1^0\frac{{\rm d}w}{w},
\end{equation*}
with $A_k^0=A_k$. The $A^0_j$ only depend on the compatible framing $g$ and not on the extension $\tau$.
\end{Lemma}

\begin{Definition}\label{defFT}
We will call the \smash{$A^0_j$} appearing in the lemma above the formal type of $(E,\nabla,g)$, and we will call $\Lambda:=A_1^0$ the exponent of formal monodromy. If we write the formal type as~${A^0={\rm d}Q(w)+\Lambda\frac{{\rm d}w}{w}}$, with $Q(w)$ a diagonal matrix with entries in $w^{-1}\mathbb{C}\bigl[w^{-1}\bigr]$, we will say that $(Q,\Lambda)$ specifies the formal type.

\end{Definition}

\begin{Definition} Fix a formal type $(Q,\Lambda)$. Let $q_i(w)$ denote the $i$-th diagonal component of the leading term \smash{$\frac{A_k}{(k-1)w^{k-1}}$} of $-Q(w)$, and let $q_{ij}(w):=q_i(w)-q_j(w)$. Now let ${\rm e}^{{\rm i}\theta} \in S^1$ and $r_\theta$ the ray going from $w=0$ to ${\rm e}^{{\rm i}\theta}$. We will say that $r_\theta$ is an anti-Stokes ray if $q_{ij}(w)<0$ on $r_\theta$ for some ordered pair $(i,j)$. Furthermore, we will say $r_\theta$ is a Stokes ray if $\textnormal{Re}(q_{ij}(w))=0$ on $r_{\theta}$ for some $i$ and $j$.
\end{Definition}

For simplicity, and because it will be the case of interest to us, we will further restrict to the case where $\text{rank}(E)=2$. Notice that in this case, if $Q$ has a pole of order $k-1$, then there are~${2k-2}$ anti-Stokes rays ($k-1$ associated to the ordered pair $(1,2)$ and the other $k-1$ associated to $(2,1)$) and $2k-2$ Stokes rays. In Figure~\ref{fig7}, we illustrate a specific case for $k=3$, which will be of interest for us in our application to wild harmonic bundles.

\begin{figure}
 \centering
 \includegraphics[width=5cm]{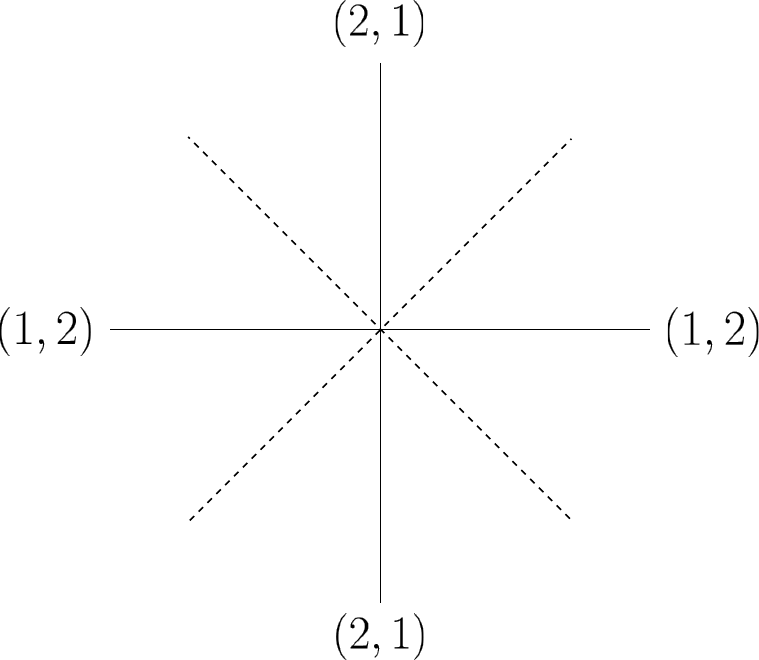}
 \caption{We illustrate the case where we put $q_1=-q_2=-\frac{1}{2w^2}$. The bold rays denote the anti-Stokes rays, while the dotted rays denote the Stokes rays.} \label{fig7}
\end{figure}

Now choose one of the anti-Stokes rays as the first anti-Stokes ray, and let $r_{\theta_i}$ be the anti-Stokes rays for $i=1,2,\dots ,2k-2$, numbered in a counterclockwise manner. With this choice, we let $\operatorname{Sect}_i$ be the sector from $r_{\theta_i}$ to $r_{\theta_{i+1}}$. Furthermore, let $\widehat{\operatorname{Sect}}_i$ be the extended sector from the ray through $\theta_i - \frac{\pi}{2k-2}$ to the ray through $\theta_{i+1} + \frac{\pi}{2k-2}$. Notice that while the sectors $\operatorname{Sect}_i$ are determined by anti-Stokes rays, the extended sectors \smash{$\widehat{\operatorname{Sect}}_i$} are determined by Stokes rays.

\begin{Theorem}[{\cite[Theorem 3.1]{B}}] Let $(E,\nabla,g)\to (X,D)$ be a compatibly framed connection with formal type given by $A^0={\rm d}Q+\Lambda \frac{{\rm d}w}{w}$. Furthermore, let $\tau$ be an extension of $g$ to a local holomorphic framing near $D$, and \smash{$\widehat{F}$} the formal gauge transformation from Lemma {\rm\ref{FGT}}. Then for each $i=1,2,\dots ,2k-2$ and for a sufficiently small disk $B$ centered at $w=0$, there is a~unique invertible matrix \smash{$\Sigma_i\bigl(\widehat{F}\bigr)$} of holomorphic functions defined on \smash{$B\cap \widehat{\textnormal{Sect}}_i$}, such that in the sectorial frame \smash{$\tau \cdot \Sigma_i\bigl(\widehat{F}\bigr)$} we have that
$
 \nabla={\rm d} + {\rm d}Q +\Lambda \frac{{\rm d}w}{w}$.
Each \smash{$\Sigma_i\bigl(\widehat{F}\bigr)$} is asymptotic\footnote{If $\widehat{F}=\sum_{j=0}^{\infty}F_jw^j$ with $F_j\in \operatorname{End}\bigl(\mathbb{C}^2\bigr)$, then the statement that $\Sigma_i\bigl(\widehat{F}\bigr)$ is asymptotic to $\widehat{F}$ as $w\to 0$ means that for each $n\in \mathbb{N}$, we have \smash{$\big|\Sigma_i\bigl(\widehat{F}\bigr)(w)-\sum_{j=0}^n F_jw^j\big|=\mathcal{O}\bigl(|w|^{n+1}\bigr)$}.} to \smash{$\widehat{F}$} as $w \to 0$ along \smash{$B\cap \widehat{\textnormal{Sect}}_i$}.

\end{Theorem}

Now let $(E,\nabla,g)\to (X,D)$ be a compatibly framed connection. If $w$ is a holomorphic coordinate vanishing at $D$, then $g$ extends to a holomorphic frame $\tau$ near $D$ where $\nabla$ has the form
\begin{equation*}
 \nabla={\rm d}+A={\rm d}+ {\rm d}Q+\Lambda \frac{{\rm d}w}{w}+ \text{regular holomorphic terms},
\end{equation*}
and $Q$ and $\Lambda$ specify the formal type of $(E,\nabla,g)$. Strictly speaking, we only know that $g$ has an extension where the leading coefficient of ${\rm d}Q$ is diagonal; however, under our assumption that the leading coefficient is generic, it is easy to see that we can extend the frame $g$ to one where the whole singular part of $\nabla$ is diagonal. This diagonalized singular part must coincide with the formal type ${\rm d}Q+\Lambda \frac{{\rm d}w}{w}$ (see the proof of \cite[Lemma 1]{B2}).

Let $\widehat{F}$ be the formal gauge transformation such that in the frame $\tau \cdot \widehat{F}$ the connection has the form
$\nabla={\rm d}+{\rm d}Q+\Lambda \frac{{\rm d}w}{w}$.

With $\Sigma_i\bigl(\widehat{F}\bigr)$ we can define sectorial frames of flat sections of $\nabla$ in the following way: fix a~branch of the logarithm with branch cut along one of the anti-Stokes rays determined by the leading term of $-Q$. We will call that anti-Stokes ray $r_{\theta_1}$, and we will number the rest in a~counterclockwise manner. With these choices, we get a frame of flat sections for $\nabla={\rm d}+A$ on~$\widehat{\operatorname{Sect}}_i$ by writing $\Phi_i=\tau \cdot \Sigma_i\bigl(\widehat{F}\bigr)w^{-\Lambda} {\rm e}^{-Q}$. We will use the following convention for the $w^{-\Lambda}$: it uses the choice of the branch of the logarithm if $i\neq 1,2k-2$; it uses the analytically continued branch from $\operatorname{Sect}_1$ to $\operatorname{Sect}_{2k-2}$ for $\Phi_1$; and it uses the analytically continued branch from $\operatorname{Sect}_{2k-2}$ to $\operatorname{Sect}_{1}$ for $\Phi_{2k-2}$.

\begin{Definition} \label{defStokesmatrices} For $i\neq 2k-2$, we will denote the transition function between the flat frame~$\Phi_{i}$ to the flat frame $\Phi_{i+1}$ on $\widehat{\textnormal{Sect}}_{i}\cap\widehat{\textnormal{Sect}}_{i+1}$ by $S_{i}$; and for $i=2k-2$ we will denote the transition function between the flat frame $\Phi_{2k-2}$ to the flat frame $\Phi_1 \cdot M_0$ on~${\widehat{\textnormal{Sect}}_{2k-2}\cap\widehat{\textnormal{Sect}}_{1}}$ by~$S_{2k-2}$, where $M_0={\rm e}^{-2\pi {\rm i} \Lambda}$ is the formal monodromy in the counterclockwise manner.
Hence, for $i\neq 2k-2$, we have \smash{$S_{i}={\rm e}^Qw^{\Lambda}\Sigma_{i}\bigl(\widehat{F}\bigr)^{-1}\Sigma_{i+1}\bigl(\widehat{F}\bigr)w^{-\Lambda}{\rm e}^{-Q}$}, and for $i=2k-2$, we have~\smash{$S_{2k-2}={\rm e}^Qw^{\Lambda}\Sigma_{2k-2}\bigl(\widehat{F}\bigr)^{-1}\Sigma_{1}\bigl(\widehat{F}\bigr)w^{-\Lambda}{\rm e}^{-Q}M_0$}. We will call the matrices $S_i$ the Stokes matrices.
\end{Definition}

\begin{Remark} In the last expression for $S_{2k-2}$, we are abusing notation: the $w^{\Lambda}$ on the left is analytically continued from $\operatorname{Sect}_{2k-2}$ to $\operatorname{Sect}_1$, while the one appearing on the right is analytically continued from $\operatorname{Sect}_1$ to $\operatorname{Sect}_{2k-2}$.
\end{Remark}

The following proposition is a well-known property of Stokes data.

\begin{Proposition} The Stokes matrices $S_i$ are constant and unipotent.
\end{Proposition}

Furthermore, since we are working with a meromorphic connection $\nabla$ on $E\to \CP$ with poles along $D=\infty$, monodromy of a $\nabla$ along a loop around $D=\infty$ must be trivial, which implies the following:

\begin{Proposition}\label{MonodromyStokes} The Stokes matrices satisfy the relation $S_1S_2\cdots S_{2k-2}M_0^{-1}=1$.
\end{Proposition}
\end{subsubsection}

\begin{subsubsection}{Associated compatibly framed flat bundles}

In this section, we will associate a $\mathbb{C}^*$-family of ``compatibly framed, filtered flat bundles'' to an element of $\mathcal{H}^{{\rm fr}}$. By this we mean a tuple $\bigl(\mathcal{P}_*\mathcal{E}^{\xi},\nabla^{\xi}, \tau_*^{\xi}\bigr)$ for each $\xi\in \mathbb{C}^*$, where $\bigl(\mathcal{P}_*\mathcal{E}^{\xi},\nabla^{\xi}\bigr) \to \bigl(\CP,\infty\bigr)$ is a filtered flat bundle,\footnote{Given $\xi \in \mathbb{C}^*$ and $\bigl(E,\overline{\partial}_E,\theta,h,g\bigr)\in \mathcal{H}^{{\rm fr}}$ with corresponding parameters $m\in \mathbb{C}$ and $m^{(3)}\in \bigl(-\frac12,\frac12\bigr]$ as in~\eqref{singHiggs} and \eqref{singparabolic}, we will see that the induced $a$-parabolic flat bundle $\bigl(\mathcal{P}^h_a\mathcal{E}^{\xi},\nabla^{\xi}\bigr)$ has a parabolic structure that depends on the choices of $m$, \smash{$m^{(3)}$} and $\xi$. Since we will vary these parameters in the sequel, it will be more convenient to consider filtered flat bundles $\bigl(\mathcal{P}_*^h\mathcal{E}^{\xi},\nabla^{\xi}\bigr)$, rather than a specific $a$-truncation.} and \smash{$\tau_a^{\xi}$} is a compatible frame for the bundle with meromorphic flat connection $(\mathcal{P}_a\mathcal{E}^{\xi},\nabla^{\xi})$ for each $a\in \mathbb{R}$.

Let $\bigl(E,\overline{\partial_E},\theta,h,g\bigr)\in \mathcal{H}^{{\rm fr}}$, let $a\in \mathbb{R}$, and let $w$ be a holomorphic coordinate related to $z$ by $w=\frac{1}{z}$. By following the argument given in \cite[Sections 7 and 8]{BB04}, for $\xi \in \mathbb{C}^*$ we can find a~holomorphic frame \smash{$\tau_a^{\xi}$} of \smash{$\mathcal{E}^{\xi}$} in a neighborhood of $z=\infty$ of the following form:
\begin{equation}\label{twistcompframe}
 \tau_a^{\xi}(w)=(e_1,e_2)\cdot g_{\xi}(w) |w|^{(m^{(3)}+2\xi \overline{m})H}w^{N(a)}\exp \left(\frac{\overline{\xi}H}{2w^2}-\frac{\xi H}{2\overline{w}^2}\right),
\end{equation}
where
\begin{itemize}\itemsep=0pt
 \item $(e_1,e_2)$ is an extension of the frame given by $g$, to an ${\rm SU}(2)$ framing in neighborhood of~$\infty$, satisfying the properties of Definition \ref{maindef}.
 \item $g_{\xi}(w)$ is a gauge transformation defined in a neighborhood of $w=0$ that gauges away the regular $(0,1)$ part of $\overline{\partial}_E+\xi \theta^{\dagger_h}$, and such that $g_{\xi}(0)=1$.
 \item \smash{$N(a):=\bigl[\begin{smallmatrix} n_1(a) & 0 \\
 0 & n_2(a)
 \end{smallmatrix}\bigr]$} with $n_i(a)\in \mathbb{Z}$. This term will ensure that the parabolic weight~associated to $\tau_{i,a}^{\xi}$ lies in $ (a-1,a]$. Here $n_i(a)$ is the unique integer such that $(-1)^{i+1}\bigl[m^{(3)}+2\operatorname{Re}(\xi\overline{m})\bigr] + n_i(a)\in (a-1,a]$ for $i=1,2$.
\end{itemize}

We now use the holomorphic frame $\tau_a^{\xi}$ to perform a holomorphic extension of the bundle~${\mathcal{E}^{\xi}\to \CP \setminus \{\infty\}}$, to a holomorphic bundle over $\CP$. This holomorphic bundle acquires an $a$-parabolic structure at $\infty$, given by the growth conditions of the holomorphic framing \smash{$\tau_a^{\xi}$}. It is easy to check that this bundle is precisely $\mathcal{P}_a^h\mathcal{E}^{\xi}\to \CP$ from Definition \ref{phdef}. Furthermore, we have the following:

\begin{Proposition} In the frame $\tau_a^{\xi}$ the connection $\nabla^{\xi}$ acquires the following form:
\begin{equation}\label{singform}
 \nabla^{\xi}={\rm d}-\bigl(\xi^{-1}+\overline{\xi}\bigr)H\frac{{\rm d}w}{w^3}+\Lambda(\xi)\frac{{\rm d}w}{w}+ \text{holomorphic $(1,0)$ terms},
\end{equation}
where
\begin{equation}\label{expformon}
 \Lambda(\xi)=\begin{bmatrix}-\xi^{-1}m +m^{(3)}+\xi \overline{m} +n_1(a) & 0\\
 0& \xi^{-1}m -m^{(3)}-\xi \overline{m} +n_2(a)
 \end{bmatrix}.
\end{equation}
In particular, $\bigl(\mathcal{P}^h_a\mathcal{E}^{\xi},\nabla^{\xi}\bigr)\to \bigl(\CP,\infty\bigr)$ is an $a$-parabolic flat bundle, and $\tau_a^{\xi}|_{\infty}$ is a compatible frame in the sense of the previous section.

\begin{proof}
Basically the same argument as \cite[Sections 7 and 8]{BB04}. The idea is that the estimates of the entries of $g_{\xi}$ from \cite{BB04}, and the fact that the terms of the connection form $A$ of $\nabla^{\xi}$ in the frame~$\tau_{a,\xi}$ must satisfy $\overline{\partial}A=0$, force the singular part of $\nabla^{\xi}$ to be the one written above.
\end{proof}

\end{Proposition}
From the compatibly framed $a$-parabolic flat bundle \smash{$\bigl(\mathcal{P}^h_a\mathcal{E}^{\xi},\nabla^{\xi},\tau_a^{\xi}\bigr)$}, we get the filtered bundle~${\bigl(\mathcal{P}^h_*\mathcal{E}^{\xi},\nabla^{\xi},\tau_{*}^{\xi}\bigr)}$, where as $\bigl(\mathcal{P}^h_*\mathcal{E}^{\xi},\nabla^{\xi}\bigr)$ is obtained as explained in Section \ref{filtbunsec}, while $\tau_a^{\xi}$ for each $a\in \mathbb{R}$ is determined by \eqref{twistcompframe}.

Notice that while the exponent of formal monodromy $\Lambda(\xi)$ of $\nabla^{\xi}$ depends on $a\in \mathbb{R}$, the formal monodromy does not
\begin{equation*}
 \exp (2\pi {\rm i} \Lambda)=\begin{bmatrix} \exp \bigl(2 \pi {\rm i}\bigl(-\xi^{-1}m +m^{(3)}+\xi \overline{m}\bigr)\bigr) & 0\\
 0 & \exp \bigl(2\pi {\rm i}\bigl(\xi^{-1}m -m^{(3)}-\xi \overline{m}\bigr)\bigr)\\
 \end{bmatrix}.
\end{equation*}

The expression of the diagonal entries greatly resembles the formula of the electric twistor coordinate of the Ooguri--Vafa space (see \eqref{OV electric}). In the next sections, we will see which one of the two diagonal entries is the ``right one'' to pick.
\end{subsubsection}
\begin{subsubsection}{Stokes data of the associated compatibly framed flat bundles}

We will now study the Stokes data of the $\mathbb{C}^*$-family of compatibly framed filtered flat bundles~${\bigl(\mathcal{P}^h_*\mathcal{E}^{\xi},\nabla^{\xi}, \tau_{*}^{\xi}\bigr)}$ for $\xi \in \mathbb{C}^*$. Strictly speaking, so far it only makes sense to associate Stokes data to \smash{$\bigl(\mathcal{P}^h_a\mathcal{E}^{\xi},\nabla^{\xi}, \tau_{a}^{\xi}\bigr)$} for some fixed $a\in \mathbb{R}$. In the last section, we saw that while the exponent of formal monodromy depends on $a$, it does in such a way that the formal monodromy does not depend on $a$. Hence, it makes sense to talk about the formal monodromy associated to~${\bigl(\mathcal{P}^h_*\mathcal{E}^{\xi},\nabla^{\xi}, \tau_{*}^{\xi}\bigr)\to \bigl(\CP,\infty\bigr)}$. The following proposition says that the same holds for the Stokes matrices:

\begin{Proposition}\label{SM} The Stokes matrices associated to $\bigl(\mathcal{P}^h_a\mathcal{E}^{\xi},\nabla^{\xi}, \tau_{a}^{\xi}\bigr)\to \bigl(\CP,\infty\bigr)$ do not depend on the choice of $a\in \mathbb{R}$.
\end{Proposition}

Before proving the proposition, we first prove the following lemmas.
\begin{Lemma}\label{twistasymp}
 Let $\bigl(E,\overline{\partial}_E,\theta,h,(e_1,e_2)|_{\infty}\bigr)\in \mathcal{H}^{{\rm fr}}$ and let $\bigl(\mathcal{P}^h_a\mathcal{E}^{\xi},\nabla^{\xi}, \tau_{a}^{\xi}\bigr)\to \bigl(\CP,\infty\bigr)$ be as before. Consider the holomorphic coordinate $w$ related to $z$ by $w=\frac{1}{z}$, and fix a branch of the logarithm in the $w$-plane. Using this fixed branch of the logarithm, let $\Phi_i(w,\xi)$ denote the frame of flat sections of $\nabla^{\xi}$ defined on the extended sector \smash{$\widehat{\textnormal{Sect}}_i$}, as in Section {\rm\ref{GMD}}. Furthermore, consider the matrix function ${\rm e}^{Q(w,\xi)}$, where $Q$ is a $2\times 2$ diagonal matrix with entries
\begin{gather*}
 Q_1(w,\xi):= -\xi^{-1}\left(\frac{1}{2w^2}-m\operatorname{Log}(w)\right)-{\rm i}m^{(3)}\operatorname{Arg}(w)-\xi\left(\frac{1}{2\overline{w}^2}-\overline{m}\overline{\operatorname{Log}(w)}\right),\\
 Q_2(w,\xi):=\xi^{-1}\left(\frac{1}{2w^2}-m\operatorname{Log}(w)\right)+{\rm i}m^{(3)}\operatorname{Arg}(w)+\xi\left(\frac{1}{2\overline{w}^2}
 -\overline{m}\overline{\operatorname{Log}(w)}\right).
\end{gather*}
 If $\Phi_i(w,\xi)=(e_1,e_2)\cdot A_i(w,\xi)$, then the matrix $A_i(w,\xi)$ satisfies
\smash{$A_i(w,\xi)\cdot {\rm e}^{-Q(w,\xi)} \to 1$} as $ w\to 0$, \smash{$w\in \widehat{\operatorname{Sect}}_i$},
where $Q(w,\xi)$ uses the same branch of the $\operatorname{Log}$ and $\operatorname{Arg}$ as $\Phi_i$.
\end{Lemma}

\begin{proof}
The formal type of the connection $\nabla^{\xi}$ is given is this case by
$
 -\bigl(\xi^{-1}+\overline{\xi}\bigr)H\frac{{\rm d}w}{w^3}+\Lambda\frac{{\rm d}w}{w}$,
where
\begin{equation*}
 \Lambda=\begin{bmatrix}-\xi^{-1}m +m^{(3)}+\xi \overline{m} +n_1(a) & 0\\
 0& \xi^{-1}m -m^{(3)}-\xi \overline{m} +n_2(a)
 \end{bmatrix}.
\end{equation*}
By following the recipe from Section \ref{GMD}, we find that
\begin{equation*}
 \Phi_i(w,\xi)=\tau_{a}^{\xi}\cdot\Sigma_i\bigl(\widehat{F}(\xi)\bigr) w^{-\Lambda(\xi)}\exp \left(-\bigl(\xi^{-1}+\overline{\xi}\bigr)\frac{H}{2w^2}\right).
\end{equation*}
On the other hand, we have that
\begin{equation*}
 \tau_{a}^{\xi}(w)=(e_1,e_2)\cdot g_{\xi}(w) |w|^{(m^{(3)}+2\xi \overline{m})H}w^{N(a)}\exp \left(\frac{\overline{\xi}H}{2w^2}-\frac{\xi H}{2\overline{w}^2}\right),
\end{equation*}
so that
\begin{align*}
 A_i(w,\xi)={}& g_{\xi}(w) |w|^{(m^{(3)}+2\xi \overline{m})H}w^{N(a)}\exp \left(\frac{\overline{\xi}H}{2w^2}-\frac{\xi H}{2\overline{w}^2}\right)\Sigma_i\bigl(\widehat{F}\bigr) w^{-\Lambda(\xi)}\\
 &\times\exp \left(-\bigl(\xi^{-1}+\overline{\xi}\bigr)\frac{H}{2w^2}\right)\\
:={}&g_{\xi}(w)B_i(w,\xi).
\end{align*}

Since $g_{\xi}(w)\to 1$ as $w \to 0$, it is enough to prove that $B_i(w,\xi){\rm e}^{-Q(w,\xi)}\to 1$ as $w\to 0$, \smash{$w\in \widehat{\operatorname{Sect}}_i$}. A simple computation shows that the diagonal entries of $B_i$ are the following:
\begin{equation*}
 B_{i}(w,\xi)_{jj}=\Sigma_{i}\bigl(\widehat{F}\bigr)_{jj}{\rm e}^{Q_j(w,\xi)},
\end{equation*}
while the off-diagonal terms have the form
\begin{equation*}
 B_{i}(w,\xi)_{jk}=\Sigma_{i}\bigl(\widehat{F}\bigr)_{jk}(w^{n(a,j)-n(a,k)}|w|^{(2(-m^{(3)}-2\xi \overline{m})H_{kk})}\exp \left(\left(-\frac{\overline{\xi}}{w^2}+\frac{\xi}{\overline{w}^2}\right)H_{kk}\right){\rm e}^{Q_k(w,\xi)}.
\end{equation*}

On the other hand, by the proof of \cite[Lemma 1]{B2}, we have the following estimates:
\[
\Sigma_i\bigl(\hat{F}\bigr)_{jj}-1=\mathcal{O}(|w|) \qquad \text{and}\qquad \Sigma_i\bigl(\hat{F}\bigr)_{jk}=\mathcal{O}(|w|) \quad \text{for $j\neq k$}.
\]
 We also have the following:
 \[
 \bigl|n(a,j)-n(a,k)+\bigl(2\bigl(-m^{(3)}-2\operatorname{Re}(\xi \overline{m})\bigr)H_{kk}\bigr)\bigr|<1,
 \]
 since that quantity is the difference of the two parabolic weights in the range $(a-1,a]$. Hence find that
\begin{equation*}
 B_i(w,\xi)_{jj}{\rm e}^{-Q_j(w,\xi)} \to 1, \qquad
 B_i(w,\xi)_{jk}{\rm e}^{-Q_k(w,\xi)} \to 0,
\end{equation*}
so we get what we wanted.
\end{proof}

\begin{Lemma}\label{UniqueA} The asymptotics of the Lemma {\rm\ref{twistasymp}} uniquely characterize the frame of flat sections $\Phi_i(w,\xi)$.
\end{Lemma}
\begin{proof}Suppose that we have two frames of flat sections $\Phi_i(w,\xi)=(e_1,e_2)\cdot A(w,\xi)$ and $\widetilde{\Phi}_i(w,\xi)=(e_1,e_2)\cdot \widetilde{A}(w,\xi)$ defined on \smash{$\widehat{\operatorname{Sect}}_i$}, and satisfying the asymptotics of the previous lemma.

We then have that $A(w,\xi)=\widetilde{A}(w,\xi)S$ for some constant matrix $S$, since $S$ is the transition function between flat sections. From the asymptotic conditions, we then have that
\begin{gather*}
 A(w,\xi){\rm e}^{-Q(w,\xi)}=\widetilde{A}{\rm e}^{-Q(w,\xi)}{\rm e}^{Q(w,\xi)}S{\rm e}^{-Q(w,\xi)} \to 1 \qquad \text{as} \quad w \to 0,\quad w \in \widehat{\operatorname{Sect}}_i,\\
 \widetilde{A}{\rm e}^{-Q(w,\xi)}\to 1\qquad \text{as} \quad w \to 0,\quad w \in \widehat{\operatorname{Sect}}_i,
\end{gather*}
so we conclude that
\begin{equation*}
 {\rm e}^{Q(w,\xi)}S{\rm e}^{-Q(w,\xi)}=\begin{bmatrix} s_{11} & {\rm e}^{Q_1-Q_2}s_{12}\\
 {\rm e}^{Q_2-Q_1}s_{21} & s_{22}\\
 \end{bmatrix}\to 1 \qquad \text{as} \quad w \to 0,\quad w \in \widehat{\operatorname{Sect}}_i.
\end{equation*}

On the sector $\widehat{\operatorname{Sect}}_i$, there are two subsectors separated by a Stokes ray. This in particular implies that if in one of the subsectors $\operatorname{Re}(Q_i-Q_j)>0$, then on the other subsector we have~${\operatorname{Re}(Q_i-Q_j)<0}$. Hence, the limit above forces $s_{21}=s_{12}=0$ and $s_{11}=s_{22}=1$. It then follows that $\Phi_i=\widetilde{\Phi}_i$.
\end{proof}

\begin{proof}[Proof of Proposition \ref{SM}] By the previous two lemmas, we see that the frames of flat sections~$\Phi_i$ satisfy the same asymptotics of the first lemma independently of the choice of $a\in \mathbb{R}$; and these asymptotics uniquely characterize them. Hence, the frames $\Phi_i$ do not depend on~${a \in \mathbb{R}}$, so we conclude that the Stokes matrices also do not depend on $a\in \mathbb{R}$.
\end{proof}
\end{subsubsection}
\begin{subsubsection} {Definition of the twistor coordinates}\label{deftwistcoords}

In the following, we will assume the fact that the Stokes data associated to $\bigl(\mathcal{P}^h_*\mathcal{E}^{\xi},\nabla^{\xi}, \tau_{*}^{\xi}\bigr)\to \bigl(\CP,\infty\bigr)$ depends holomorphically on the parameter $\xi \in \mathbb{C}^*$. This will be proved later in the next section. We will also assume throughout the section that $m\in \mathbb{C}^{*}$, where $m$ is the parameter appearing in \eqref{singHiggs}.

We will study how the Stokes matrices and formal monodromy associated to $\bigl(\mathcal{P}^h_*\mathcal{E}^{\xi},\nabla^{\xi}, \tau_{*}^{\xi}\bigr)\to \bigl(\CP,\infty\bigr)$ varies when we vary the parameter $\xi \in \mathbb{C}^*$. The behavior of the Stokes matrices with respect to the twistor parameter will motivate the definition of the twistor coordinates for framed wild harmonic bundles.

Let us recall the setting of Section \ref{GMD} in the specific case we are working with. We now have a compatibly framed rank 2 flat vector bundle $\bigl(\mathcal{P}^h_a\mathcal{E}^{\xi},\nabla^{\xi}, \tau_{a}^{\xi}\bigr)\to \bigl(\CP,\infty\bigr)$. The compatible frame $\tau_{a,\xi}$ extends to a local holomorphic frame near $z=\infty$, where $\nabla^{\xi}$ has the following form (with respect to the holomorphic coordinate $w=\frac1z$):
\begin{equation*}
 \nabla^{\xi}={\rm d} -\bigl(\xi^{-1}+\overline{\xi}\bigr)H\frac{{\rm d}w}{w^3}+\Lambda(\xi)\frac{{\rm d}w}{w} + \text{holomorphic $(1,0)$ terms},
\end{equation*}
where
\begin{equation}\label{formmon}
 \Lambda(\xi)=\begin{bmatrix}-\xi^{-1}m +m^{(3)}+\xi \overline{m} +n_1(a) & 0\\
 0& \xi^{-1}m -m^{(3)}-\xi \overline{m} +n_2(a)
 \end{bmatrix}.
\end{equation}

The formal type of the connection is given in our case by
\begin{equation}\label{formtype}
 -\bigl(\xi^{-1}+\overline{\xi}\bigr)H\frac{{\rm d}w}{w^3}+\Lambda(\xi)\frac{{\rm d}w}{w},
\end{equation}
so that in the notation of Section \ref{GMD} we have that $-Q(w,\xi)=-\bigl(\xi^{-1}+\overline{\xi}\bigr)\frac{H}{2w^2}=\operatorname{diag}(q_1,q_2)$ determines the Stokes and anti-Stokes rays. In this case, we have $2$ anti-Stokes rays where $q_{12}=-\bigl(\xi^{-1}+\overline{\xi}\bigr)\frac{1}{w^2}<0$ and $2$ anti-Stokes rays where $q_{21}=\bigl(\xi^{-1}+\overline{\xi}\bigr)\frac{1}{w^2}<0$. Similarly, we have~$4$ Stokes rays corresponding to $\operatorname{Re}(q_{12})=0$ (see Figure~\ref{fig7} from Section~\ref{GMD}).

For the definition of Stokes data, we need to fix a branch of the Log and a labeling of the sectors determined by anti-Stokes rays. The reason for the choices made below will become apparent in subsequent sections, where we start comparing and matching the twistor coordinates for harmonic bundles with the twistor coordinates of the Ooguri--Vafa space. The choices of branch and labeling will depend on the value of $m \in \mathbb{C}^*$ in the following way:
\begin{itemize}\itemsep=0pt
 \item Let $\operatorname{Arg}(w)$ denote the argument function with values in $[0,2\pi)$. Given $m\in \mathbb{C}^*$, we denote by $\operatorname{Arg}_m(w)$ the argument function with values in $[-\frac 12 \operatorname{Arg}(m), -\frac12 \operatorname{Arg}(m)+2\pi)$. Furthermore, we denote by $\operatorname{Log}_m(w)$ the branch of the logarithm that uses $\operatorname{Arg}_m(w)$.

 \item For $\xi=m$, the ray from $0$ to \smash{${\rm e}^{-\frac{1}{2}{\rm i}\operatorname{Arg}(m)}$} is one of the anti-Stokes lines. We will denote this anti-Stokes line by $r_{\theta_1}(\xi=m)$ and denote the others by $r_{\theta_i}(\xi=m)$ for $i=2,3,4$ in a counterclockwise manner. As before, we let $\operatorname{Sect}_i(\xi=m)$ be the sector going from $r_{\theta_i}$ to~$r_{\theta_{i+1}}$ and denote by \smash{$\widehat{\operatorname{Sect}}_i(\xi=m)$} the extended sectors (see Figure~\ref{fig2}).

 \begin{figure}[ht] \centering
 \includegraphics[width=13cm]{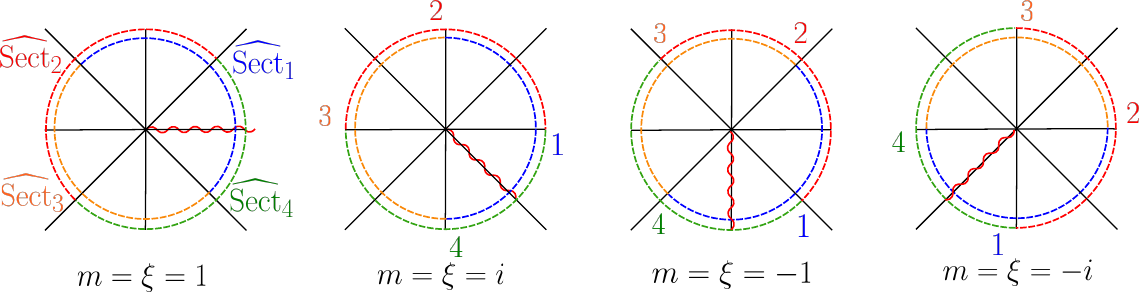}
 \caption{How the labeling of the sectors varies with $m$ while keeping $m=\xi$. The wavy red line denotes the place of the branch cut of $\operatorname{Arg}_m$.}
 \label{fig2}
 \end{figure}

 \item When $\xi$ varies from $\xi=m$, the anti-Stokes rays move continuously. Given a fixed $m \in \mathbb{C}^*$, we will denote by $\operatorname{Sect}_i(\xi)$ the sector obtained by varying $\xi$ starting from $\xi=m$ and with~${\xi \in \mathbb{C}^*\setminus \{\xi \in \mathbb{C}^* \mid -\xi^{-1}{\rm i}m>0 \}}$. We will denote the corresponding extended sectors by \smash{$\widehat{\operatorname{Sect}}_i(\xi)$} (see Figure~\ref{fig4}).

 \begin{figure}[ht] \centering
 \includegraphics[width=13cm]{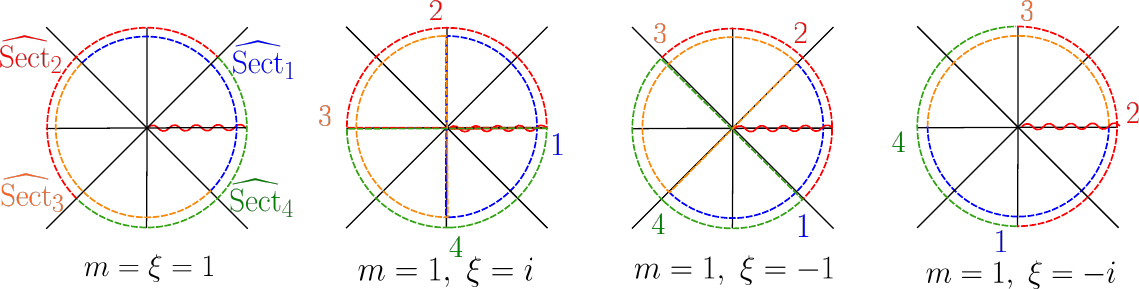}
 \caption{How the labeling of the sectors change with $m=1$ fixed, while we vary $\xi$ from $\xi=1$ in a~counterclockwise manner. The wavy red line denotes the place of the branch cut of $\operatorname{Arg}_m$.}
 \label{fig4}
 \end{figure}
\end{itemize}

Now we set our conventions for the sectorial frames of flat sections $\Phi_i(\xi)$ of $\nabla^{\xi}$ used to define the Stokes matrices:
\begin{itemize}\itemsep=0pt
 \item To define the sectorial flat frames $\Phi_i(\xi)$, we must first fix a branch of the logarithm. For a~fixed $m\in \mathbb{C}^*$ and $\xi=m$, we will use $\operatorname{Log}_m$ and the convention explained in the paragraph before Definition \ref{defStokesmatrices} (i.e., it uses $\operatorname{Log}_m(w)$ if $i\neq 1,4$; it uses the analytically continued branch of $\operatorname{Log}_m(w)$ from $\operatorname{Sect}_1(\xi)$ to $\operatorname{Sect}_{4}(\xi)$ if $i=1$; and it uses the analytically continued branch of $\operatorname{Log}_m(w)$ from $\operatorname{Sect}_{4}$ to $\operatorname{Sect}_{1}$ for $i=4$). The branch that we use to define $\Phi_i(\xi)$ for $\xi\neq m$ varies continuously with the sector \smash{$\widehat{\operatorname{Sect}_i}(\xi)$}. For example, for $m=1$, $\Phi_1(\xi)$ uses the branch with argument taking values in $\bigl[-\frac12\operatorname{Arg}(\xi) -\frac{\pi}{4}, -\frac 12 \operatorname{Arg}(\xi) +\frac{7\pi}{4}\bigr)$ (see Figure~\ref{fig4}).
 \item With the above choice of branch, $\Phi_i(\xi)$ denotes the flat frame on \smash{$\widehat{\operatorname{Sect}_i}(\xi)$} specified by the asymptotics from Lemma \ref{twistasymp} (which uniquely characterizes the $\Phi_i(\xi)$ by Lemma $\ref{UniqueA}$).
 \item We define the Stokes matrices $S_i(\xi)$ using the $\Phi_i(\xi)$ as in Definition \ref{defStokesmatrices}.
\end{itemize}

 Having defined the required Stokes matrices, recall that by Proposition \ref{MonodromyStokes}, we have the following relation for the Stokes data:
\begin{equation}\label{trivialmonodromy}
 S_1(\xi)S_2(\xi)S_3(\xi)S_4(\xi)M_0^{-1}(\xi)=1,
\end{equation}
where $M_0(\xi)={\rm e}^{-2\pi {\rm i} \Lambda(\xi)}$.

We denote by $a$, $b$, $c$ and $d$, the nontrivial off-diagonal elements of the $2\times 2$ unipotent matrices~$S_1$, $S_2$, $S_3$ and $S_4$, respectively, and $M_0=\operatorname{diag}\bigl(\mu,\mu^{-1}\bigr)$. With our conventions $S_1$ turns out to be lower triangular (and hence so is $S_3$, while $S_2$ and $S_4$ are upper triangular), and we get the following relations among the Stokes matrices elements from the matrix relation of equation~\eqref{trivialmonodromy}:
\begin{equation}\label{Stokes relations}
 1+bc =\mu,\qquad
 \mu d+ b =0,\qquad
 a+\mu^{-1} c =0,\qquad
 ab+1 =\mu^{-1}.
\end{equation}

For now, we will consider $a(\xi)$, $b(\xi)$, $c(\xi)$, $d(\xi)$ as analytic functions on the half plane $\mathbb{H}_m$, where $\mathbb{H}_m:=\bigl\{ \xi \in \mathbb{C}^* \mid \operatorname{Re}\bigl(\xi^{-1}m\bigr)>0\bigr\}$ (recall the assumption made at the beginning of the section).

If we denote by $s_i$ the flat section of the frame $\Phi_i$ that is exponentially decreasing on \smash{$\widehat{\operatorname{Sect}}_{i-1}\cap \widehat{\operatorname{Sect}}_i$} (with $i$ taken mod $4$), we have that $s_i=\Phi_i\cdot (1,0)^{t}$ for $i=1,3$ and $s_i=\Phi_i \cdot (0,1)^{t}$ for $i=2,4$. They satisfy the following relations among themselves, obtained by using the Stokes matrices and how they relate the flat frames $\Phi_i$ for $i=1,2,3,4$,
\begin{gather}
 s_3=s_1+as_2 \qquad \text{on} \quad \widehat{\operatorname{Sect}_1}\cap\widehat{\operatorname{Sect}_2},\nonumber\\
 s_4=s_2+bs_3 \qquad \text{on} \quad \widehat{\operatorname{Sect}_2}\cap\widehat{\operatorname{Sect}_3},\nonumber\\
 \mu s_1= s_3 +cs_4 \qquad\text{on} \quad \widehat{\operatorname{Sect}_3}\cap\widehat{\operatorname{Sect}_4},\nonumber\\
 \mu^{-1}s_2=s_4+d \mu s_1 \qquad \text{on} \quad \widehat{\operatorname{Sect}_4}\cap\widehat{\operatorname{Sect}_1}.\label{Stokes relations 3}
\end{gather}

Now consider the rays in the $\xi$ plane determined by
\begin{equation*}
 l_{\pm}(-2{\rm i}m)=\left\{\xi\in \mathbb{C}^* \,\Big|\, \pm \frac{ -2{\rm i}m}{\xi}<0\right\}.
\end{equation*}
Note that these rays, together with $0$, form the boundary of the half planes $\mathbb{H}_{m}$ and $\mathbb{H}_{-m}$ (while we could write $l_{\pm}(-{\rm i}m)$ instead of $l_{\pm}(-2{\rm i}m)$, the reason we write $-2{\rm i}m$ instead of just $-{\rm i}m$ is because $-2{\rm i}m$ will be identified with the base parameter of the Ooguri--Vafa space). If we start at $\mathbb{H}_m$ and move counterclockwise, we cross $l_{+}(-2{\rm i}m)$ into $\mathbb{H}_{-m}$, and then we cross $l_{-}(-2{\rm i}m)$ into $\mathbb{H}_{m}$ (see Figure~\ref{fig6}).

\begin{figure}
 \centering
 \includegraphics[width=7cm]{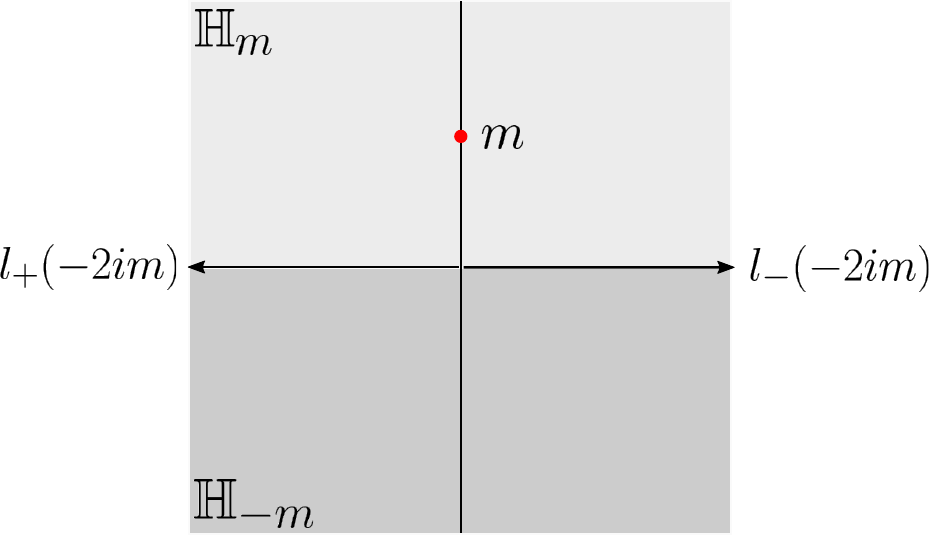}
 \caption{Given $m \in \mathbb{C}^*$, the figure above shows the configuration of the half-planes $\mathbb{H}_{\pm m}$ and the rays~$l_{\pm}(-2{\rm i}m)$.}
 \label{fig6}
\end{figure}

We can analytically continue the Stokes matrices elements from $\mathbb{H}_m$ to $\mathbb{H}_{-m}$ through \linebreak $l_+({-}2{\rm i}m)$, while still maintaining the same relations $\eqref{Stokes relations}$ as before.\footnote{From the arguments in the next section, where we show the holomorphic dependence of the flat frames $\Phi_i(\xi)$ in $\xi$; it will be clear that given any $\xi_0\in \mathbb{C}^*$, there is a small sector centered at the ray from $0$ to $\xi_0$ where $\Phi_i(\xi)$ is holomorphic for $\xi$ in the sector. The opening of this sector will be uniform in $\xi_0\in \mathbb{C}^*$, so we can perform the analytic continuations of Stokes data mentioned above.} The Stokes matrix elements are then holomorphic functions on $\mathbb{C}^*-l_-(-2{\rm i}m)$. Now suppose that we further analytically continue from~$\mathbb{H}_{-m}$ to $\mathbb{H}_m$ by going through $l_-(-2{\rm i}m)$. Notice that in the process of doing this, the labeling of our sectors experience monodromy: \smash{$\widehat{\operatorname{Sect}}_1$} gets interchanged with \smash{$\widehat{\operatorname{Sect}}_3$}, and~\smash{$\widehat{\operatorname{Sect}}_2$} gets interchanged with \smash{$\widehat{\operatorname{Sect}}_4$} (see Figure~\ref{fig4}). We will denote the Stokes matrix elements on~${\mathbb{C}^*-l_+(-2{\rm i}m)}$ obtained this way by $\widehat{a}$, $\widehat{b}$, $\widehat{c}$, $\widehat{d}$. On $\mathbb{H}_{-m}$ they coincide with the previous ones. However, if we denote by $\widehat{\Phi}_i$ the corresponding fundamental solutions on the corresponding sector, we have the following relations for $\xi \in \mathbb{H}_{m}$
\begin{equation}\label{Stokes relations 2}
 \widehat{\Phi}_1=\Phi_3\cdot M_0^{-1},\qquad
 \widehat{\Phi}_2=\Phi_4\cdot M_0^{-1},\qquad
 \widehat{\Phi}_3=\Phi_1,\qquad
 \widehat{\Phi}_4=\Phi_2.
\end{equation}

Now we finally define the electric and magnetic twistor coordinates. We will verify that the magnetic coordinate has the appropriate jumps on $l_{\pm}(-2{\rm i}m)$, matching with the corresponding jumps of the Ooguri--Vafa coordinates (see Proposition \ref{jumpsOV}). For the definition of the magnetic coordinate, we will need to assume for now that the Stokes matrix element given by $b(\xi)$ does not vanish for $\xi \in \mathbb{H}_{-m}$ (this will be shown in Section \ref{nonvanishing}).

\begin{Definition} Let $\bigl(E,\overline{\partial}_E,\theta,h,g\bigr)\in \mathcal{H}^{{\rm fr}}$ with parameters $m$ and \smash{$m^{(3)}$} describing the corresponding singularities. For $\xi \in \mathbb{C}^*$, the electric twistor coordinate is then defined by
\begin{equation}\label{defetc}
\mathcal{X}_e\bigl(\bigl(E,\overline{\partial}_E,\theta,h,g\bigr),\xi\bigr):=\mu^{-1}(\xi)=\exp\bigl(-2\pi {\rm i} \bigl(\xi^{-1}m -m^{(3)}-\overline{m}\xi\bigr)\bigr).
\end{equation}
\end{Definition}

\begin{Definition} Let $\bigl(E,\overline{\partial_E},\theta,h,g\bigr)\in \mathcal{H}^{{\rm fr}}$ with parameters $m\neq 0$ and \smash{$m^{(3)}$} describing the corresponding singularities. The magnetic twistor coordinate is defined by
\begin{equation*}
 \mathcal{X}_m\bigl(\bigl(E,\overline{\partial_E},\theta,h,g\bigr),\xi\bigr):=\begin{cases}
 a(\xi) & \text{for} \ \xi \in \mathbb{H}_{m}, \\
 -\dfrac{1}{b(\xi)} & \text{for} \ \xi \in \mathbb{H}_{-m},
 \end{cases}
\end{equation*}
where $a(\xi)$ and $b(\xi)$ are the corresponding Stokes matrix elements of the associated compatibly framed filtered flat bundle \smash{$\bigl(\mathcal{P}^h_*\mathcal{E}^{\xi},\nabla^{\xi}, \tau_{*}^{\xi}\bigr)\to \bigl(\CP,\infty\bigr)$} $($holomorphic on $\mathbb{C}^*\setminus l_{-}(-2{\rm i}m))$.
\end{Definition}

\begin{Remark} Notice that the subscript $m$ of $\mathcal{X}_m$ refers to ``magnetic'', and not to the complex parameter $m$. We hope this notation does not cause confusion.
\end{Remark}
\begin{Lemma} $\mathcal{X}_e\bigl(\bigl(E,\overline{\partial_E},\theta,h,g\bigr),\xi\bigr)$ and $\MC\bigl(\bigl(E,\overline{\partial_E},\theta,h,g\bigr),\xi\bigr)$ descend to functions on $\mathfrak{X}^{{\rm fr}}$.

\end{Lemma}

\begin{proof}
If we have an isomorphism $\varphi\colon \bigl(E_1,\overline{\partial}_{E_1},\theta_1,h_1,g_1\bigr)\to \bigl(E_2,\overline{\partial}_{E_2},\theta_2,h_2,g_2\bigr)$, it is easy to check that it induces an isomorphism between the associated compatibly framed filtered flat bundles $\bigl(\mathcal{P}_{*}^{h_i}\mathcal{E}^{\xi}_i,\nabla^{\xi}_i, \tau_{*,i}^{\xi}\bigr)\to \bigl(\CP,\infty\bigr)$ for $i=1,2$. It then follows that both compatibly framed filtered flat bundles have the same Stokes matrices and formal monodromy. Hence, the twistor coordinates descend to $\mathfrak{X}^{{\rm fr}}$.
\end{proof}

We will often just write $\mathcal{X}_e(\xi)$ and $\MC(\xi)$, with the understanding that they depend on elements of $\mathfrak{X}^{{\rm fr}}$.

We now prove the following analog of Proposition \ref{jumpsOV}.

\begin{Theorem} For a fixed $\bigl(E,\overline{\partial}_E,\theta,h,g\bigr) \in \mathcal{H}^{{\rm fr}}$ with $m\neq 0$, the magnetic twistor coordinate~$\mathcal{X}_m(\xi)$ has the following jumps when we vary the twistor parameter $\xi$
\begin{gather*}
 \mathcal{X}_m(\xi)^{+}=\mathcal{X}_m(\xi)^{-}(1-\mathcal{X}_e(\xi))^{-1} \qquad \text{along} \quad \xi \in l_{+}(-2{\rm i}m),\\
 \mathcal{X}_m(\xi)^{+}=\mathcal{X}_m(\xi)^{-}\bigl(1-\mathcal{X}_e(\xi)^{-1}\bigr) \qquad \text{along} \quad \xi \in l_{-}(-2{\rm i}m),
\end{gather*}
where the $+$ or $-$ on the coordinate denotes the clockwise or counterclockwise limit to the ray, respectively.
\end{Theorem}

\begin{proof}
Because of the relations \eqref{Stokes relations} and the analytic continuation that we chose, we automatically have the relation
\smash{$
 a\bigl(1-\mu^{-1}\bigr)^{-1}=-\frac{1}{b}$} along $ \xi \in l_+(-2{\rm i}m)$.

Hence, we have the following relation:
\begin{equation*}
 \mathcal{X}_m(\xi)^{+}=\mathcal{X}_m(\xi)^{-}\bigl(1-\mu^{-1}\bigr)^{-1}=\mathcal{X}_m(\xi)^{-}(1-\mathcal{X}_e(\xi))^{-1} \qquad \text{along} \quad \xi \in l_{+}(-2{\rm i}m).
\end{equation*}

To check the other jump, notice that by \eqref{Stokes relations 2} we have that on $\mathbb{H}_m$: $s_1=\widehat{s}_3$, $s_2=\widehat{s}_4$, $s_3=\mu \widehat{s}_1$, and $s_4=\mu^{-1} \widehat{s}_2$, so that along $l_-(-2{\rm i}m)$
\begin{equation*}
 a^{+}=\left(\frac{s_3\wedge s_1}{s_2\wedge s_1}\right)^{+}
 = \left(\frac{\mu \widehat{s_1}\wedge \widehat{s_3}}{\widehat{s_4}\wedge \widehat{s_3}}\right)^{+}=\left(\frac{\mu \widehat{s_1}\wedge \widehat{s_3}}{\widehat{s_4}\wedge \widehat{s_3}}\right)^{-}=\left(\mu \frac{s_1\wedge s_3}{s_4\wedge s_3}\right)^{-}=c^-.
\end{equation*}

On the other hand, from \eqref{Stokes relations} we see that
$
 c=-\frac{1}{b}(1-\mu)$,
so that
$
 a^{+}=c^-=-\frac{1}{b^-}(1-\mu)$ along $\xi \in l_-(-2{\rm i}m)$.

Hence, we get the following jump:
\begin{equation*}
 \mathcal{X}_m(\xi)^{+}=\mathcal{X}_m(\xi)^{-}(1-\mu)=\mathcal{X}_m(\xi)^{-}\bigl(1-\mathcal{X}_e(\xi)^{-1}\bigr) \qquad \text{along} \quad \xi \in l_{-}(-2{\rm i}m).
\end{equation*}
This shows what we wanted.
\end{proof}

We then conclude that the magnetic coordinate has the same jumps as the magnetic coordinate of the Ooguri--Vafa space.

\end{subsubsection}

\subsection[Twistor coordinates Part 2: Holomorphic dependence in the twistor parameter]{Twistor coordinates Part 2:\\ Holomorphic dependence in the twistor parameter}

In this section, we prove the holomorphic dependence of the Stokes data with respect to the twistor parameter $\xi \in \mathbb{C}^*$. This will show that $\mathcal{X}_{m}(\xi)$ depends holomorphically on $\xi \in \mathbb{C}^*-l_{\pm}(-2{\rm i}m)$. We will use the procedure of deformation of irregular values, which varies the compatibly framed flat bundle, while keeping the Stokes data the same. The idea of ``deformation of irregular values'' goes back to Jimbo--Miwa--Ueno \cite{JMU}. Some other references can be found in~\cite{B} and \cite{M6}; while a reference where this is applied to the twistor family of meromorphic bundles associated to a wild harmonic bundle can be found in \cite[Chapters 3, 4, 9 and 11]{M}.

\begin{subsubsection}{Deformation of irregular values}\label{defirr}

One of the main challenges in trying to show that Stokes data depends holomorphically on the twistor parameter, is the fact that we are looking at Stokes data associated to meromorphic flat bundles \smash{$\bigl(\mathcal{P}^h_a\mathcal{E}^{\xi},\nabla^{\xi},\tau^{\xi}_{a}\bigr)$}, with formal type given by
\begin{equation}\label{formaltype}
 -\frac{1+|\xi|^2}{\xi}\frac{H{\rm d}w}{ w^3}+\Lambda(\xi)\frac{{\rm d}w}{w}.
\end{equation}
The cubic pole term has non-holomorphic dependence in $\xi$, so it is not clear that Stokes data is going to depend holomorphically on $\xi$. One first step to solve this is the following construction of ``deformation of irregular values''. We will follow the reference \cite[Section 4.5]{M}.\footnote{The procedure of deformation of irregular values explained in \cite{M} is actually more general than what we will explain here. In our particular case, the leading term of equation \eqref{formaltype} fails to be holomorphic in $\xi$ by the real factor $1+|\xi|^2$, so our deformation will not change the Stokes sectors.}

This procedure of deformation of irregular values will take the data of \smash{$\bigl(\mathcal{P}^h_a\mathcal{E}^{\xi},\nabla^{\xi},\tau^{\xi}_{a}\bigr)$}, and produce a new compatibly framed flat bundle \smash{$\bigl(\mathcal{Q}_a\mathcal{E}^{\xi},\nabla^{\xi},\nu^{\xi}_{a}\bigr)$} such that
\begin{itemize}\itemsep=0pt
 \item The new framing at $\infty$ extends to a holomorphic framing where $\nabla^{\xi}$ has the form
 \begin{equation*}
 \nabla^{\xi}={\rm d}-\frac{H{\rm d}w}{ \xi w^3}+\Lambda(\xi)\frac{{\rm d}w}{w}+\text{holomorphic $(1,0)$ terms}.
 \end{equation*}
 \item The Stokes data associated to \smash{$\bigl(\mathcal{Q}_a\mathcal{E}^{\xi},\nabla^{\xi},\nu^{\xi}_{a}\bigr)$} is the same as the Stokes data of $\smash{\bigl(\mathcal{P}^h_a\mathcal{E}^{\xi},\nabla^{\xi}},\allowbreak\smash{\tau^{\xi}_{a}}\bigr)$.
\end{itemize}
We will explain how the construction works in our setting:

To start, let $U_{\infty}$ denote a neighborhood of $\infty$ in $\CP$, and consider the covering of $U_{\infty}$ by the extended Stokes sectors $\mathcal{S}_i:=\widehat{\operatorname{Sect}}_i$ (with $i=1,2,3,4$). On each such sector, we have the holomorphic frames \smash{$\tau^{\xi}_a \cdot \Sigma_{i}\bigl(\widehat{F}(\xi)\bigr)$}, where the connection $\nabla^{\xi}$ has the form
\begin{equation}\label{canonical form}
 \nabla^{\xi}={\rm d}-\frac{1+|\xi|^2}{\xi}\frac{H{\rm d}w}{ w^3}+\Lambda(\xi)\frac{{\rm d}w}{w}.
\end{equation}

Here we recall that $\widehat{F}(\xi)$ is the formal gauge transformation that takes the connection to the form \eqref{canonical form} and $\widehat{F}(w=0)=1$.

Now let \smash{$\nu_{i}:=\tau^{\xi}_a \cdot \Sigma_{i}\bigl(\widehat{F}\bigr)\cdot \exp \bigl(\frac{-\overline{\xi}}{2w^2}H\bigr)$}. In this new sectorial frame the connection has the form
\begin{equation*}
 \nabla^{\xi}={\rm d}-\frac{H{\rm d}w}{\xi w^3}+\Lambda(\xi)\frac{{\rm d}w}{w}.
\end{equation*}

With these sectorial frames we cannot a priori define an extension of the holomorphic flat bundle $\bigl(\mathcal{E}^{\xi},\nabla^{\xi}\bigr)\to \CP \setminus \{\infty\}$ to a meromorphic flat bundle over $\CP$. However, if~$\smash{\pi\colon\widetilde{\CP}(\infty)}\to \CP$ denotes the real blowup of $\CP$ at $\infty$, we have the following.

\begin{Proposition} Using the sectorial frames $\nu_i$, we can extend the holomorphic flat bundle~${\bigl(\mathcal{E}^{\xi},\nabla^{\xi}\bigr)\to \CP \setminus \{\infty\}}$ to a meromorphic flat bundle \smash{$\bigl(\widetilde{\mathcal{E}}^{\xi},\widetilde{\nabla}^{\xi}\bigr)\to \bigl(\widetilde{\CP}(\infty),\pi^{-1}(\infty)\bigr)$}.

\begin{proof}
Away from $\infty$ the holomorphic sections of the bundle and the connection are the same as before. On the other hand, on the sectorial neighborhoods $\overline{\mathcal{S}_i}$ in the blowup, we have the following transition functions between the frames $\nu_{S_i}$ and $\nu_{S_{i+1}}$:\footnote{This follows from the fact that the $\nu_{S_i}$ differ from the flat frames defining Stokes data by the factor~\smash{$w^{-\Lambda(\xi)}{\rm e}^{-\frac{H}{2\xi w^2}}$}.}
\begin{equation*}
 \begin{bmatrix}
 1 & cw^{\Lambda_1-\Lambda_2}\exp \bigl(-\frac{H}{\xi w^2}\bigr)\\
 0 & 1\\
 \end{bmatrix}
\qquad \text{or}\qquad
 \begin{bmatrix}
 1 & 0\\
 cw^{\Lambda_2-\Lambda_1}\exp \bigl(\frac{H}{\xi w^2}\bigr) & 1\\
 \end{bmatrix},
\end{equation*}
depending on whether $\exp \bigl(-\frac{H}{\xi w^2}\bigr)$ or $\exp \bigl(\frac{H}{\xi w^2}\bigr)$ is exponentially decreasing on $\mathcal{S}_i\cap \mathcal{S}_{i+1}$. Here~``$c$'' denotes the non-trivial Stokes matrix element of the Stokes matrix associated to $\mathcal{S}_i\cap \mathcal{S}_{i+1}$.

The transition functions are holomorphic functions on $\overline{\mathcal{S}_i}\cap \overline{\mathcal{S}_{i+1}}$ (in the sense of being holomorphic in the real blowup, see \cite[Section 3.1.3]{M}) and hence defines a holomorphic bundle over~\smash{$\widetilde{\CP}(\infty)$}. The connection $\nabla^{\xi}$ clearly extends to a meromorphic connection on \smash{$\widetilde{\CP}(\infty)$} with poles along $\pi^{-1}(\infty)$.
\end{proof}

\end{Proposition}

Now we will explain how to actually get a meromorphic flat bundle over $\CP$. Let $\chi_{i}$ be a~partition of unity subordinate to the extended Stokes sectors $\mathcal{S}_i$. We will furthermore pick them such that for any differential operator $D$, we have that $D\chi_{i}=\mathcal{O}\bigl(|w|^{N}\bigr)$ for some natural number~$N$ (depending on the differential operator). We define a $C^{\infty}$ frame on $U_{\infty}\setminus \{\infty\}$ by~$
 \nu:=\sum_i \nu_{i}\chi_{i}$.

The fact that this actually turns out to be a frame follows from the form of the transition functions between the $\nu_{i}$, and the fact that the $\chi_{i}$ are non-negative.

We have the following lemmas concerning the smooth frame $\nu$.

\begin{Lemma}[{\cite[Lemma 3.1.15]{M}}] Let $C_i$ be the matrix such that $\nu=\nu_i\cdot(1+C_i)$. Then if~$Z=\overline{\mathcal{S}}_i\cap \pi^{-1}(\infty)$, we have that $C_i$ goes to $0$ faster than any polynomial as we get near to $Z$.

\end{Lemma}

\begin{Lemma}[{\cite[Lemma 3.1.16]{M}}] If $A$ is such that $\overline{\partial}_{\mathcal{E}^{\xi}}\nu=\nu \cdot A$, then for each sector $\mathcal{S}$ as before and $Z=\overline{\mathcal{S}}\cap \pi^{-1}(\infty)$, we have that $A$ goes to $0$ faster than any polynomial as we get near to $Z$.

\end{Lemma}
In particular, the last lemma implies the following.

\begin{Corollary}\label{exthol}
$A$ descends to give a smooth matrix of functions over a neighborhood of $\infty$ that vanishes faster than any polynomial as $z\to \infty$.
\end{Corollary}

We are now ready to define the deformed compatibly framed flat meromorphic bundle $\bigl(\mathcal{Q}_a\mathcal{E}^{\xi},\nabla^{\xi},\nu^{\xi}_{a}\bigr)$.

\begin{Proposition} Let $\mathcal{Q}_a\mathcal{E}^{\xi}\to \CP$ denote the $C^{\infty}$ bundle defined by extending the bundle~${\mathcal{E}^{\xi} \to \CP \setminus \{\infty\}}$ to a bundle over $\CP$ by using the $\nu$ frame. Then
\begin{itemize}\itemsep=0pt
 \item The holomorphic structure of $\mathcal{E}^{\xi} \to \CP \setminus \{\infty\}$ extends to a holomorphic structure on~$\mathcal{Q}_a\mathcal{E}^{\xi}$.
 \item The connection $\nabla^{\xi}$ becomes a meromorphic connection on $\mathcal{Q}_a\mathcal{E}^{\xi}$, with poles at $\infty$.
 \item The frame $\nu$ used to define the extension gives a compatible framing at $\infty$. The compatible framing extends to a local holomorphic framing where the connection has the form
 \begin{equation}\label{holform}
 \nabla^{\xi}= {\rm d}-\frac{H{\rm d}w}{\xi w^3}+\Lambda(\xi)\frac{{\rm d}w}{w} + \text{holomorphic $(1,0)$ terms}.
 \end{equation}
\end{itemize}
We denote the compatibly framed meromorphic flat bundle that we obtained by \smash{$\bigl(\mathcal{Q}_a\mathcal{E}^{\xi},\nabla^{\xi},\nu^{\xi}_{a}\bigr)$}.
\end{Proposition}

\begin{proof}
The fact that the holomorphic structure of $\mathcal{E}^{\xi}$ extends to a holomorphic structure on~$\mathcal{Q}_a\mathcal{E}^{\xi}$ follows from Corollary \ref{exthol}.

On the other hand, it is easy to check that because of our conditions on the $\chi_i$, we have that in the frame $\nu$
\begin{equation*}
 \nabla^{\xi}= {\rm d}-\frac{H{\rm d}w}{\xi w^3}+\Lambda(\xi)\frac{{\rm d}w}{w} + \text{regular terms}.
\end{equation*}
Since $\nu$ is holomorphic only up to terms that decrease faster than any polynomial (see Corollary~\ref{exthol}), one also easily checks that we can find a gauge transformation $g$ such that $g(w=0)=\text{Id}$, $\nu \cdot g$ is holomorphic, and in the holomorphic frame $\nu \cdot g$ we get the expression \eqref{holform}.

Hence, we obtain the required compatibly framed meromorphic bundle \smash{$\bigl(\mathcal{Q}_a\mathcal{E}^{\xi},\nabla^{\xi},\nu^{\xi}_{a}\bigr)$}. Notice that the compatible frame is specified by either $\nu$, $\nu \cdot g$, or the $\nu_i$.
\end{proof}

\begin{Proposition}The Stokes data associated to \smash{$\bigl(\mathcal{Q}_a\mathcal{E}^{\xi},\nabla^{\xi},\nu^{\xi}_{a}\bigr)$} is the same as the Stokes data of \smash{$\bigl(\mathcal{P}^h_a\mathcal{E}^{\xi},\nabla^{\xi},\tau^{\xi}_{a}\bigr)$}.

\end{Proposition}
In the statement of the proposition, it is assumed that we are using the holomorphic coordinate $w=\frac{1}{z}$ vanishing at $z=\infty$, and the same branch of the Log to the define the Stokes data.

\begin{proof}
It is clear that the Stokes sectors defined by $-\xi^{-1}\frac{H}{2w^2}$ agree with the ones defined by~\smash{$-\bigl(\xi^{-1}+\overline{\xi}\bigr)\frac{H}{2w^2}$}.

To construct the associated Stokes data of \smash{$\bigl(\mathcal{Q}_a\mathcal{E}^{\xi},\nabla^{\xi},\nu^{\xi}_{a}\bigr)$}, we would first extend $\nu_a^{\xi}$ to a local holomorphic framing $\widetilde{\nu}$ around $w=0$, where the connection has the form in \eqref{holform}. Then we would consider the frames of flat sections $\widetilde{\Phi}_i$ on the extended Stokes sectors $\mathcal{S}_i$, where
\begin{equation*}
 \widetilde{\Phi}_i=\widetilde{\nu}\cdot \Sigma_i\bigl(\widehat{G}\bigr)w^{-\Lambda(\xi)}\exp \left(-\xi^{-1}\frac{H}{2w^2}\right),
\end{equation*}
and where $\widehat{G}$ and $\Sigma_i\bigl(\widehat{G}\bigr)$ satisfy the required properties given in Section \ref{GMD}.

On the other hand, we have the frames of flat sections $\Phi_i$ on $\mathcal{S}_i$ defined by
\begin{equation*}
 \Phi_i=\tau^{\xi}_a \cdot \Sigma_{i}\bigl(\widehat{F}\bigr)w^{-\Lambda(\xi)} \exp \left(-\frac{\xi^{-1}+\overline{\xi}}{2w^2}H\right)=\nu_{i}\cdot w^{-\Lambda(\xi)}\exp \left(-\xi^{-1}\frac{H}{2w^2}\right).
\end{equation*}

We claim that $\Phi_i=\widetilde{\Phi}_i$. In order to see this, let $C$ be the (constant) matrix relating both frames of flat sections (i.e., $\Phi_i=\widetilde{\Phi}_i \cdot C$) and notice that both $\nu_i$ and $\widetilde{\nu}$ go to the compatible framing \smash{$\nu_a^{\xi}$} as $w\to 0$. Using this and that fact that \smash{$\Sigma_i\bigl(\widehat{G}\bigr)\to 1$} as $w\to 0$, we conclude that we must have{\samepage
\begin{equation*}
 w^{-\Lambda(\xi)}\exp \left(-\xi^{-1}\frac{H}{2w^2}\right)C\exp \left(\xi^{-1}\frac{H}{2w^2}\right)w^{\Lambda(\xi)} \to 1\qquad \text{as} \quad w\to 0,\quad w\in \mathcal{S}_i.
\end{equation*}
The same argument at the end of Lemma \ref{UniqueA} then let us conclude that $C=1$, so that $\Phi_i=\widetilde{\Phi}_i$.}

Since the frames of flat sections $\Phi_i$ define the Stokes matrices of \smash{$\bigl(\mathcal{P}^h_a\mathcal{E}^{\xi},\nabla^{\xi},\tau^{\xi}_{a}\bigr)$}, and the frames \smash{$\widetilde{\Phi}_i$} define the Stokes matrices of \smash{$\bigl(\mathcal{Q}_a\mathcal{E}^{\xi},\nabla^{\xi},\nu^{\xi}_{a}\bigr)$}, we conclude what we want.
\end{proof}
\end{subsubsection}

\begin{subsubsection}{Gluing together the deformed flat meromorphic bundles}

Now we would like to address the issue of how the deformed meromorphic flat bundles $\bigl(\mathcal{Q}_a\mathcal{E}^{\xi},\nabla^{\xi}\bigr)$ glue together as a family parametrized by $\xi \in \mathbb{C}^{*}$.

Given the holomorphic bundle \smash{$\mathcal{E}^{\xi}=\bigl(E,\overline{\partial}_E+\xi \theta^{\dagger_h}\bigr)$} over $\CP\setminus \{\infty\}$ and given the canonical projection $p\colon\CP \times \mathbb{C}^*\to \CP$, we can form the holomorphic bundle $\mathcal{E}=\bigl(p^*E,\overline{\partial}_E+\xi \theta^{\dagger_h}+\overline{\partial}_{\xi}\bigr)$ over $\bigl(\CP\setminus \{\infty\} \bigr)\times \mathbb{C}^*$. Let $\xi_0 \in \mathbb{C}^*$, and $U(\xi_0)$ a small neighborhood of $\xi_0$. We would like to perform an extension of $(\mathcal{E},\nabla)\to \bigl(\CP\setminus \{\infty\}\bigr)\times U(\xi_0)$ to a filtered flat meromorphic bundle over $\CP\times U(\xi_0)$, with poles along $\{\infty\}\times U(\xi_0)$. Here \smash{$\nabla\colon \mathcal{E}\to \mathcal{E}\otimes \Omega^1_{(\CP\setminus \{\infty\})\times U(\xi_0)/U(\xi_0)}$} is~considered as a relative connection (it does not differentiate in the $\xi$ direction), and $\nabla^{\xi}=\nabla|_{\mathcal{E}^{\xi}}$.

\begin{Theorem} For $U(\xi_0)$ small enough, the holomorphic bundle with relative flat connection $(\mathcal{E},\nabla)\to \bigl(\CP \setminus \{\infty\}\bigr)\times U(\xi_0)$ extends to a filtered bundle with a $($relative$)$ flat meromorphic connection \smash{$\bigl(\mathcal{Q}_*^{(\xi_0)}\mathcal{E},\nabla\bigr)\to \bigl(\CP\times U(\xi_0),\{\infty\}\times U(\xi_0)\bigr)$}. Furthermore, we have that \smash{$\bigl(\mathcal{Q}_a^{(\xi_0)}\mathcal{E},\nabla\bigr)|_{\CP\times \{\xi\} }= \bigl(\mathcal{Q}_a\mathcal{E}^{\xi},\nabla^{\xi}\bigr)$}.
\end{Theorem}

\begin{proof}
This follows from \cite[Proposition 9.2.1, Corollary 9.2.5, and Theorem 11.1.2]{M}.
\end{proof}

\begin{Theorem} \label{holStokes} Let $\nu$ be a holomorphic section of \smash{$\mathcal{Q}_a^{(\xi_0)}\mathcal{E}|_{\{\infty\}\times U(\xi_0)}$} such that $\nu(\xi)$ is a compatible framing for $\bigl(\mathcal{Q}_a\mathcal{E}^{\xi},\nabla^{\xi}\bigr)$. Then the sectorial frames of flat sections $\Phi_i$ used to build the Stokes data associated to $\bigl(\mathcal{Q}_a\mathcal{E}^{\xi},\nabla^{\xi},\nu(\xi)\bigr)$ vary holomorphically in $\xi$. Hence the Stokes data itself also varies holomorphically in $\xi$.
\end{Theorem}

\begin{proof}
Given $\nu(\xi)$, we can extend it to a holomorphic local trivialization of \smash{$\mathcal{Q}_a^{(\xi_0)}$} such that in that local frame we have
\begin{equation*}
 \nabla^{\xi}= {\rm d} -\frac{H}{\xi}\frac{{\rm d}w}{w^3}+\Lambda(\xi)\frac{{\rm d}w}{w} + A(w,\xi){\rm d}w,
\end{equation*}
where $A(w,\xi)$ is holomorphic in both variables. The result now follows from \cite[Lemma~7 and Corollary~8]{B2}.
\end{proof}

This does not prove the holomorphic dependence of $\mathcal{X}_m(\xi)$ in $\xi$, since we do not know that the compatible frames \smash{$\nu_a^{\xi}$} from Section \ref{defirr} glue together holomorphically in $\xi$.
\end{subsubsection}

\begin{subsubsection}{Proof of the holomorphic dependence}\label{proofhol}

Pick $\xi_0 \in \mathbb{C}^*$ and let \smash{$\bigl(\mathcal{Q}^{(\xi_{0})}_{a}\mathcal{E},\nabla\bigr)$} be as before. We also fix some $i=1,2,3,4$ and pick $U(\xi_0)$ small enough so that there is a sector $S_i$ centered at $z=\infty$ such that
\begin{itemize}\itemsep=0pt
 \item \smash{$S_i\subset \widehat{\operatorname{Sect}}_i(\xi)$} for $\xi \in U(\xi_0)$.
 \item $S_i$ contains the Stokes ray and two anti-Stokes rays in the interior of \smash{$\widehat{\operatorname{Sect}}_i(\xi)$}.
\end{itemize}

Furthermore, let $\nu(\xi)$ be a holomorphic frame of \smash{$\mathcal{Q}_a^{(\xi_0)}\mathcal{E}|_{\{\infty\}\times U(\xi_0)}$} such that $\nu(\xi)$ is a compatible frame for $\bigl(\mathcal{Q}_a\mathcal{E}^{\xi},\nabla^{\xi}\bigr)$ for each $\xi \in U(\xi_0)$ (the fact that such a frame exists follows from the third and fourth statement of \cite[Theorem 11.1.2]{M}). Then (after a possible reordering of the elements of the frame) $\nu(\xi)$ extends to a holomorphic local frame of \smash{$\mathcal{Q}_a^{(\xi_0)}\mathcal{E}$}, where
\begin{equation*}
 \nabla^{\xi}= {\rm d} -\frac{H}{\xi}\frac{{\rm d}w}{w^3}+\Lambda(\xi)\frac{{\rm d}w}{w} + A(w,\xi){\rm d}w
\end{equation*}
as in Theorem \ref{holStokes}. We can then construct for $\xi\in U(\xi_0)$ the frame of flat sections $\widetilde{\Phi}_i(\xi)$ on $S_i$ used to build the Stokes data of $\bigl(\mathcal{Q}_a\mathcal{E}^{\xi},\nabla^{\xi},\nu(\xi)\bigr)$. On the other hand, we have the frame of flat sections $\Phi_i(\xi)$ of \smash{$\bigl(\mathcal{Q}_a\mathcal{E}^{\xi},\nabla^{\xi},\nu_a^{\xi}\bigr)$} used to build the magnetic twistor coordinate. We clearly have that $\widetilde{\Phi}_i(\xi)=\Phi_i(\xi)\cdot C_i(\xi)$ for some matrix $C_i(\xi)$ that depends only on $\xi \in U(\xi_0)$.

\begin{Lemma}\label{diag} The matrix $C_i(\xi)$ is diagonal.

\end{Lemma}

\begin{proof}
For each $\xi \in U(\xi_0)$, the frames $\nu_a^{\xi}$ and $\nu(\xi)$ are compatible frames for the same meromorphic flat bundle $\bigl(\mathcal{Q}_a\mathcal{E}^{\xi},\nabla^{\xi}\bigr)$, with associated irregular type
\smash{$
 -\frac{H}{\xi}\frac{{\rm d}w}{w^3}+\Lambda(\xi)\frac{{\rm d}w}{w}$}.
Hence, we must have that \smash{$\nu(\xi)|_{\infty}=\nu_a^{\xi}|_{\infty}\cdot D(\xi)$} for some diagonal matrix $D(\xi)$ depending only on $\xi$.

By applying an argument like the one found in Lemma \ref{UniqueA}, we have that
\begin{equation*}
 w^{-\Lambda(\xi)}\exp \left(-\xi^{-1}\frac{H}{2w^2}\right)C_i(\xi)\exp \left(\xi^{-1}\frac{H}{2w^2}\right)w^{\Lambda(\xi)} \to D(\xi) \qquad \text{as} \quad w\to 0,\quad w\in S.
\end{equation*}
By the choice of sector $S$ containing a Stokes ray for each $\xi \in U(\xi_0)$, we have that the off-diagonal entries of $C_i(\xi)$ must be $0$, and the diagonal entries must match the diagonal entries of~$D(\xi)$. Hence, we conclude what we want.
\end{proof}

Next, we will need the following lemma.

\begin{Lemma} \label{Lemma} The asymptotics in Lemma {\rm\ref{twistasymp}} of the frames of flat sections $\Phi_i(\xi)$ used in the definition of $\mathcal{X}_m(\xi)$ hold uniformly in $\xi \in U(\xi_0)$ for $U(\xi_0)$ small enough and bounded.
\end{Lemma}

\begin{proof}
The proof of this is in Appendix \ref{AC}.
\end{proof}

\begin{Theorem} The coordinate $\MC(\xi)$ depends holomorphically on $\xi$.
\end{Theorem}

\begin{proof}
It is enough to prove that the frames of flat sections $\Phi_i(\xi)$ depend holomorphically on~$\xi$.

We know by Lemma \ref{diag} that $\widetilde{\Phi}_i(\xi)=\Phi_i(\xi)\cdot C_i(\xi)$ for $C_i(\xi)$ diagonal. Furthermore, we have that $\overline{\partial}_{\xi}\widetilde{\Phi}_i(\xi)=0$, so if we show that $\overline{\partial}_{\xi}C_i(\xi)=0$, then we would be able to conclude what we want.

In the following, we will use the same notation of Lemma \ref{twistasymp}. By Lemma \ref{Lemma}, we then have that
\smash{$A_i(w,\xi)\cdot {\rm e}^{-Q(\xi,w)} \to 1$} uniformly in $\xi \in U(\xi_0)$ as $ w\to 0$, $w\in S$,
where $\Phi_i=(e_1,e_2)\cdot A_i(w,\xi)$.

Since $\widetilde{\Phi}_i(\xi)$ is holomorphic in $\xi$ and the frame $(e_1,e_2)$ is also holomorphic in $\xi$ (since it does not depend on $\xi$), we conclude that the matrix $A_i(w,\xi)\cdot C_i(\xi)$ is holomorphic in $\xi$. Furthermore, ${\rm e}^{Q(\xi)}$ is diagonal and also holomorphic in $\xi$, so we have that $\smash{A_i(w,\xi)\cdot C_i(\xi)\cdot {\rm e}^{-Q(\xi,w)}}=\smash{A_i(w,\xi)\cdot {\rm e}^{-Q(\xi,w)}\cdot C_i(\xi)}$ and \smash{$\overline{\partial}_{\xi}\bigl(A_i(w,\xi)\cdot C_i(\xi)\cdot {\rm e}^{-Q(\xi,w)}\bigr)=0$}.

By shrinking $U(\xi_0)$ if necessary, we can assume that $C_i(\xi)$ is a bounded function of $\xi$. Hence we get that
\smash{$
A_i(w,\xi) \cdot {\rm e}^{-Q(\xi)} \cdot C_i(\xi) \to C_i(\xi)
$}
uniformly in $\xi \in U(\xi_0)$ as $w\to 0$ along $S_i$. Since~\smash{$A_i(w,\xi) \cdot {\rm e}^{-Q(\xi)} \cdot C_i(\xi)$} is holomorphic in $\xi$, we conclude that $C_i(\xi)$ must be a holomorphic function in $\xi$, so $\Phi_i(\xi)$ must depend holomorphically on $\xi$.
\end{proof}
\end{subsubsection}

\subsection{Twistor coordinates Part 3: Asymptotics in the twistor parameter}

In this section, we compute the asymptotics of $\MC(\xi)$ as $\xi \to 0$ and $\xi \to \infty$. More precisely, we will verify the asymptotics have a formula similar to \eqref{aft}. Along the way, we also show that the reality condition for $\MC(\xi)$ holds (see \eqref{rc}).

The plan for computing the asymptotics will be the following:
\begin{itemize}\itemsep=0pt
 \item First, we express the sectorial flat sections used in the definition of Stokes data, in a~convenient way for studying the asymptotics as $\xi \to 0$. This is done in Lemma \ref{asympform}.
 \item Then, we apply Lemma \ref{asympform} to the actual computation of the asymptotics of $\MC(\xi)$ as~${\xi \to 0}$ in Theorem \ref{a0}.
 \item Finally, we prove the reality condition for $\MC(\xi)$, and use it to compute the asymptotics as $\xi\to \infty$. This is Theorem~\ref{rct} (resp.~Corollary \ref{rcc}).
\end{itemize}

\begin{Notation} For the rest of this section, we will use the following notation: $\operatorname{Arg}(w)$, $\operatorname{Arg}_m(w)$ and $\operatorname{Log}_m(w)$ are defined as in Definition \ref{deftwistcoords}; $\operatorname{Arg}_P(w)$ and $\operatorname{Log}_P(w)$ are the principal branches (i.e., with $\operatorname{Arg}_P(w)\in (-\pi,\pi]$). When we switch coordinates to $z=1/w$, we denote by $\operatorname{Log}_p(z)$ and $\operatorname{Arg}_p(z)$ the branches with $[-\pi,\pi)$. We then have the relations $\operatorname{Log}_{p}(z)=-\operatorname{Log}_P(w)$ and~${\operatorname{Arg}_{p}(z)=-\operatorname{Arg}_P(w)}$.
\end{Notation}

\begin{subsubsection}{Preliminaries for the asymptotic computation}

Here we develop some of the preliminary notation and computations that we will need for the main asymptotic computation of the twistor magnetic coordinate.

We first start with some results that will help us understand the asymptotic behaviour of the exponentially decreasing flat sections of $\nabla^{\xi}$ along certain special curves.

\begin{Definition}We fix the quadratic differential $\phi:=\bigl(z^2+2m \bigr){\rm d}z^2$ over $\CP$, with $m\in \mathbb{C}^*$. Given a phase ${\rm e}^{{\rm i}\theta} \in S^1$, a WKB curve\footnote{This terminology comes from \cite{GMN2}; a more common terminology is that of horizontal trajectory for ${\rm e}^{-2{\rm i}\theta}\phi$.} with phase ${\rm e}^{{\rm i}\theta}$ is a parametrized curve $\gamma$ in $\CP$ such that its velocity vector $\dot{\gamma}$ satisfies $
 \phi (\dot{\gamma})= {\rm e}^{2{\rm i}\theta}$.
\end{Definition}

For the rest of the section, we fix:
\begin{itemize}\itemsep=0pt
 \item $m\in \mathbb{C}^*$.
 \item An element $\bigl(E,\overline{\partial}_E,\theta,h,(e_1,e_2)\bigr)\in \mathcal{H}^{{\rm fr}}$ with $\operatorname{Det}(\theta)=-\bigl(z^2+2m\bigr){\rm d}z^2$.
 \item A WKB curve $\gamma$ with phase ${\rm e}^{{\rm i} \operatorname{Arg}(m)}$ not running into the zeroes of $\phi$.
 \item A frame $(\eta_1,\eta_2)$ of eigenvectors of $\theta$ along the curve $\gamma$, with growth near $w=\frac{1}{z}=0$ of the form $|\eta_i|_h=\mathcal{O}(|w|^{a_i})$ for some $a_i \in \mathbb{R}$. For example, we could take the frame compatible with the parabolic structure, or a normalized frame.
\end{itemize}

With this data fixed, we consider the flatness equation of the pullback connection $\gamma^*\nabla^{\xi}$
\begin{equation*}
 \frac{{\rm d}}{{\rm d}t}{\rm d}t +\xi^{-1}\gamma^*\theta + \gamma^*A + \xi \gamma^* \theta ^{\dagger_h}=0,
\end{equation*}
where $A$ denotes the connection form of $D\bigl(\overline{\partial}_E,h\bigr)$ in the frame $(\eta_1,\eta_2)$.

\begin{Proposition} \label{prelimasymp} Let $a,a_0\in \mathbb{R}$ and let $M(t):=\exp \bigl(-\int_{a}^{t}\gamma^*A_{\operatorname{diag}}\bigr)$ denote the diagonal gauge transformation that gauges away the diagonal part of $\gamma^*A$.\footnote{Namely, if $\gamma^*A_{\operatorname{diag}}$ denotes the diagonal components of $\gamma^*A$ in the frame $(\eta_1,\eta_2)$, then one can check that in the new frame $(\eta_1,\eta_2)\cdot M(t)$ the connection form $M(t)^{-1}\gamma^*AM(t)+M^{-1}{\rm d}M(t)$ is off-diagonal.} Furthermore, let
\begin{equation*}
 \lambda_i(t,\xi)= -\xi^{-1}\gamma^*\theta_{ii} - \xi \gamma^* \theta^{\dagger_h}_{ii},
\end{equation*}
and assume that either $\textnormal{Re}(\lambda_1(t,\xi))>\textnormal{Re}(\lambda_2(t,\xi))$ or $\textnormal{Re}(\lambda_1(t,\xi))<\textnormal{Re}(\lambda_2(t,\xi))$ holds for all~${t\in \mathbb{R}}$ and for small enough $\xi \in \mathbb{H}_m$. Then there is a neighborhood $U_0$ of $\xi=0$ such that for~${\xi \in U_0\cap \mathbb{H}_m}$, we can write an exponentially decreasing flat section of $\gamma^*\nabla^{\xi}$ for $t\in (a_0,\infty)$ with the following form:
\begin{itemize}\itemsep=0pt
 \item In the case $\textnormal{Re}(\lambda_1(t,\xi))>\textnormal{Re}(\lambda_2(t,\xi))$, we have
 \begin{equation*}
 s(t,\xi)=\exp \left(\int_a^{t}\lambda_2(\tau,\xi){\rm d}\tau\right)(E_1(t,\xi)\eta_1 + (M_{22}(t)+E_2(t,\xi))\eta_2).
 \end{equation*}
 \item And in the case $\textnormal{Re}(\lambda_1(t,\xi))<\textnormal{Re}(\lambda_2(t,\xi))$,
 \begin{equation*}
 s(t,\xi)=\exp \left(\int_a^{t}\lambda_1(\tau,\xi){\rm d}\tau\right)((M_{11}(t)+E_1(t,\xi))\eta_1 + E_2(t,\xi)\eta_2).
 \end{equation*}
\end{itemize}
In both cases $E_i(t,\xi)$ satisfies
\begin{itemize}\itemsep=0pt
 \item For fixed $t \in (a_0,\infty)$ we have that $E_i(t,\xi) \to 0$ as $\xi \to 0$.
 \item $E_i(t,\xi)\to 0$ as $t \to \infty$ uniformly in $\xi \in \mathbb{H}_m\cap U_{0}$.
\end{itemize}
\end{Proposition}
\begin{proof}
In Appendix \ref{AD}. The proof of this uses classical techniques in the theory of ordinary differential equations depending of parameters.
\end{proof}

Let us now state the setting for our asymptotic computation:
\begin{itemize}\itemsep=0pt
 \item Let \smash{$\lambda:=\sqrt{z^2+2m} {\rm d}z$} be the square root of $\phi$ with branch cut given by the line segment between \smash{$z=\pm\sqrt{-2m}$}. $\lambda$ has power series expansion centered at $z=\infty$ given by $(z+m/z+\cdots ){\rm d}z=\bigl(-1/w^3 -m/w +\cdots \bigr){\rm d}w$.
 \item We let $\gamma(t)$ be a WKB curve with phase ${\rm e}^{{\rm i} \operatorname{Arg}(m)}$, such that $\lambda(\dot{\gamma})={\rm e}^{{\rm i}\operatorname{Arg}(m)}$. Notice that we can take, and will take, $\gamma(t)$ of the form $\gamma(t)=g(t){\rm e}^{{\rm i}\operatorname{Arg}(m)/2}$, where $g(t)$ is a certain real valued function satisfying $g(t)\to \pm \infty$ (or $\mp \infty)$ as $t\to \pm \infty$. With this choice, $\gamma(t)$ crosses the branch cut of $\lambda$ (see Figure~\ref{fig3}).\footnote{The trajectory structure of the WKB curves decomposes the sphere into four half-planes and one stripe, and here we are taking one particular WKB curve lying in the stripe.}
 \begin{figure}[ht]
 \centering
 \includegraphics[width=7cm]{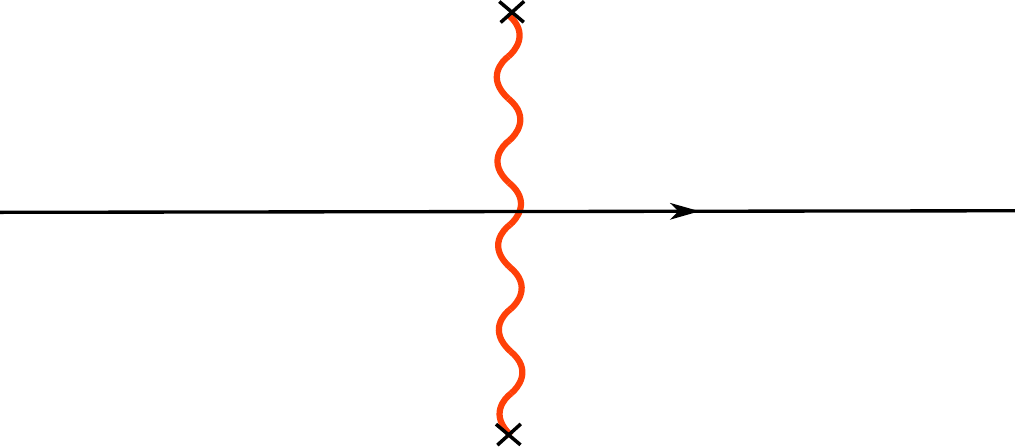}
 \caption{The crosses denote $\pm \sqrt{-2m}$, the wavy red line denotes the branch cut of $\lambda$, and the horizontal black line the WKB path.}
 \label{fig3}
 \end{figure}
 \item For $\xi \in \mathbb{H}_m$ and for sufficiently big $t$, $\gamma(t)$ lies either on \smash{$\widehat{\operatorname{Sect}}_1(\xi)\cap\widehat{\operatorname{Sect}}_4(\xi)$} or in $\smash{\widehat{\operatorname{Sect}}_2(\xi)}\cap\smash{\widehat{\operatorname{Sect}}_3(\xi)}$. We orient $\gamma$ such that it lies in \smash{$\widehat{\operatorname{Sect}}_1(\xi)\cap\widehat{\operatorname{Sect}}_4(\xi)$} for sufficiently big $t$.
 \item Let \smash{$\widetilde{\lambda}$} be a branch of a square root of $\phi$ defined in a neighborhood of $\gamma$, such that ${\widetilde{\lambda}(\dot{\gamma})=\lambda(\dot{\gamma})}$ for big enough $t$. With our choices, we then have that $\operatorname{Re}\bigl(\frac{1}{\xi}\widetilde{\lambda}(\dot{\gamma})\bigr)>0$ for~${\xi \in \mathbb{H}_m}$.
 \item
 We pick a frame $(\eta_1,\eta_2)$ of ${\rm SL}(2,\mathbb{C})$ eigenvectors of $\theta$ along $\gamma$, such that $(\eta_1,\eta_2) \to (e_1,e_2)$ as $t\to \infty$, and $(\eta_1,\eta_2)\to (e_2,-e_1)$ as $t \to -\infty$. In the frame $(\eta_1,\eta_2)$ we have
\smash{$
 \theta =\bigl[\begin{smallmatrix} \widetilde{\lambda} & 0\\
 0 & -\widetilde{\lambda}
 \end{smallmatrix}\bigr]$}.

 \item We pick $t_0<t_1$, such that $t_0$ is a time before the crossing with the branch cut of $\lambda$, and~$ t_1$ is a time after the crossing with the branch cut of $\lambda$.
\end{itemize}

We now state the lemma that will allow us to compute the asymptotics of the magnetic twistor coordinate:

\begin{Lemma}\label{asympform} In the previous setting, let $s_i(z,\xi)$ be the flat sections from Section {\rm\ref{deftwistcoords}}. They satisfy the following:

For small enough $\xi \in \mathbb{H}_m$
\begin{align*}
 s_1(\gamma(t),\xi)={}& \exp \left(\int_{\gamma(t_1)}^{\gamma(t)}\left(-\frac{1}{\xi}\widetilde{\lambda}+\mathcal{O}(\xi)\right)\right)\\
 &\times
 \bigl((M_{11}^{(1)}(\gamma(t))+E^{(1)}_{1}(\gamma(t),\xi))\eta_1 + E^{(1)}_{2}(\gamma(t),\xi)\eta_2\bigr)\beta_1(\gamma(t_1),\xi),\\
 s_3(\gamma(t),\xi)={}& \exp \left(\int_{\gamma(t_0)}^{\gamma(t)}\left(\frac{1}{\xi}\widetilde{\lambda}+\mathcal{O}(\xi)\right)\right)\\
 &\times\bigl((M_{22}^{(3)}(\gamma(t))+E^{(3)}_{2}(\gamma(t),\xi))\eta_2 + E^{(3)}_{1}(\gamma(t),\xi)\eta_1\bigr)\beta_3(\gamma(t_0),\xi).
\end{align*}

For small enough $\xi \in \mathbb{H}_{-m}$
\begin{align*}
 s_2(\gamma(t),\xi)={}& \exp \left(\int_{\gamma(t_1)}^{\gamma(t)}\left(\frac{1}{\xi}\widetilde{\lambda}+\mathcal{O}(\xi)\right)\right)\\
 &\times
 \bigl((M_{22}^{(2)}(\gamma(t))+E^{(2)}_{2}(\gamma(t),\xi))\eta_2 + E^{(2)}_{1}(\gamma(t),\xi)\eta_1\bigr)\beta_2(\gamma(t_1),\xi),\\
 s_4(\gamma(t),\xi)={}& \exp \left(\int_{\gamma(t_0)}^{\gamma(t)}\left(-\frac{1}{\xi}\widetilde{\lambda}+\mathcal{O}(\xi)\right)\right)\\
 &\times
 \bigl((M_{11}^{(4)}(\gamma(t))+E^{(4)}_{1}(\gamma(t),\xi))\eta_1 + E^{(4)}_{2}(\gamma(t),\xi)\eta_2\bigr)\beta_4(\gamma(t_0),\xi).
\end{align*}
In the above, $\beta_i(\gamma(t_j),\xi)$ are normalization constants that ensure that the flat sections have the correct asymptotics as $t \to \infty$ $($resp.\ $t\to -\infty)$ for $i=1,2$ $($resp.\ $i=3,4)$, and where \smash{$M^{(i)}_j(t)$} and \smash{$E^{(i)}_j(t,\xi)$} have the same role and properties as in Proposition {\rm\ref{prelimasymp}}.

\end{Lemma}

\begin{proof} We will give a proof for $s_1(\gamma(t),\xi)$, since the others follow similarly.

By Proposition \ref{prelimasymp} and the setting of our computation, we can write an exponentially decreasing flat section for $\gamma^*\nabla^{\xi}$ by
\begin{equation*}
 s(t,\xi)=\exp \left(\int_{\gamma(t_1)}^{\gamma(t)}\left(-\frac{1}{\xi}\widetilde{\lambda}+\mathcal{O}(\xi)\right)\right)((M_{11}(t)+E_1(t,\xi))\eta_1 + E_2(t,\xi)\eta_2),
\end{equation*}
where we are denoting \smash{$\mathcal{O}(\xi)=-\xi \theta^{\dagger_h}_{11}$}. Note that because of Theorem \ref{FT1}, we have that \smash{$\theta^{\dagger_h}_{11}=\overline{\theta_{11}}^{t} + \phi$}, where $\phi$ is exponentially decreasing as $z\to \infty$.

By our definition of $\gamma(t)$, when $t\to \infty$ we have that $\gamma(t)$ lies in $\widehat{\operatorname{Sect}}_1(\xi)\cap \widehat{\operatorname{Sect}}_4(\xi)$. Since exponentially decreasing flat sections along such sectors are uniquely determined up to scaling, we have that $s_1(\gamma(t),\xi)=c(\xi,\gamma(t_1))s(t,\xi)$ for some number $c(\xi,\gamma(t_1))$.

On the other hand, it is not hard to check that
\smash{$
 \operatorname{Arg}_m(\gamma(t))\to \operatorname{Arg}_m\bigl(\sqrt{-m^{-1}}\bigr)-\frac{\pi}{2}$} as $t\to \infty$,
where \smash{$\sqrt{-m^{-1}}$} uses the principal branch in the $w$-coordinate. Furthermore, because of our conventions, the following holds (see Figure~\ref{fig5})
\smash{$
 \operatorname{Arg}_m\bigl(\sqrt{-m^{-1}}\bigr)=\operatorname{Arg}_P\bigl(\sqrt{-m^{-1}}\bigr)$}.

\begin{figure}
 \centering
 \includegraphics[width=14cm]{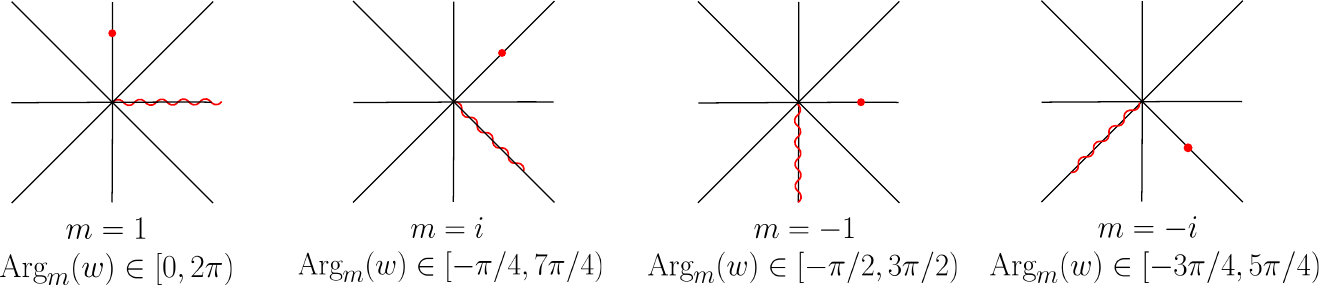}
 \caption{The pictures show that $\operatorname{Arg}_m\bigl(\sqrt{-m^{-1}}\bigr)=\operatorname{Arg}_P\bigl(\sqrt{-m^{-1}}\bigr)$ for $m\in \{1,{\rm i},-1,-{\rm i}\}$, but similar pictures hold for any $m$. The red dot in the pictures denotes $\sqrt{-m^{-1}}$, where $\sqrt{w}$ uses the principal branch. The wavy red line denotes the branch cut of $\operatorname{Arg}_m$.}
 \label{fig5}
\end{figure}

Let $F(w)$ be the antiderivative of $-\lambda=\bigl(\frac{1}{w^3} +m/w +\mathcal{O}(w)\bigr){\rm d}w$ of the form $-\frac{1}{2w^2}+m\operatorname{Log}_{m}(w)+\mathcal{O}(w)$. Furthermore, let \smash{$\widetilde{F}(w)$} be the antiderivative of the $\mathcal{O}(\xi)$ term of the form~\smash{$-\frac{1}{2\overline{w}^2} +\overline{m}\overline{\operatorname{Log}_{m}(w)}+\mathcal{O}(|w|)$}. Now define the following:
\begin{gather*}
 \beta_{1,{\xi^{-1}}}(t_1):=\exp \left(\frac{1}{\xi}F\left(\frac{1}{\gamma(t_1)}\right)\right),\\
 \beta_{1,{\xi^{0}}}(t_1):=\exp \left(- {\rm i}m^{(3)}\left(\operatorname{Arg}_P\bigl(\sqrt{-m^{-1}}\bigr)-\frac{\pi}{2}\right)+\int_{t_1}^{\infty}\gamma^*A_{11}\right),\\
 \beta_{1,\xi}(t_1):=\exp \left(\xi \widetilde{F}\left(\frac{1}{\gamma(t_1)}\right)\right),\qquad
 \beta_1(t_1,\xi):=\beta_{1,{\xi^{-1}}}(t_1)\beta_{1,{\xi^{0}}}(t_1)\beta_{1,\xi}(t_1).
\end{gather*}

It is then easy to check that $\beta_1(t_1,\xi)s(t,\xi)$ has the same asymptotics as $t\to \infty$ as $s_1(\gamma(t),\xi)$ (see Lemma \ref{twistasymp}). So $s_1(\gamma(t),\xi)=\beta_1(t_1,\xi)s(t,\xi)$.
\end{proof}

For completeness, and because we will need it in the asymptotic computation for the magnetic twistor coordinate, we write the other normalization constants $\beta_i(t_j,\xi)$. These are
\begin{align*}
 \beta_2(t_1,\xi):={}&\exp \left(-\frac{1}{\xi}F\left(\frac{1}{\gamma(t_1)}\right) +{\rm i}m^{(3)}\left(\operatorname{Arg}_P\bigl(\sqrt{-m^{-1}}\bigr)-\frac{\pi}{2}\right)\right.\\
 &\left.+\int_{t_1}^{\infty}\gamma^*A_{22} - \xi \widetilde{F}\left(\frac{1}{\gamma(t_1)}\right)\right),\\
 \beta_3(t_0,\xi):={}&\exp \left(\frac{1}{\xi}F\left(\frac{1}{\gamma(t_0)}\right)- {\rm i}m^{(3)}\left(\operatorname{Arg}_P\bigl(\sqrt{-m^{-1}}\bigr)+\frac{\pi}{2}\right)+{\rm i}\pi\right.\\
&\left.+\int_{t_0}^{-\infty}\gamma^*A_{22} + \xi \widetilde{F}\left(\frac{1}{\gamma(t_0)}\right)\right),\\
 \beta_4(t_0,\xi):={}&\exp \left(-\frac{1}{\xi}F\left(\frac{1}{\gamma(t_0)}\right) +{\rm i}m^{(3)}\left(\operatorname{Arg}_P\bigl(\sqrt{-m^{-1}}\bigr)+\frac{\pi}{2}\right)\right.\\
& \left.+\int_{t_0}^{-\infty}\gamma^*A_{11} - \xi \widetilde{F}\left(\frac{1}{\gamma(t_0)}\right)\right),
\end{align*}
where the extra ${\rm i}\pi$ in $\beta_3$ comes because of the fact that $\eta_2 \to -e_1$ instead of $e_1$ as $t\to -\infty$.

We should justify why the integrals of the form
\begin{equation}\label{finitenessint}
 \int_{t_j}^{\pm \infty} \gamma^*A_{ii}
\end{equation}
are finite. Let $w=r{\rm e}^{{\rm i}\theta}$, and write $A_{ii}=A_{ii,r}{\rm d}r+A_{ii,\theta}{\rm d}\theta$. By our choice of $\gamma(t)$, we have that~${{\rm d}\theta(\dot{\gamma})=0}$. On the other hand, $A_{ii,r}$ is regular at $w=0$ and ${\rm d}r(\dot{\gamma})\sim t^{-3/2}$ as $t \to \pm \infty$. Hence, the integrals of the above form converge.
\end{subsubsection}

\begin{subsubsection}[Computing the asymptotics when xi to 0]{Computing the asymptotics when $\boldsymbol{\xi \to 0}$}\label{asymp0}

We are finally ready to start computing the asymptotics of the magnetic twistor coordinate.

\begin{Theorem} \label{a0} Consider the magnetic twistor coordinate $\MC(\xi)$ associated to $\bigl[\bigl(E,\overline{\partial}_E,\theta,h,\allowbreak(e_1,e_2)\bigr)\bigr]\in \mathfrak{X}^{{\rm fr}}$. Then with the setting and notations of the previous section, we have
\begin{equation*}
 \MC(\xi)\sim_{\xi \to 0} A(\xi),
\end{equation*}
where
\begin{equation*}
A(\xi):=\exp \left(-\frac{1}{\xi}m\operatorname{Log}_p\left(-\frac{m}{2e}\right)+{\rm i}m^{(3)}\operatorname{Arg}_p(-m)+{\rm i}\pi
+\int_{-\infty}^{\infty}\gamma^*A_{11}\right).
\end{equation*}

\end{Theorem}

\begin{proof}
First assume that $\xi \in \mathbb{H}_m$. We pick $\xi$ to be small enough so that the intervals where the expressions for $s_1$ and $s_3$ of Lemma \ref{asympform} hold, overlap for some $t\in (t_0,t_1)$. By using Lemma~\ref{asympform}, we have that
\begin{align*}
 s_3(\gamma(t),\xi)\wedge s_1(\gamma(t),\xi)
 ={}& \exp \left(\int_{\gamma(t_0)}^{\gamma(t_1)}\left(\frac{1}{\xi}\widetilde{\lambda}+\mathcal{O}(\xi)\right)\right)\beta_1(\gamma(t_1),\xi)\beta_3(\gamma(t_0),\xi)\\
 &\times\bigl((M^{(1)}_{11}+E_{1}^{(1)})(M_{22}^{(3)}+E_{2}^{(3)})-E_1^{(3)}E_2^{(1)}\bigr)\eta_2\wedge \eta_1.
\end{align*}
It is then easy to see that the following asymptotic computation holds
\begin{gather*}
 \MC(\xi)=a(\xi)=\frac{s_3(\gamma(t),\xi)\wedge s_1(\gamma(t),\xi)}{s_2(\gamma(t),\xi)\wedge s_1(\gamma(t),\xi)}\\
 \phantom{ \MC(\xi)=}{}\sim_{\xi \to 0, \xi \in \mathbb{H}_m}\frac{\exp \bigl(\frac{1}{\xi}\int_{\gamma(t_0)}^{\gamma(t_1)}\widetilde{\lambda}\bigr)\beta_{1,\xi^{-1}}(t_1)\beta_{1,\xi^{0}}
 (t_1)\beta_{3,\xi^{-1}}(t_0)\beta_{3,\xi^{0}}(t_0)M_{11}^{(1)}(\gamma(t))}{s_2(\gamma(t),\xi)\wedge s_1(\gamma(t),\xi)}\\
\phantom{\qquad\sim_{\xi \to 0, \xi \in \mathbb{H}_m}}{} \times M_{22}^{(3)}(\gamma(t))\eta_2\wedge \eta_1.
\end{gather*}

To continue the computation, we use the following lemma.

\begin{Lemma}\label{auxasylemma} The following holds:
\begin{equation*}
 \exp \left(\frac{1}{\xi}\int_{\gamma(t_0)}^{\gamma(t_1)}\widetilde{\lambda}\right)\beta_{1,\xi^{-1}}(t_1)\beta_{3,\xi^{-1}}(t_0)
 =\exp \left(-\frac{1}{\xi}m\operatorname{Log}_p\left(-\frac{m}{2e}\right)\right).
\end{equation*}

\end{Lemma}

\begin{proof}
By deforming the path $\gamma$ to a path from $\gamma(t_0)$ to $\gamma(t_1)$ passing through $z=\sqrt{-2m}$ (where $\sqrt{-2m}$ uses $\operatorname{Arg}_p(z)$), we get that
\begin{equation}\label{auxasy}
 \exp \left(\frac{1}{\xi}\int_{\gamma(t_0)}^{\gamma(t_1)}\widetilde{\lambda}\right)=\exp \left(\frac{1}{\xi}\int_{\sqrt{-2m}}^{\gamma(t_1)}\lambda -\frac{1}{\xi}\int_{\gamma(t_0)}^{\sqrt{-2m}}\lambda\right).
\end{equation}
To compute this quantity, we use the following antiderivative $\Lambda(z)$ of \smash{$\lambda=\sqrt{z^2+2m}{\rm d}z$}
\begin{align*}
 \Lambda(z):=&{}\frac{z}{2}\sqrt{z^2+2m} + m\operatorname{Log}_p\bigl(z+\sqrt{z^2+2m}\bigr) -\frac{m}{2}-m\operatorname{Log}_p(2)\\
={}&\frac{z^2}{2}+m\operatorname{Log}_p(z)+\mathcal{O}\left(\frac{1}{z}\right).
\end{align*}
 Notice that with our conventions, $\operatorname{Log}_P(1/\gamma(t_i))=\operatorname{Log}_m(1/\gamma(t_i))$, so that is easy to check that~${F(1/\gamma(t_j))=-\Lambda(\gamma(t_j))}$. Hence, $\beta_{i,\xi^{-1}}(t_j)=\exp \bigl(-\xi^{-1}\Lambda(\gamma(t_j))\bigr)$.

Applying this to \eqref{auxasy}, we get
\begin{align*}
 \exp \left(\frac{1}{\xi}\int_{\gamma(t_0)}^{\gamma(t_1)}\widetilde{\lambda}\right)\beta_{1,\xi^{-1}}(t_1)\beta_{3,\xi^{-1}}(t_0)&=\exp \left(-\frac{2}{\xi}\Lambda\bigl(\sqrt{-2m}\bigr)\right)\\
 &=\exp \left(-\frac{1}{\xi}\bigl(2m\operatorname{Log}_p\bigl(\sqrt{-2m}\bigr) -m -2m\operatorname{Log}_p(2)\bigr)\right)\\
 &=\exp \left(-\frac{1}{\xi}m\operatorname{Log}_p\left(-\frac{m}{2e}\right)\right).\tag*{\qed}
\end{align*}\renewcommand{\qed}{}
\end{proof}

Using Lemma \ref{auxasylemma}, we get the following asymptotics as $\xi \to 0$, $\xi \in \mathbb{H}_m$
\begin{gather*}
 a(\xi)=\frac{s_3(\gamma(t),\xi)\wedge s_1(\gamma(t),\xi)}{s_2(\gamma(t),\xi)\wedge s_1(\gamma(t),\xi)}\\
 \hphantom{a(\xi)}
 \cong_{\xi \to 0, \xi \in \mathbb{H}_m}\frac{\exp \bigl(-\frac{1}{\xi}m\operatorname{Log}_p\bigl(-\frac{m}{2e}\bigr)-2{\rm i}m^{(3)}\operatorname{Arg}_P\bigl(\sqrt{-m^{-1}}\bigr)+{\rm i}\pi\bigr)\exp \left(\int_{t}^{\infty}\gamma^*A_{11}\right)}{s_2(\gamma(t),\xi)\wedge s_1(\gamma(t),\xi)}\\
 \qquad\times\exp \bigl(\int_{t}^{-\infty}\gamma^*A_{22}\bigr)\eta_2\wedge \eta_1\\
\qquad\qquad =\frac{\exp \bigl(-\frac{1}{\xi}m\operatorname{Log}_p\bigl(-\frac{m}{2e}\bigr)+{\rm i}m^{(3)}\operatorname{Arg}_p(-m)+{\rm i}\pi\bigr)\exp \bigl(\int_{t}^{\infty}\gamma^*A_{11}\bigr)}{s_2(\gamma(t),\xi)\wedge s_1(\gamma(t),\xi)}\\
 \qquad\qquad \phantom{=}{}\times\exp \left(\int_{t}^{-\infty}\gamma^*A_{22}\right)\eta_2\wedge \eta_1.
\end{gather*}

The last expression might seem to depend on $t$, but it actually does not. To show this, notice that since $(\eta_1,\eta_2)$ is an ${\rm SL}(2,\mathbb{C})$ frame and $D={\rm d}+A$ preserves the volume form, we have that~${\operatorname{Tr}(\gamma^*A)=0}$. This implies that $\exp \bigl(\int_{t}^{\infty}\gamma^*A_{11}\bigr)\exp \bigl(\int_{t}^{-\infty}\gamma^*A_{22}\bigr)=\exp \bigl(\int_{-\infty}^{\infty}\gamma^*A_{11}\bigr)$. On the other hand, we have that $s_2\wedge s_1$ is a flat section of \smash{$\operatorname{Det}\bigl(\gamma^*\nabla^{\xi}\bigr)$}, which in the $\mathrm{SL}(2,\mathbb{C})$ frame given by $\eta_1\wedge \eta_2$ along the WKB curve has the form
\begin{equation*}
\operatorname{Det}\bigl(\gamma^*\nabla^{\xi}\bigr)=\frac{{\rm d}}{{\rm d}t}{\rm d}t +\xi^{-1}\operatorname{Tr}(\gamma^* \theta)+\operatorname{Tr}(\gamma^* A) + \xi \operatorname{Tr}\bigl(\gamma^* \theta^{\dagger_h}\bigr) = \frac{{\rm d}}{{\rm d}t}{\rm d}t.
\end{equation*}

This tells us that we can write
$
s_2(\gamma(t),
\xi)\wedge s_1(\gamma(t),\xi)=c\eta_1\wedge \eta_2
$
for some constant $c$. In fact, because of the asymptotics of the flat sections, we see that $c=-1$. Hence, putting the results together we get the following asymptotics:
\begin{align*}
& \frac{s_3(\gamma(t),\xi)\wedge s_1(\gamma(t),\xi)}{s_2(\gamma(t),\xi)\wedge s_1(\gamma(t),\xi)}
\\
& \qquad{}\cong_{\xi \to 0, \xi \in \mathbb{H}_m} \exp \left(-\frac{1}{\xi}m\operatorname{Log}_p\Bigl(-\frac{m}{2e}\Bigr)+{\rm i}m^{(3)}\operatorname{Arg}_p(-m)+{\rm i}\pi
 +\int_{-\infty}^{\infty}\gamma^*A_{11}\right).
\end{align*}

So far we have computed the asymptotics as $\xi \to 0$ with $\xi \in \mathbb{H}_m$. Let us compute the asymptotics when $\xi \to 0$ with $\xi \in \mathbb{H}_{-m}$, and see that they match with the previous asymptotics.

Now we need to compute the asymptotics of
\begin{equation*}
-\frac{1}{b(\xi)}=-\frac{s_3(\gamma(t),\xi)\wedge s_2 (\gamma(t),\xi)}{s_4(\gamma(t),\xi)\wedge s_2(\gamma(t),\xi)}\qquad \text{where} \quad \xi \in \mathbb{H}_{-m}.
\end{equation*}

Following a similar computation from before using the relevant expressions for the $s_i$ and $\beta_i$ for $i=2,4$, we get the following:
\begin{equation*}
 \MC(\xi)\cong_{\xi \to 0, \xi \in \mathbb{H}_{-m}} \frac{-1}{\exp \bigl(\frac{1}{\xi}m \operatorname{Log}_p\bigl(-\frac{m}{2e}\bigr)-{\rm i}m^{(3)}\operatorname{Arg}_p(-m)-\int_{-\infty}^{\infty}\gamma^*A_{11}\bigr)}=A(\xi).
\end{equation*}

Hence, the asymptotics agree on $\mathbb{H}_m$ and $\mathbb{H}_{-m}$, and we proved what we want.
\end{proof}
\end{subsubsection}
\begin{subsubsection}[The reality condition and the asymptotics as xi to infty]{The reality condition and the asymptotics as $\boldsymbol{\xi \to \infty}$}

In this section, we prove the reality condition for the magnetic twistor coordinate (recall \eqref{rc}). Once we have the reality condition, the asymptotics when $\xi \to \infty$ automatically follow from the asymptotics when $\xi \to 0$.

The reality condition states that the following equality must hold:
\begin{equation*}
\mathcal{X}_m(\xi)=\overline{\mathcal{X}_m\bigl(-1/\overline{\xi}\bigr)^{-1}} \qquad \text{for} \quad \xi \in \mathbb{C
 }^*/l_{\pm}(-2{\rm i}m).
\end{equation*}

Let $\bigl(E,\overline{\partial}_E, \theta, h,g\bigr) \in \mathcal{H}^{{\rm fr}}$ be a compatibly framed wild harmonic bundle. We will denote by $(e_1,e_2)$ an extension of $g$ to a ${\rm SU}(2)$ frame in a neighborhood $U_{\infty}$ of $z=\infty$, where the singularity of the associated flat connection $\nabla^{\xi}=\xi^{-1}\theta + D\bigl(\overline{\partial}_E,h\bigr) + \xi \theta^{\dagger}$ has the appropriate form. We write in this frame $\nabla^{\xi}={\rm d}+A(\xi)$.

We can consider the conjugate bundle \smash{$\overline{E}\to \overline{\CP}$}, with the induced connection \smash{$\overline{\nabla^{\xi}}$}. In the induced frame $\{\overline{e}_i\}$ the connection has the form \smash{$\overline{\nabla^{\xi}}={\rm d}+\overline{A(\xi)}$}. Furthermore, we can consider the dual bundle \smash{$\overline{E}^{*}\to \overline{\CP}$} with the induced connection \smash{$\overline{\nabla^{\xi}}^*$}. In the dual frame $\{\overline{e}^*_i\}$, the connection \smash{$\overline{\nabla^{\xi}}^{*}$} takes the form \smash{${\rm d}-\overline{A(\xi)}^t={\rm d}+A\bigl(-1/\overline{\xi}\bigr)$}.

Consider now the associated compatibly framed meromorphic flat bundle $\bigl(\mathcal{PE}_a^{\xi},\nabla^{\xi},\tau_a^{\xi}\bigr)$ and a fundamental solution $(y_1,y_2)=\Phi_i(\xi)$ of $\nabla^{\xi}$ on \smash{$\widehat{\operatorname{Sect}}_i$}, with the corresponding asymptotics determined by the compatible frame. If we denote by $e^{Q(\xi)}=\operatorname{diag}\bigl(e^{Q_1(\xi)},e^{Q_2(\xi)}\bigr)$, where
\begin{gather*}
 e^{Q_1(\xi)}=\exp \left(-\xi^{-1}\left(\frac{1}{2w^2}-m\operatorname{Log}(w)\right)-{\rm i}m^{(3)}\operatorname{Arg}(w)
 -\xi\left(\frac{1}{2\overline{w}^2}-\overline{m}\overline{\operatorname{Log}(w)}\right)\right),\\
 e^{Q_2(\xi)}=\exp \left(\xi^{-1}\left(\frac{1}{2w^2}-m\operatorname{Log}(w)\right)+{\rm i}m^{(3)}\operatorname{Arg}(w)
 +\xi\left(\frac{1}{2\overline{w}^2}-\overline{m}\overline{\operatorname{Log}(w)}\right)\right),
\end{gather*}
then we know that by Theorem \ref{twistasymp} that $(y_1,y_2)\cdot e^{-Q(\xi)} \to (e_1,e_2)$ as $w\to 0$ along $\widehat{\operatorname{Sect}}_i$. Here, the $Q_i(\xi)$ are defined using the same branch of the logarithm and argument as the one used to define the flat frames $\Phi_i(\xi)$ in Section \ref{deftwistcoords}. We then have the following lemma.

\begin{Lemma}\label{dual asympotics}
If $(e_1,e_2)\cdot B=\Phi_i=(y_1,y_2)$, then \smash{$(\overline{e}_1^*,\overline{e}_2^*)\cdot \bigl(B^{-1}\bigr)^{\dagger}:=\overline{\Phi}_i^{*}= (\overline{y}_1^{*},\overline{y}_2^{*})$} is a flat frame for \smash{$\bigl(\overline{\mathcal{PE}_a^{\xi}}^*,\overline{\nabla^{\xi}}^*\bigr)$}. Furthermore, \smash{$(\overline{y}_1^{*},\overline{y}_2^{*})\cdot e^{-Q(-1/\overline{\xi})}\to (\overline{e}_1^*,\overline{e}_2^*)$} as $w\to 0$ along \smash{$\widehat{\textnormal{Sect}}_i(\xi)$}.
\end{Lemma}

\begin{proof}
Since $\Phi_i$ is a flat frame, and in the frame $(e_1,e_2)$ we have that $\nabla^{\xi}={\rm d}+A(\xi)$, then~$B$ must satisfy
$
 B^{-1}{\rm d}B+B^{-1}AB=0$.
It then follows that
\smash{$
 -B^{\dagger}A^{\dagger}\bigl(B^{-1}\bigr)^{\dagger}+B^{\dagger}{\rm d}\bigl(B^{-1}\bigr)^{\dagger}=0$}.
Since in the frame $(\overline{e}_1^*,\overline{e}_2^*)$, we have that \smash{$\overline{\nabla^{\xi}}^*={\rm d}-A^{\dagger}$} and conclude that \smash{$(\overline{e}_1^*,\overline{e}_2^*)\cdot \bigl(B^{-1}\bigr)^{\dagger}$} must be a~flat frame for \smash{$\overline{\nabla^{\xi}}^*$}.

To check the asymptotics, just notice that $(y_1,y_2)\cdot e^{-Q(\xi)} \to (e_1,e_2)$ implies that $Be^{-Q(\xi)}\to 1$. On the other hand, we have that
\begin{align*}
 (\overline{y}_1^{*},\overline{y}_2^{*})\cdot e^{-Q(-1/\overline{\xi})}&=(\overline{e}_1^*,\overline{e}_2^*)\cdot \bigl(B^{-1}\bigr)^{\dagger}e^{-Q(-1/\overline{\xi})}=(\overline{e}_1^*,\overline{e}_2^*)\cdot \bigl(B^{-1}\bigr)^{\dagger}e^{\overline{Q(\xi)}}\\
 &=(\overline{e}_1^*,\overline{e}_2^*)\cdot \bigl(e^{Q(\xi)}B^{-1}\bigr)^{\dagger}=(\overline{e}_1^*,\overline{e}_2^*)\cdot \bigl(\bigl(Be^{-Q(\xi)}\bigr)^{-1}\bigr)^{\dagger},
\end{align*}
so we conclude that \smash{$(\overline{y}_1^{*},\overline{y}_2^{*})\cdot e^{-Q(-1/\overline{\xi})}\to (\overline{e}_1^*,\overline{e}_2^*)$} as $w\to 0$ along \smash{$\widehat{\operatorname{Sect}}_i(\xi)$}.
\end{proof}

 We will need two more easy lemmas in order to show the reality condition of the magnetic coordinate.

 \begin{Lemma}
 If on \smash{$\widehat{\textnormal{Sect}}_i\!\cap\! \widehat{\textnormal{Sect}}_{i+1}$}, we have $\Phi_{i+1}\!=\!\Phi_i \cdot S_i$, then we have that \smash{$\overline{\Phi}_{i+1}^{*}\!=\!\overline{\Phi}_i^{*}\cdot \bigl(S_i^{-1}\bigr)^{\dagger}$}.
 \end{Lemma}

 \begin{proof}
 If we write $(e_1,e_2)\cdot B_i=\Phi_i$ and $(e_1,e_2)\cdot B_{i+1}=\Phi_{i+1}$, then we have that $B_iS_i=B_{i+1}$. Hence, we get that \smash{$\bigl(B_i^{-1}\bigr)^{\dagger}\bigl(S_i^{-1}\bigr)^{\dagger}=\bigl(B_{i+1}^{-1}\bigr)^{\dagger}$}, and then by the previous lemma
 \begin{equation*}
 \overline{\Phi}_{i+1}^{*}=(\overline{e}_1^*,\overline{e}_2^*)\cdot \bigl(B_{i+1}^{-1}\bigr)^{\dagger}=(\overline{e}_1^*,\overline{e}_2^*)\cdot \bigl(B_i^{-1}\bigr)^{\dagger}\bigl(S_i^{-1}\bigr)^{\dagger}=\overline{\Phi}_{i}^{*}\cdot\bigl(S_i^{-1}\bigr)^{\dagger}.\tag*{\qed}
 \end{equation*}\renewcommand{\qed}{}
 \end{proof}

\begin{Lemma}
The local bundle map defined by $e_i\to \overline{e_i}^{*}$, gives an isomorphism $\smash{\bigl(\overline{\mathcal{PE}_a^{\xi}}^*,\overline{\nabla^{\xi}}^*\bigr)}|_{U_\infty}\allowbreak\cong \smash{\bigl(\mathcal{PE}_a^{-\frac{1}{\overline{\xi}}},\nabla^{-\frac{1}{\overline{\xi}}}\bigr)}|_{U_\infty}$, where $U_\infty$ is a small neighborhood of $z=\infty$ where the frame $(e_1,e_2)$ is defined. In particular, we have the following correspondence between flat frames: for $i=1,2$ and~${\xi \in \mathbb{H}_m}$, the flat frame $\overline{\Phi}_i^{*}(\xi)$ goes to the flat frame $\Phi_{i+1}\bigl(-1/\overline{\xi}\bigr)$ of $\nabla^{-1/\overline{\xi}}$ on $\widehat{\textnormal{Sect}}_{i+1}\bigl(-1/\overline{\xi}\bigr)$; and for~${i=2,3}$ and $\xi \in \mathbb{H}_{-m}$, to the flat frame $\Phi_{i-1}\bigl(-1/\overline{\xi}\bigr)$ of \smash{$\nabla^{-1/\overline{\xi}}$} on \smash{$\widehat{\textnormal{Sect}}_{i-1}\bigl(-1/\overline{\xi}\bigr)$}.
\end{Lemma}

\begin{proof}

For the first statement, just notice that in the frame $(\overline{e}_1^{*},\overline{e}_2^*)$, the induced connection looks like \smash{$\overline{\nabla^{\xi}}^*={\rm d}+A\bigl(-1/\overline{\xi}\bigr)$}. On the other hand, the flat frame $\Phi_{i}\bigl(-1/\overline{\xi}\bigr)$ on \smash{$\widehat{\operatorname{Sect}}_{i}\bigl(-1/\overline{\xi}\bigr)$} is uniquely characterized by the asymptotic condition \smash{$\Phi_{i}\bigl(-1/\overline{\xi}\bigr)e^{-Q(-1/\overline{\xi})}\to (e_1,e_2)$} when~${w\to 0}$ along \smash{$\widehat{\operatorname{Sect}}_{i}\bigl(-1/\overline{\xi}\bigr)$}. If $i=1,2$ and $\xi \in \mathbb{H}_m$, we have that \smash{$\widehat{\operatorname{Sect}_i}(\xi)=\widehat{\operatorname{Sect}}_{i+1}\bigl(-1/\overline{\xi}\bigr)$} and that the flat frame $\overline{\Phi}_i^{*}$ of \smash{$\overline{\nabla^{\xi}}^*$} goes to a flat frame of \smash{$\nabla^{-\frac{1}{\overline{\xi}}}$} satisfying the corresponding asymptotic condition on \smash{$\widehat{\operatorname{Sect}}_{i+1}\bigl(-1/\overline{\xi}\bigr)$}. We then conclude that for $i=1,2$ and $\xi \in \mathbb{H}_m$ we have $\overline{\Phi}^*_i(\xi)=\Phi_{i+1}\bigl(-1/\overline{\xi}\bigr)$. The other case similarly follows.
\end{proof}

\begin{Theorem}\label{rct}
The magnetic twistor coordinate satisfies the following reality condition:
\begin{equation*}
 \mathcal{X}_m(\xi)=\overline{\mathcal{X}_m\bigl(-1/\overline{\xi}\bigr)^{-1}} \qquad \text{for} \quad \xi \in \mathbb{C
 }^*/l_{\pm}(-2{\rm i}m).
\end{equation*}
\end{Theorem}

\begin{proof}
Assume first that $\xi \in \mathbb{H}_m$. We then want to relate $b\bigl(-1/\overline{\xi}\bigr)$ and $a(\xi)$.

By the previous two lemmas, we have that
\begin{align*}
\Phi_2\bigl(-1/\overline{\xi}\bigr)\cdot S_2\bigl(-1/\overline{\xi}\bigr)&=\Phi_3\bigl(-1/\overline{\xi}\bigr)=\overline{\Phi}^*_2(\xi)=\overline{\Phi}^*_1(\xi)\cdot \bigl(S_1(\xi)^{-1}\bigr)^{\dagger}\\
&=\Phi_2\bigl(-1/\overline{\xi}\bigr)\cdot\bigl(S_1^{-1}(\xi)\bigr)^{\dagger}.
\end{align*}

From this, we conclude that \smash{$S_2\bigl(-1/\overline{\xi}\bigr)=\bigl(S_1^{-1}(\xi)\bigr)^{\dagger}$}, so that the Stokes matrix elements are related by \smash{$b\bigl(-1/\overline{\xi}\bigr)=-\overline{a(\xi)}$}. Hence, we have that for $\xi \in \mathbb{H}_{m}$
\begin{equation*}
 \mathcal{X}_m(\xi)=a(\xi)=-\overline{b\bigl(-1/\overline{\xi}\bigr)}=\overline{\mathcal{X}_m\bigl(-1/\overline{\xi}\bigr)^{-1}}.
\end{equation*}

Similarly, if $\xi \in \mathbb{H}_{-m}$, we need to compare $b(\xi)$ and $a\bigl(-1/\overline{\xi}\bigr)$. By the above two lemmas, we have that
\begin{align*}
 \Phi_1\bigl(-1/\overline{\xi}\bigr)\cdot S_1\bigl(-1/\overline{\xi}\bigr)&=\Phi_2\bigl(-1/\overline{\xi}\bigr)=\overline{\Phi}^*_3(\xi)=\overline{\Phi}^*_2(\xi)\cdot \bigl(S_2(\xi)^{-1}\bigr)^{\dagger}\\
 &=\Phi_1\bigl(-1/\overline{\xi}\bigr)\cdot\bigl(S_2^{-1}(\xi)\bigr)^{\dagger},
\end{align*}
so that $S_1\bigl(-1/\overline{\xi}\bigr)=S_2^{-1}(\xi)^{\dagger}$ and hence the relation among the Stokes matrix elements is \smash{$a\bigl(-1/\overline{\xi}\bigr)=-\overline{b(\xi)}$}. We then conclude that for $\xi \in \mathbb{H}_{-m}$
\begin{equation*}
 \mathcal{X}_m(\xi)=-\frac{1}{b(\xi)}=\frac{1}{\overline{a\bigl(-1/\overline{\xi}\bigr)}}=\overline{\mathcal{X}_m\bigl(-1/\overline{\xi}\bigr)^{-1}}.
\end{equation*}

This proves the reality condition.
\end{proof}

\begin{Corollary} \label{rcc}Using the same notation from Section {\rm\ref{asymp0}}, we have that
\begin{equation*}
 \MC(\xi)\sim_{\xi \to \infty} \overline{A\bigl(-1/\overline{\xi}\bigr)^{-1}}.
\end{equation*}

\end{Corollary}

\begin{proof}
By the reality condition, we have that the asymptotics of $\MC(\xi)$ as $\xi \to \infty$ match the asymptotics of \smash{$\overline{\mathcal{X}_m\bigl(-1/\overline{\xi}\bigr)^{-1}}$} as $\xi \to \infty$. But by Proposition \ref{a0}, we have that
\begin{equation*}
\overline{\mathcal{X}_m\bigl(-1/\overline{\xi}\bigr)^{-1}}\sim_{\xi \to \infty} \overline{A\bigl(-1/\overline{\xi}\bigr)^{-1}},
\end{equation*}
so the result follows.
\end{proof}

\end{subsubsection}

\subsection{Non-vanishing of Stokes data}\label{nonvanishing}

In this section, we show that given $\bigl(E,\overline{\partial}_E,\theta,h,g\bigr) \in \mathcal{H}^{{\rm fr}}$ with $m\neq 0$, the coordinate $\MC\bigl(\bigl(E,\overline{\partial}_E,\allowbreak\theta, h,g\bigr),\xi\bigr)$ is actually well defined for $\xi \in \mathbb{C}^* -l_{\pm}(-2{\rm i}m)$ (i.e., that $b(\xi)\neq 0$ for $\xi \in \mathbb{H}_{-m}$).

\begin{subsubsection}{The case of trivial Stokes data}

We begin with the following lemma dealing with the case of trivial Stokes data.

\begin{Lemma}\label{m=0} Consider \smash{$\bigl(E,\overline{\partial}_E,\theta,h,g\bigr) \in \mathcal{H}^{\textnormal{fr}}$} and its associated compatibly framed filtered flat bundle $\bigl(\mathcal{P}^{h}_*\mathcal{E}^{\xi},\nabla^{\xi},\tau_*^{\xi}\bigr)$ for $\xi \in \mathbb{C}^*$. Suppose that the Stokes data associated to $\bigl(\mathcal{P}^{h}_*\mathcal{E}^{\xi},\nabla^{\xi},\tau_*^{\xi}\bigr)$ is trivial $($i.e., $S_i=1$ for $i=1,2,3,4$ and $M_0=1)$. Then the parameter $m$ specifying the singularity of the Higgs field of $\bigl(E,\overline{\partial}_E,\theta,h,g\bigr)$ must be $0$.
\end{Lemma}

\begin{proof}We divide the proof into two cases:
\begin{itemize}\itemsep=0pt
 \item Suppose first that the parabolic weights of \smash{$\bigl(\mathcal{P}^{h}_{1/2}\mathcal{E}^{\xi},\nabla^{\xi},\tau_{1/2}^{\xi}\bigr)$} lie in $\bigl(-\frac12,\frac12\bigr)\subset \mathbb{R}$. Con\-sider~the framed wild harmonic bundle $\bigl(E_0,\overline{\partial}_{E_0},\theta_0,h_0,g_0\bigr) \in \mathcal{H}^{{\rm fr}}$ given in Example \ref{trivWHB}. Let us give a description of the associated filtered bundle \smash{$\bigl(\mathcal{P}^h_*\mathcal{E}^{\xi}_0,\nabla^{\xi}_0,\tau_{*,0}^{\xi}\bigr)$}. To give a description of this filtered bundle, it is enough to specify the parabolic flat bundle \smash{$\bigl(\mathcal{P}^{h}_{1/2}\mathcal{E}^{\xi}_0,\nabla^{\xi}_0, \tau_{1/2,0}^{\xi}\bigr)$}. Since $m^{(3)}=m=0$, the parabolic filtration is trivial in this case, with the holomorphic frames describing the extension of $\mathcal{E}^{\xi}$ to $\mathcal{PE}_0^{\xi}$ given by
 \[
 \tau_0^{\xi}=g_0\cdot \exp \left(\left(\frac{\overline{\xi}}{2w^2}-\frac{\xi}{2\overline{w}^{2}}\right)H\right).
 \]
 Furthermore, \smash{$g_0\cdot \exp \bigl(\bigl(\frac{z^2\overline{\xi}}{2}-\frac{\overline{z}^{2}\xi}{2}\bigr)H\bigr)$} gives a global holomorphic trivialization of \smash{$\mathcal{P}^{h}_{1/2}\mathcal{E}_0^{\xi}$}, so that \smash{$\mathcal{P}^{h}_{1/2}\mathcal{E}_0^{\xi}\cong \mathcal{O}\oplus \mathcal{O}$}. In this global trivialization, the connection has the following form:
\begin{equation*}
 \nabla^{\xi}_0={\rm d}-\bigl(\xi^{-1}+\overline{\xi}\bigr)H\frac{{\rm d}w}{w^3}.
\end{equation*}

This connection has a frame of flat sections defined over $\mathbb{C} \subset \CP$, given by $g_0\cdot \exp \smash{\bigl(\bigl(-\frac{z^2\xi^{-1}}{2}}\allowbreak-\smash{\frac{\overline{z}^{2}\xi}{2}\bigr)H\bigr)}$ satisfying the appropriate asymptotics with respect to $\tau_0^{\xi}$ (and $g_0$). From this fact, we see that the formal monodromy and the associated Stokes matrices are trivial.
Now going back to \smash{$\bigl(\mathcal{P}^{h}_{1/2}\mathcal{E}^{\xi},\nabla^{\xi},\tau_{1/2}^{\xi}\bigr)$}, the fact that it comes from a harmonic bundle implies that \smash{$\operatorname{pdeg}\bigl(\mathcal{P}^{h}_{1/2}\mathcal{E}^{\xi}\bigr)=0$}, while the assumption on the parabolic weights imply that \smash{$\operatorname{deg}\bigl(\mathcal{P}^{h}_{1/2}\mathcal{E}^{\xi}\bigr)=0$}. Since it also has trivial Stokes data, the Riemann--Hilbert correspondence for framed flat bundles presented in $\cite{BB04}$ implies that \smash{$\bigl(\mathcal{P}^{h}_{1/2}\mathcal{E}^{\xi},\nabla^{\xi},\tau_{1/2}^{\xi}\bigr)$} is isomorphic to \smash{$\bigl(\mathcal{P}^{h}_{1/2}\mathcal{E}^{\xi}_0,\nabla^{\xi}_0, \tau_{1/2,0}^{\xi}\bigr)$} (as compatibly framed flat bundles). In particular, \smash{$\bigl(\mathcal{P}^{h}_{1/2}\mathcal{E}^{\xi},\nabla^{\xi}\bigr)$} is isomorphic with \smash{$\bigl(\mathcal{P}^{h}_{1/2}\mathcal{E}^{\xi}_0,\nabla^{\xi}_0\bigr)$}, so that \smash{$\bigl(E,\overline{\partial}_E,\theta,h\bigr)$} is isomorphic to $\bigl(E_0,\overline{\partial}_{E_0},\theta_0,h_0\bigr)$. Hence,
\begin{equation*}
 -\bigl(z^2+2m\bigr){\rm d}z^2=\operatorname{Det}(\theta)=\operatorname{Det}(\theta_0)=-z^2{\rm d}z^2,
\end{equation*} which implies that $m=0$.
\item The parabolic weights of \smash{$\bigl(\mathcal{P}^{h}_{1/2}\mathcal{E}^{\xi},\nabla^{\xi},\tau_{1/2}^{\xi}\bigr)$} are both equal to $1/2$: in this case, the condition on the parabolic weights and the fact that \smash{$\operatorname{pdeg}\bigl(\mathcal{P}^{h}_{1/2}\mathcal{E}^{\xi}\bigr)=0$} implies that \smash{$\operatorname{deg}\bigl(\mathcal{P}^{h}_{1/2}\mathcal{E}^{\xi}\bigr)=1$}. Hence, the exponent of formal monodromy $\Lambda(\xi)$ of \smash{$\bigl(\mathcal{P}^{h}_{1/2}\mathcal{E}^{\xi},\nabla^{\xi},\tau_{1/2}^{\xi}\bigr)$} must satisfy $\operatorname{Tr}(\Lambda(\xi))\allowbreak=-1$ (see \cite{B}), which together with the fact that the formal monodromy is trivial implies that
\begin{equation*}
 \Lambda(\xi)=\begin{bmatrix} n & 0 \\
 0 & -n-1\\ \end{bmatrix}
\end{equation*}
for some $n \in \mathbb{Z}$.
Now consider the degree $1$ bundle $\mathcal{O}(-n)\oplus \mathcal{O}(n+1)\to \CP$. In the $z$ coordinate of~${\mathbb{C}\subset \CP}$, and the usual trivialization over that neighborhood, we consider the connection
$
 \nabla= {\rm d} +\bigl(\xi^{-1}+\overline{\xi}\bigr)Hz{\rm d}z$.
In the neighborhood $\CP \setminus\{0\}$ with coordinate $w=\frac{1}{z}$ and its usual trivialization, we have that
\begin{equation*}
 \nabla= {\rm d} -\bigl(\xi^{-1}+\overline{\xi}\bigr)H\frac{{\rm d}w}{w^3} +\begin{bmatrix} n & 0\\
 0 & -n-1 \\ \end{bmatrix}\frac{{\rm d}w}{w}.
\end{equation*}
So $\nabla$ is a meromorphic connection with the appropriate irregular part. The trivialization over $\CP\setminus\{0\}$ clearly gives a compatible framing $\tau$ over $z=\infty$, so we get a compatibly framed meromorphic flat bundle $(\mathcal{O}(-n)\oplus \mathcal{O}(n+1),\nabla,\tau)$.
It is easy to check that the~Stokes matrices and formal monodromy associated to $(\mathcal{O}(-n)\oplus \mathcal{O}(n+1),\nabla,\tau)$ are trivial. Hence, since $(\mathcal{O}(-n)\oplus \mathcal{O}(n+1),\nabla,\tau)$ and \smash{$\bigl(\mathcal{P}^{h}_{1/2}\mathcal{E}^{\xi},\nabla^{\xi},\tau_{1/2}^{\xi}\bigr)$} have the same formal type~and~same Stokes data, we get that \smash{$\bigl(\mathcal{P}^{h}_{1/2}\mathcal{E}^{\xi},\nabla^{\xi},\tau_{1/2}^{\xi}\bigr)$} is isomorphic to $(\mathcal{O}(-n)\oplus \mathcal{O}({n+1}),\allowbreak\nabla,\tau)$ as framed flat bundles (the proof of this fact follows part of the argument of the proof of Lemma \ref{secondlemma}, for example).
On the other hand, \smash{$\bigl(\mathcal{P}^{h}_{1/2}\mathcal{E}^{\xi},\nabla^{\xi},\tau_{1/2}^{\xi}\bigr)$} has the trivial filtration as a parabolic bundle with parabolic weights both equal to $\frac12$, so $(\mathcal{O}(-n)\oplus \mathcal{O}(n+1),\nabla,\tau)$ also gets this parabolic structure. Now notice that $\nabla$ preserves the line bundles $\mathcal{O}(-n)$ and $\mathcal{O}(n+1)$, and with the induced parabolic structures on the line bundles we have that
\begin{gather*}
 \operatorname{pdeg}(\mathcal{O}(-n))=-n-\frac{1}{2},\qquad
 \operatorname{pdeg}(\mathcal{O}(n+1))=n+1-\frac{1}{2}=n+\frac{1}{2}.
\end{gather*}
This shows that no matter what $n\in \mathbb{Z}$ is, we have that \smash{$\bigl(\mathcal{P}^{h}_{1/2}\mathcal{E}^{\xi},\nabla^{\xi}\bigr)$} is unstable. But \smash{$\bigl(\mathcal{P}^{h}_{1/2}\mathcal{E}^{\xi},\nabla^{\xi}\bigr)$} comes from a harmonic bundle, so it must also be polystable by Theorem~\ref{HMFB}. This contradiction shows that \smash{$\bigl(\mathcal{P}^{h}_{1/2}\mathcal{E}^{\xi},\nabla^{\xi},\tau_{1/2}^{\xi}\bigr)$} cannot have parabolic weights equal to~$1/2$ if it has trivial monodromy data. We conclude that the only case that occurs is the previous case.
\end{itemize}

Hence, we conclude that if $\bigl(E,\overline{\partial}_E,\theta,h,g\bigr)$ has trivial Stokes data, we must necessarily have $m=0$.
\end{proof}

\end{subsubsection}

\begin{subsubsection}[Proof of the non-vanishing of Stokes data when m neq 0]{Proof of the non-vanishing of Stokes data when $\boldsymbol{m\neq 0}$}

Recall that our Stokes data is made out of the $2\times 2$ unipotent Stokes matrices $S_i(\xi)$ for $i=1,2,3,4$ and formal monodromy $M_0={\rm e}^{-2\pi {\rm i} \Lambda(\xi)}$ that must satisfy the relation $S_1S_2S_3S_4M_0^{-1}=1$. We will label the off-diagonal non-trivial complex numbers of $S_1$, $S_2$, $S_3$ and $S_4$ by $a$, $b$, $c$ and~$d$, respectively, as in Section \ref{deftwistcoords}. Hence, $\MC(\xi)$ is defined using $a(\xi)$ and $-1/b(\xi)$ using the conventions from Section \ref{deftwistcoords}.

Because of the reality condition satisfied by the Stokes data, we have that for $\xi \in \mathbb{H}_{-m}$ the equality \smash{$-\overline{a\bigl(-1/\overline{\xi}\bigr)}=b(\xi)$} holds. Hence, if we want to show that $\MC(\xi)$ is well defined, it is enough to show that $a(\xi)$ does not vanish for $\xi \in \mathbb{H}_m$. To show this, we will need the following lemma, whose proof is easy:

\begin{Lemma}\label{nvm} Let $U_{\pm}$ be the set of upper $($lower, respectively$)$ unipotent $2\times 2$ matrices, and $T\subset {\rm SL}(2,\mathbb{C})$ the subset of diagonal matrices. Then the set
\begin{equation*}
 \mathcal{M}=\bigl\{(S_1,S_2,S_3,S_4,M_0)\in (U_{-}\times U_{+})^2\times T \mid S_1S_2S_3S_4M_0^{-1}=1\bigr\}
\end{equation*}
is a complex $2$-dimensional manifold. Furthermore, if we denote by $\mathcal{S}$ the subset of $\mathcal{M}$ defined~by
\begin{equation*}
 \mathcal{S}=\{(S_1,S_2,S_3,S_4,M_0) \in \mathcal{M} \mid S_1\neq S_2 \ \text{and either}\ S_1=1 \ \text{or}\ S_2=1\},
\end{equation*}
then $\mathcal{S}$ is a complex $1$-dimensional submanifold with $2$ components. The components are determined by whether $S_1=1$ $($i.e., $a=0)$ or $S_2=1$ $($i.e., $b=0)$.
\end{Lemma}

\begin{Proposition} \label{nonvan}Let $\bigl(E,\overline{\partial}_E,\theta, h, g\bigr)\in \mathcal{H}^{\textnormal{fr}}$ with associated parameters $m^{(3)}$ and $m\neq 0$. Then for $\xi \in \mathbb{H}_m$ we have that $a\bigl(\bigl[\bigl(E,\overline{\partial}_E,\theta, h, g\bigr)\bigr],\xi\bigr)\neq 0$.
\end{Proposition}

\begin{proof}
Let us first see what condition we get if $a(\xi)=0$ for $\xi \in \mathbb{H}_m$. If $a(\xi)=0$, then by the relation $1+a(\xi)b(\xi)=\mu^{-1}(\xi)$ from \eqref{Stokes relations} we conclude that $\mu^{-1}=1$. In particular, by \eqref{defetc} we must have
\smash{$
 \xi^{-1}m -m^{(3)}-\overline{m}\xi=k
$}
for some $k \in \mathbb{Z}$.

Solving for $\xi$, we get two solutions $\xi^{\pm}_k$ for each $k \in \mathbb{Z}$ given by
\begin{equation*}
 \xi_k^{\pm}=\left(\frac{-\bigl(k+m^{(3)}\bigr)\pm \sqrt{\bigl(k+m^{(3)}\bigr)^2+4|m|^2}}{2|m|^2}\right)m.
\end{equation*}
The values $\xi_k^{\pm}$ are the only possible values of the twistor parameter for which we could have $a(\xi)=0$ \big(for our fixed $\bigl(E,\overline{\partial}_E,\theta, h, g\bigr) \in \mathcal{H}^{{\rm fr}}$\big). In particular, the solutions contained in $\mathbb{H}_m$ are given by
\begin{equation*}
 \xi_k^{+}=\left(\frac{-\bigl(k+m^{(3)}\bigr)+ \sqrt{\bigl(m^{(3)}+k\bigr)^2+4|m|^2}}{2|m|^2}\right)m,
\end{equation*}
and from the formula we see that $\xi_k^{+} \to 0$ as $k\to \infty$, and that $\xi_k^{+} \to \infty$ as $k \to -\infty$. Furthermore, all $\xi_k^{+}$ are contained in the ray determined by $0$ and $m$.

Because of the asymptotic behavior of $a(\xi)$ as $\xi \to 0$ from Section \ref{asymp0}, we see that $a\bigl(\xi_k^{+}\bigr)$ cannot be $0$ for all $k$, otherwise the asymptotics would not hold. Hence, $a\bigl(\xi_k^{+}\bigr)\neq 0$ for some $k$ sufficiently big. Similarly, we get from the asymptotics as $\xi \to \infty$ that $a\bigl(\xi_k^{+}\bigr)\neq 0$ for some $k$ sufficiently negative.

We will now show that the fact that $a\bigl(\xi_k^{+}\bigr)\neq 0$ for at least one $k \in \mathbb{Z}$ implies that $a(\xi_n^{+})\neq 0$ for all $n\in \mathbb{Z}$. To show this, we will need the following lemma:

\begin{Lemma} Let $\bigl(\hspace{-0.3pt}E,\overline{\partial}_E,\theta, h, g\bigr)\!\in\! \mathcal{H}^{\textnormal{fr}}$ with associated parameters $m^{(3)}$ and $m\neq 0$, and~\mbox{consi\-der} the curve \smash{$\xi\bigl(\widetilde{m^{(3)}}\bigr)\colon\mathbb{R} \to \mathbb{H}_m$} given by
\begin{equation*}
 \xi\bigl(\widetilde{m^{(3)}}\bigr)=\left(\frac{-\bigl(k+\widetilde{m^{(3)}}\bigr)+ \sqrt{\bigl(\widetilde{m^{(3)}}+k\bigr)^2+4|m|^2}}{2|m|^2}\right)m.
\end{equation*}
Furthermore, let $\bigl(\mathcal{P}_*^{h}\mathcal{E}^{\xi},\nabla^{\xi},\tau_*^{\xi}\bigr)$ denote the associated compatibly framed flat bundle of $\bigl(E,\overline{\partial}_E,\allowbreak\theta, h, g\bigr)$ for $\xi \in \mathbb{C}^*$, and $\mathcal{S}$, $\mathcal{M}$ be as in Lemma {\rm\ref{nvm}}. Then we can find a continuous map $\tau\colon \mathbb{R} \to \mathcal{S}\subset \mathcal{M}$ such that \smash{$\tau\bigl(m^{(3)}+n\bigr)$} is the Stokes data associated \smash{$\bigl(\mathcal{P}_*^{h}\mathcal{E}^{\xi(m^{(3)}+n)},\nabla^{\xi(m^{(3)}+n)},\tau_*^{\xi(m^{(3)}+n)}\bigr)$} for all $n \in \mathbb{Z}$.

\end{Lemma}

\begin{proof}[of lemma:] Let $\bigl(\mathcal{P}_*^{h}\mathcal{E}^0,\theta\bigr)$ be the associated filtered Higgs bundle to $\bigl(E,\overline{\partial}_E,\theta,h\bigr)$, and let $(E,\theta)=\bigl(\mathcal{P}_*^{h}\mathcal{E}^0,\theta\bigr)|_{\CP \setminus \{\infty\}}$. Near $z=\infty$, we choose a splitting of $(E,\theta)$ in eigenlines, such that $\theta$ has the form \eqref{singHiggs}. We then have that \smash{$m^{(3)}$} specifies the filtered/parabolic structure of~$\bigl(\mathcal{P}_*^{h}\mathcal{E}^0,\theta\bigr)$ (see, for example, the construction of Lemma \ref{C1}). By varying \smash{$m^{(3)}$} in $\mathbb{R}$, we thus vary the parabolic structure and hence the filtered Higgs bundle of parabolic degree $0$ we obtain from $(E,\theta)$ by following Lemma \ref{C1}. We will denote by \smash{$\bigl(\mathcal{P}_*^{h}\mathcal{E}^0\bigl(\widetilde{m^{(3)}}\bigr),\theta\bigr)$} the filtered Higgs bundles of parabolic degree $0$ that we obtain from $(E,\theta)$ with parabolic structure specified~by \smash{$\widetilde{m^{(3)}}\in \mathbb{R}$}.

Notice that the filtered/parabolic structure only depends on \smash{$\widetilde{m^{(3)}}$} mod $1$, so if \smash{$h\bigl(\widetilde{m^{3}}\bigr)$} denotes a curve of harmonic metrics adapted to \smash{$\bigl(\mathcal{P}_*^{h}\mathcal{E}^0\bigl(\widetilde{m^{(3)}}\bigr),\theta\bigr)$}, we get by the uniqueness of harmonic metrics that for all $n \in \mathbb{Z}$, we have \smash{$h\bigl(\widetilde{m^{(3)}}+n\bigr)=c(n)h\bigl(\widetilde{m^{(3)}}\bigr)$} for some constant $c(n)>0$.

We denote by \smash{$\bigl(E,\overline{\partial}_E,\theta,h\bigl(\widetilde{m^{(3)}}\bigr)\bigr)$} the corresponding curve of wild harmonic bundles, and by
\begin{equation*}
 \sigma\bigl(\widetilde{m^{(3)}}\bigr)=\bigl[\bigl(E\bigl(\widetilde{m^{(3)}}\bigr),\overline{\partial}_{E(\widetilde{m^{(3)}})},\theta,h\bigl(\widetilde{m^{(3)}}\bigr),g\bigl(\widetilde{m^{(3)}}\bigr)\bigr)\bigr]
\end{equation*}
the corresponding curve of equivalence classes of framed wild harmonic bundles obtained by applying the construction of Proposition \ref{constructionframe}. Notice that \smash{$g\bigl(\widetilde{m^{(3)}}\bigr)$} satisfies that $\smash{g\bigl(\widetilde{m^{(3)}}+n\bigr)}=\smash{g\bigl(\widetilde{m^{(3)}}\bigr)\cdot \sqrt{c(n)^{-1}}}$ and that \[\bigl(E\bigl(\widetilde{m^{(3)}}\bigr),\overline{\partial}_{E(\widetilde{m^{(3)}})}\bigr)
=\bigl(E\bigl(\widetilde{m^{(3)}}+n\bigr),\overline{\partial}_{E(\widetilde{m^{(3)}}+n)}\bigr)
\]
 for all $n \in \mathbb{Z}$. Hence, by Example \ref{exmorph} we have that \smash{$\sigma\bigl(\widetilde{m^{(3)}}\bigr)=\sigma\bigl(\widetilde{m^{(3)}}+n\bigr)$} for all $n \in \mathbb{Z}$.

Now consider the continuous map \smash{$\tau\bigl(\widetilde{m^{3}}\bigr)\colon \mathbb{R} \to \mathcal{M}$}, where \smash{$\tau\bigl(\widetilde{m^{(3)}}\bigr)$} is the Stokes data of the associated compatibly framed flat bundle
\[
\bigl[\bigl(\mathcal{P}^{h(\widetilde{m^{(3)}})}_*\mathcal{E}^{\xi(\widetilde{m^{(3)}})},\nabla^{\xi(\widetilde{m^{(3)}})},
\tau_*^{\xi(\widetilde{m^{(3)}})}\bigr)\bigr] \qquad \text{of} \quad \sigma\bigl(\widetilde{m^{(3)}}\bigr).
\]
By Lemma $\ref{m=0}$, we have that $\tau$ does not go through $(1,1,1,1,1) \in \mathcal{M}$. Hence, by the choice of~\smash{$\xi\bigl(\widetilde{m^{(3)}}\bigr)$} we have that $\tau$ lands in $\mathcal{S}\subset \mathcal{M}$. Furthermore, since \smash{$\sigma\bigl(\widetilde{m^{(3)}}+n\bigr)=\sigma\bigl(\widetilde{m^{(3)}}\bigr)$} for~all ${n \in \mathbb{Z}}$, we see that the same holds for $\tau$. We conclude that $\tau$ satisfies the required properties.
\end{proof}

Going back to the proof of Proposition \ref{nonvan}, let $\tau$ be a curve like in the previous lemma, associated to our chosen $\bigl(E,\overline{\partial}_E,\theta,h,g\bigr) \in \mathcal{H}^{{\rm fr}}$. This curve must then be contained in one of the two components of $\mathcal{S}$, and it goes through all the points where $a(\xi)$ could be $0$ (i.e., through the points~${\tau\bigl(m^{(3)}+n\bigr)}$ with $n \in \mathbb{Z}$). Since the components of $\mathcal{S}$ are determined by whether $a\neq0$ or~${b\neq 0}$, and we know that \smash{$a\bigl(\bigl[\bigl(E,\overline{\partial}_E,\theta,h,g\bigr)\bigr],\xi_n^{+}\bigr)=a\bigl(\bigl[\bigl(E,\overline{\partial}_E,\theta,h,g\bigr)\bigr],\xi\bigl(m^{(3)}+n\bigr)\bigr)\neq 0$} for at least one $n$, we conclude that \smash{$a\bigl(\bigl[\bigl(E,\overline{\partial}_E,\theta,h,g\bigr)\bigr],\xi_n^{+}\bigr)\neq 0$} for all $n \in \mathbb{Z}$. Since these are all the points in $\mathbb{H}_m$ where $a\bigl(\bigl[\bigl(E,\overline{\partial}_E,\theta,h,g\bigr)\bigr],\xi\bigr)$ could be $0$, we then have that $a\bigl(\bigl[\bigl(E,\overline{\partial}_E,\theta,h,g\bigr)\bigr],\xi\bigr)\neq 0$ for all $\xi \in \mathbb{H}_m$.
\end{proof}

By the remark at the beginning of the section, we then conclude.
\begin{Corollary} Given $\bigl[\bigl(E,\overline{\partial}_E,\theta,h,g\bigr)\bigr] \in \mathfrak{X}^{{\rm fr}}$ with parameter $m\neq 0$, the magnetic twistor coordinate \smash{$\MC\bigl(\bigl[\bigl(E,\overline{\partial}_E,\theta,h,g\bigr)\bigr],\xi\bigr)$} is well defined for $\xi \in \mathbb{C}^*-l_{\pm}(-2{\rm i}m)$ and it takes values in~$\mathbb{C}^*$.
\end{Corollary}
\end{subsubsection}

\section[X\^{}\{fr\} and the Ooguri--Vafa space]{$\boldsymbol{\mathfrak{X}^{{\rm fr}}}$ and the Ooguri--Vafa space}\label{sec4}

In this section, we show that there is a natural one-to-one correspondence between a subset of~$\mathfrak{X}^{{\rm fr}}$ and the Ooguri--Vafa space $\mathcal{M}^{\text{ov}}$. The correspondence is such that the coordinates $\mathcal{X}_e(\xi)$ and~$\mathcal{X}_m(\xi)$ for $\mathfrak{X}^{{\rm fr}}$ built on the last section match the twistor coordinates $\mathcal{X}_e^{\text{ov}}(\xi)$ and $\mathcal{X}_m^{\text{ov}}(\xi)$ of~$\mathcal{M}^{\text{ov}}$.

\subsection{Matching parameters with the Ooguri--Vafa space}\label{OVP}

We start by comparing the electric twistor coordinates
\begin{gather}
 \mathcal{X}_e^{\text{ov}}(\xi)=\exp\bigl[\pi \xi^{-1}z+{\rm i}\theta_e +\pi\xi \overline{z}\bigr],\nonumber\\
 \mathcal{X}_e(\xi)=\exp \bigl[-2\pi {\rm i} \bigl(\xi^{-1}m -m^{(3)}-\overline{m}\xi\bigr)\bigr]\nonumber\\
 \phantom{ \mathcal{X}_e(\xi)}{}= \exp \bigl[\pi \xi^{-1}(-2{\rm i}m) + {\rm i}\bigl(2\pi m^{(3)}\bigr) + \pi \xi \overline{(-2{\rm i}m)}\bigr].\label{elecmatch}
\end{gather}
The above formulas suggest that $z=-2{\rm i}m$, while $\theta_e=2\pi m^{(3)}$ (mod $2\pi$).

To find the analog of the magnetic angle, we will use the asymptotics of $\MC(\xi)$ from Proposition \ref{a0}
\begin{equation*}
 \MC(\xi)\sim_{\xi \to 0}\exp \left(-\frac{1}{\xi}m\operatorname{Log}_p\left(-\frac{m}{2e}\right)+{\rm i}m^{(3)}\operatorname{Arg}_p(-m)+{\rm i}\pi +\int_{-\infty}^{\infty}\gamma^*A_{11}\right),
\end{equation*}
which we can rewrite as
\begin{align}
 \MC(\xi)\sim_{\xi \to 0}&\exp \left(-\frac{{\rm i}}{2\xi}\left((-2{\rm i}m)\operatorname{Log}_p\left(-\frac{2{\rm i}m}{4{\rm i}}\right)-(-2{\rm i}m)\right)\right.\nonumber\\
 &\left.{} +{\rm i}m^{(3)}\operatorname{Arg}_p(-m)+{\rm i}\pi +\int_{-\infty}^{\infty}\gamma^*A_{11}\right).\label{asy1}
\end{align}
Comparing the last expression with the asymptotics of $\MC^{\text{ov}}(\xi)$ from Proposition \ref{asympov}
\begin{equation}\label{asy2}
 \MC^{\text{ov}}(\xi) \sim_{\xi \to 0} \exp \biggl(-\frac{{\rm i}}{2\xi}\left(z\operatorname{Log}_p\left(\frac{z}{\Lambda}\right)-z\right) +{\rm i}\theta_m +\frac{1}{2\pi {\rm i}}\sum_{s\neq 0, s \in \mathbb{Z}}\frac{1}{s}{\rm e}^{{\rm i}s\theta_e}K_0(2\pi |sz|) \biggr)
\end{equation}
we see that if $z=-2{\rm i}m $, the $\xi^{-1}$ term matches if we pick $\Lambda = 4{\rm i}$. On the other hand, when comparing the $\xi^{0}$ term, we see that by matching the imaginary part (notice that the term with the sum of Bessel functions is real), we get that $\theta_m$ should correspond to
\begin{equation*}
 m^{(3)}\operatorname{Arg}_p(-m)+\pi +\operatorname{Im}\bigg(\int_{-\infty}^{\infty}\gamma^*A_{11}\bigg) \quad (\text{mod} 2\pi).
\end{equation*}

\begin{subsubsection}[The magnetic angle on X\^{}\{fr\}]{The magnetic angle on $\boldsymbol{\mathfrak{X}^{{\rm fr}}}$}

The matching of parameters from above will let us define an analog of the Ooguri--Vafa magnetic angle on $\mathfrak{X}^{{\rm fr}}$. Before defining it, we will need the following lemmas.

\begin{Lemma}\label{mm3inv} Suppose that \smash{$\bigl(E_i,\overline{\partial}_{E_i},\theta_i,h_i,g_i\bigr)\in \mathcal{H}^{\textnormal{fr}}$} for $i=1,2$ have parameters $m_i$ and \smash{$m^{(3)}_i$} describing the singularity of $\theta_i$ and $\overline{\partial}_{E_i}$ at $z=\infty$ $($as in \eqref{singHiggs} and \eqref{singparabolic}$)$. If there is an isomorphism between these two elements, then $m_1=m_2$ and \smash{$m^{(3)}_1=m^{(3)}_2$}. In particular, it makes sense to associate $m$ and $m^{(3)}$ to elements of $\mathfrak{X}^{\textnormal{fr}}$.

\end{Lemma}

\begin{proof} If the two elements of $\mathcal{H}^{{\rm fr}}$ are isomorphic, then
\begin{equation*}
-\bigl(z^2+2m_1\bigr){\rm d}z^2=\operatorname{Det}(\theta_1)=\operatorname{Det}(\theta_2)=-\bigl(z^2+2m_2\bigr){\rm d}z^2,
\end{equation*}
so that $m_1=m_2$. On the other hand, if we denote the compatible frames (and their extensions to local frames around $z=\infty$) by $g_1=(e_1,e_2)$, $g_2=(f_1,f_2)$, and the isomorphism by $T$, we then have
\begin{equation*}
 m^{(3)}_2f_1=\overline{\partial}_{E_2}(-2\overline{w}\partial_{\overline{w}})(f_1)|_{w=0}=T\bigl(\overline{\partial}_{E_1}(-2\overline{w}\partial_{\overline{w}})(e_1)|_{w=0}\bigr)
 =T\bigl(m^{(3)}_1e_1\bigr)=m^{(3)}_1f_1,
\end{equation*}
so that \smash{$m^{(3)}_1=m^{(3)}_2$}.
\end{proof}

\begin{Definition}\label{definition4.1} Let $m\in \mathbb{C}$ and \smash{$m^{(3)}\in \bigl(-\frac12,\frac12\bigr]$}. We will denote by \smash{$\mathfrak{X}^{\textnormal{fr}}\bigl(m,m^{(3)}\bigr)\subset \mathfrak{X}^{\textnormal{fr}}$} the set of equivalence classes whose singularity with respect to the compatible framing is described by $m$ and \smash{$m^{(3)}$}.
\end{Definition}

For the following proposition, we will need the next lemma:

\begin{Lemma} \label{UHB} Let $m\in \mathbb{C}$ and \smash{$m^{(3)}\in \bigl(-\frac12,\frac12\bigr] \subset \mathbb{R}$}. Then up to equivalence, there is a unique polystable filtered Higgs bundle $(\mathcal{E}_*,\theta)$ with $\operatorname{Tr}(\theta)=0$, $\operatorname{Det}(\theta)=-\bigl(z^2+2m\bigr){\rm d}z^2$, and parabolic weights determined by \smash{$m^{(3)}$} as follows:
\begin{itemize}\itemsep=0pt
 \item if \smash{$m^{(3)} \in \bigl(-\frac12,\frac12\bigr)$}, then for the eigenline decomposition near $\infty$ of the induced $\frac12$-parabolic Higgs bundle $(\mathcal{E}_{1/2},\theta)$, we have that \smash{$\pm m^{(3)}$} is the weight associated to the line corresponding to the eigenvalue $\pm\bigl(z+\frac{m}{z}+\cdots \bigr){\rm d}z$.
 \item if \smash{$m^{(3)}=\frac12$} then the parabolic structure of the induced $\frac12$-parabolic structure $(\mathcal{E}_{1/2},\theta)$ is the trivial filtration with weight $\frac12$.
\end{itemize}
\end{Lemma}

\begin{proof}
In Appendix \ref{af}.
\end{proof}

\begin{Proposition}\label{U1A} We have a ${\rm U}(1)$-action on $\mathfrak{X}^{\textnormal{fr}}$ given by
 \begin{equation}\label{U(1)act}
 {\rm e}^{{\rm i}\theta}\cdot \bigl[\bigl(E,\overline{\partial}_E,\theta,h,g\bigr)\bigr]=\bigl[\bigl(E,\overline{\partial}_E,\theta,h,{\rm e}^{{\rm i}\theta/2}\cdot g\bigr)\bigr],
 \end{equation}
 where if $g=(e_1,e_2)$, then \smash{${\rm e}^{{\rm i}\theta/2}\cdot (e_1,e_2)=\bigl({\rm e}^{{\rm i}\theta/2}e_1,{\rm e}^{-{\rm i}\theta/2}e_2\bigr)$}. Furthermore, for $m\neq 0$ we have that~\smash{$\mathfrak{X}^{\textnormal{fr}}\bigl(m,m^{(3)}\bigr)$} is a ${\rm U}(1)$-torsor under this action.
\end{Proposition}

\begin{proof}
First notice that $\bigl(E,\overline{\partial}_E,\theta,h,g\bigr)\in \mathcal{H}^{{\rm fr}}$ is isomorphic to $\bigl(E,\overline{\partial}_E,\theta,h,{\rm e}^{{\rm i}\pi}\cdot g\bigr)$ by the isomorphism $-\text{Id}_E$ (where $\text{Id}_E\colon E\to E$ denotes the identity map), so that the map ${\rm U}(1)\times \mathcal{H}^{{\rm fr}}\to \mathfrak{X}^{{\rm fr}}$ given by
\begin{equation*}
 {\rm e}^{{\rm i}\theta}\cdot \bigl(E,\overline{\partial}_E,\theta,h,g\bigr)=\bigl[\bigl(E,\overline{\partial}_E,\theta,h,{\rm e}^{{\rm i}\theta/2}\cdot g\bigr)\bigr]
\end{equation*}
is well defined.

On the other hand, if $T$ is an isomorphism between $\bigl(E_1,\overline{\partial}_{E_1},\theta_1,h_1,g_1\bigr)$ and $\bigl(E_2,\overline{\partial}_{E_2},\theta_2,\allowbreak h_2, g_2\bigr)$, then clearly $T$ also gives an isomorphism between \smash{$\bigl(E_1,\overline{\partial}_{E_1},\theta_1,h_1,{\rm e}^{{\rm i}\theta/2}\cdot g_1\bigr)$} and $ \smash{\bigl(E_2,\overline{\partial}_{E_2}},\allowbreak\theta_2, \smash{ h_2,{\rm e}^{{\rm i}\theta/2}\cdot g_2\bigr)}$. Hence we get an ${\rm U}(1)$-action ${\rm U}(1)\times \mathfrak{X}^{{\rm fr}}\to \mathfrak{X}^{{\rm fr}}$, defined by \eqref{U(1)act}.

Now let \smash{$\mathfrak{X}^{{\rm fr}}\bigl(m,m^{(3)}\bigr)\subset \mathfrak{X}^{{\rm fr}}$} be as in Definition~\ref{definition4.1}. Clearly, the ${\rm U}(1)$-action on $\mathfrak{X}^{{\rm fr}}$ restricts to an action on \smash{$\mathfrak{X}^{{\rm fr}}\bigl(m,m^{(3)}\bigr)$}. Let us now check that this action is freely transitive:
\begin{itemize}\itemsep=0pt
 \item The action is free. Let ${\rm e}^{{\rm i}\theta}\neq 1$. From the way Stokes data transforms under changes of compatible framing, it is easy to check that
 \begin{equation}\label{equivmc}
 {\rm e}^{{\rm i}\theta}\MC(\xi)\bigl(\bigl[\bigl(E,\overline{\partial}_E,\theta,h,g\bigr)\bigr]\bigr)=\MC(\xi)\bigl({\rm e}^{{\rm i}\theta}\cdot \bigl[\bigl(E,\overline{\partial}_E,\theta,h,g\bigr)\bigr]\bigr),
 \end{equation}
 i.e., $\MC(\xi)$ is equivariant with respect to the ${\rm U}(1)$-action on $\mathfrak{X}^{{\rm fr}}$ and the natural ${\rm U}(1)$-action on $\mathbb{C}$. On the other hand, since $\MC(\xi)$ is valued in $\mathbb{C}^*$, equation \eqref{equivmc} and the fact that ${\rm e}^{{\rm i}\theta}\neq 1$ implies that $\MC(\xi)\bigl(\bigl[\bigl(E,\overline{\partial}_E,\theta,h,g\bigr)\bigr]\bigr)\neq \MC(\xi)\bigl({\rm e}^{{\rm i}\theta}\cdot\bigl[\bigl(E,\overline{\partial}_E,\theta,h,g\bigr)\bigr]\bigr)$. Since Stokes data is an isomorphism invariant, we must have \smash{$\bigl[\bigl(E,\overline{\partial}_E,\theta,h,g\bigr)\bigr]\neq {\rm e}^{{\rm i}\theta}\cdot\bigl[\bigl(E,\overline{\partial}_E,\theta,h,g\bigr)\bigr]$}.
 \item The action is transitive: Let \smash{$\bigl[\bigl(E_i,\overline{\partial}_{E_i},\theta_i,h_i,g_i\bigr)\bigr] \in \mathfrak{X}^{{\rm fr}}\bigl(m,m^{(3)}\bigr) $} for $i=1,2$. By Lem\-ma~\ref{UHB}, we have that the underlying filtered Higgs bundles \smash{$\bigl(\mathcal{P}^{h_i}_*\mathcal{E}_i^{0},\theta_i\bigr)$} are isomorphic. But then we must have that $\bigl(E_1,\overline{\partial}_{E_1},\theta_1,h_1\bigr)$ is isomorphic to $\bigl(E_2,\overline{\partial}_{E_2},\theta_2,h_2\bigr)$. Following the construction of Proposition \ref{constructionframe}, we then find an isomorphism between \smash{$\bigl(E_1,\overline{\partial_1},\theta_1,h_1,g_1\bigr)$} and~\smash{$\bigl(E_2,\overline{\partial_2},\theta_2,h_2,\widetilde{g}_2\bigr)$} for some compatible frame $\widetilde{g}_2$. Now $g_2$ and $\widetilde{g_2}$ must be ${\rm SU}(2)$ eigenframes of $\theta_2\bigl(w^3\partial_w\bigr)|_{w=0}$ with the corresponding fixed order of the eigenvalues given by the form of the singularity, so we must have $ g_2={\rm e}^{{\rm i}\theta}\cdot\widetilde{g}_2$ for some ${\rm e}^{{\rm i}\theta} \in {\rm U}(1)$. Hence,
 \begin{equation*}
 {\rm e}^{2{\rm i}\theta}\cdot \bigl[\bigl(E_1,\overline{\partial_1},\theta_1,h_1,g_1\bigr)\bigr]= {\rm e}^{2{\rm i}\theta}\cdot \bigl[\bigl(E_2,\overline{\partial_2},\theta_2,h_2,\widetilde{g}_2\bigr)\bigr]=\bigl[\bigl(E_2,\overline{\partial_2},\theta_2,h_2,g_2\bigr)\bigr].\tag*{\qed}
 \end{equation*}
\end{itemize}\renewcommand{\qed}{}
\end{proof}

Recall that for $m\neq 0$ we have the correspondence
\begin{equation}\label{magang}
 \theta_m \iff m^{(3)}\operatorname{Arg}_p(-m)+\pi +\operatorname{Im}\left(\int_{-\infty}^{\infty}\gamma^*A_{11}\right) \quad (\text{mod} 2\pi).
\end{equation}
Because of the equivariance of $\MC(\xi)$ under the ${\rm U}(1)$-actions on $\mathfrak{X}^{{\rm fr}}$ and $\mathbb{C}^*$, we see that acting on \smash{$\bigl[\bigl(E,\overline{\partial}_E,\theta,h,g\bigr)\bigr]\in \mathfrak{X}^{{\rm fr}}\bigl(m,m^{(3)}\bigr)$} by ${\rm e}^{{\rm i}\theta}$ shifts the quantity to the right of equation \eqref{magang} by $\theta$. This motivates the following definitions:

\begin{Definition} Let $m\neq 0$. The marked point of the ${\rm U}(1)$-torsor \smash{$\mathfrak{X}^{\textnormal{fr}}\bigl(m,m^{(3)}\bigr)$} is the unique element \smash{$\bigl[\bigl(E,\overline{\partial}_E,\theta,h,g_0\bigr)\bigr]\in \mathfrak{X}^{\textnormal{fr}}\bigl(m,m^{(3)}\bigr)$} such that
\begin{equation*}
 \exp \left({\rm i}m^{(3)}\operatorname{Arg}_p(-m)+{\rm i}\pi +{\rm i}\operatorname{Im}\left(\int_{-\infty}^{\infty}\gamma^*A_{11}\right)\right)=1 .
\end{equation*}
\end{Definition}
\begin{Definition}\label{magangledef} For $m \neq 0$, the magnetic angle $\theta_m$ of \smash{$\bigl[\bigl(E,\overline{\partial}_E,\theta,h,g\bigr)\bigr]\in \mathfrak{X}^{\textnormal{fr}}\bigl(m,m^{(3)}\bigr)$} is defined to be the unique real number mod $2\pi$ such that
\smash{$
 \bigl[\bigl(E,\overline{\partial}_E,\theta,h,g\bigr)\bigr]={\rm e}^{{\rm i}\theta_m}\cdot \bigl[\bigl(E,\overline{\partial}_E,\theta,h,g_0\bigr)\bigr]$},
where $\bigl[\bigl(E,\overline{\partial}_E,\theta,h,g_0\bigr)\bigr]$ is the marked point of \smash{$\mathfrak{X}^{\textnormal{fr}}\bigl(m,m^{(3)}\bigr)$}.
\end{Definition}

Notice that our definition of marked point uses a branch of the $\operatorname{Arg}$ function (with $\operatorname{Arg}(z)\in [-\pi,\pi)$), so the magnetic angle is not a priori a global continuous function of $m$. In the following, we will compute how the magnetic angle jumps when we go around a loop in the $m$ parameter. We will see that it will match the jump of the magnetic angle of the Ooguri--Vafa space.

\begin{Definition} We will say that a map $\gamma\colon [0,1]\to \mathcal{H}^{\textnormal{fr}}$ is a continuous path, if the elements~${\gamma(t):=\bigl(E_t,\overline{\partial}_{E_t},\theta_t,h_t,g_t\bigr)}$ satisfy:
\begin{itemize}\itemsep=0pt
 \item The vector bundles $E_t\to \CP$ fit into a continuous vector bundle $E\to [0,1]\times \CP$, with $E|_{\{t\}\times \CP}=E_t$, and with trivializations having transition functions depending smoothly on the points of $\CP$ and continuously on $t\in [0,1]$.
 \item In the above trivializations, the structures $\overline{\partial}_{E_t}$, $\theta_t$ and $h_t$ vary continuously in $t$.
 \item The elements of the frames $g_t$ give a continuous section of \smash{$E|_{[0,1]\times \{\infty\}}$}.
\end{itemize}

Furthermore, we will say that a map $\gamma\colon [0,1]\to \mathfrak{X}^{\textnormal{fr}}$ is a continuous path if there is a lift to a continuous path $\widetilde{\gamma}\colon [0,1]\to \mathcal{H}^{\textnormal{fr}}$.
\end{Definition}
\begin{Remark} A way to construct such paths is, for example, by employing \cite[Proposition~4.9]{M5} and the constructions of Lemma \ref{C1} and Proposition \ref{constructionframe}.
\end{Remark}
\begin{Proposition} \label{magmon} Let $l:=\{ m\in \mathbb{C} \mid m\in \mathbb{R}_{>0}\}$, and let $\gamma\colon [0,1]\to \mathfrak{X}^{\textnormal{fr}}$ be a loop in $\mathfrak{X}^{\textnormal{fr}}$ such that the associated curve $m(\gamma(t))$ of $m$ parameters of $\gamma(t)$ gives a counterclockwise loop around~${m=0}$ with $m(\gamma(0)) \in l$. Then
\begin{equation*}
 \lim_{t\to 1}\theta_m(\gamma(t))=\theta_m(\gamma(0)) +2\pi m^{(3)}-\pi.
\end{equation*}
\end{Proposition}

\begin{proof}
We will use our usual notations for Stokes data as in Section \ref{deftwistcoords}.

The magnetic angle is defined in terms of the $\xi^0$ term of the asymptotics of the Stokes data~${a(\gamma(t),\xi)}$ as $\xi \to 0$ along $\xi \in \mathbb{H}_{m(\gamma(t))}$, or by the $\xi^{0}$ term of the asymptotics of the Stokes data $-1/b(\gamma(t),\xi)$ as $\xi \to 0$ along $\xi \in \mathbb{H}_{-m(\gamma(t))}$. We will look at the monodromy of $a(\gamma(t),\xi)$, but a~similar argument holds if we use $-1/b(\gamma(t),\xi)$.

As in Section \ref{deftwistcoords}, we denote by $\Phi_i$ the frame of flat sections of $\nabla^{\xi}$ defined on the extended sector $\widehat{\operatorname{Sect}}_i(\xi)$. As we move around the loop $\gamma(t)$, the labelings of the sectors move in a clockwise manner (recall the conventions of the $m$ dependence of the labelings of Section \ref{deftwistcoords}). After going around the loop, we get that the labelings moved in such a way that the following sectors are interchanged: $\widehat{\operatorname{Sect}}_1(\xi) \iff \widehat{\operatorname{Sect}}_3(\xi)$ and $\widehat{\operatorname{Sect}}_2(\xi) \iff \widehat{\operatorname{Sect}}_4(\xi)$. If we denote by $\widetilde{\Phi}_i$ the frame of flat sections obtained after going around the loop, it is easy to see that we get the following relations:
\begin{equation}\label{mmonodromy}
 \widetilde{\Phi}_1 = \Phi_3 \cdot M_0^{-1},\qquad
 \widetilde{\Phi}_2= \Phi_4 \cdot M_0^{-1},\qquad
 \widetilde{\Phi}_3= \Phi_1,\qquad
 \widetilde{\Phi}_4 = \Phi_2.
\end{equation}

Now recall the relations that we have among the entries of the Stokes matrices \eqref{Stokes relations} and among the flat sections \eqref{Stokes relations 3}. Using these relations and \eqref{mmonodromy}, we find that
\begin{equation*}
 \widetilde{a}=\frac{\widetilde{s}_3\wedge \widetilde{s}_1}{\widetilde{s}_2\wedge \widetilde{s}_1}= \mu^{-1}\frac{ s_1 \wedge s_3}{s_4\wedge s_3}=\mu^{-2}c=-\mu^{-1} a.
\end{equation*}

Since \smash{$-\mu^{-1}(\gamma(0),\xi)=\exp \bigl(\pi \xi^{-1}(-2{\rm i}m) + {\rm i}\bigl(2\pi m^{(3)}-\pi\bigr) + \pi \xi \overline{(-2{\rm i}m)}\bigr)$}, from the asymptotics of $-\mu^{-1}(\gamma(0),\xi)a(\gamma(0),\xi)$ as $\xi \to 0$ along $\mathbb{H}_m$ we conclude what we want.
\end{proof}

Hence, we see that the magnetic angle has the same monodromy as the usual Ooguri--Vafa magnetic angle.

\end{subsubsection}

\subsection{Matching the twistor coordinates}\label{OVT}
By taking $z= -2{\rm i}m$ and \smash{$\theta_e = 2\pi m^{(3)}$}, we clearly have that for $\bigl[\bigl(E,\overline{\partial}_E,\theta,h,g\bigr)\bigr]\in \mathfrak{X}^{{\rm fr}}\bigl(m,m^{(3)}\bigr)$, the following holds (recall \eqref{elecmatch})
$\mathcal{X}_e(\xi)\bigl(\bigl[\bigl(E,\overline{\partial}_E,\theta,h,g\bigr)\bigr]\bigr)=\mathcal{X}_e^{\text{ov}}(\xi)(z,\theta_e)$.
So the remaining question is whether the magnetic twistor coordinate $\MC(\xi)$ on $\mathfrak{X}^{{\rm fr}}$ match\-es~$\MCOV(\xi)$, under the appropriate matching of parameters.

\begin{Theorem} \label{matchingTC} Fix \smash{$\bigl[\bigl(E,\overline{\partial}_E,\theta,h,g\bigr)\bigr] \in \mathfrak{X}^{\textnormal{fr}}\bigl(m,m^{(3)}\bigr)$} and let $\MC^{\textnormal{ov}}(\xi)$ be the magnetic twistor coordinate of the Ooguri--Vafa space with $\Lambda=4{\rm i}$. Then by taking $z=-2{\rm i}m$, \smash{$\theta_e=2\pi m^{(3)}$} and $\theta_m=\theta_m\bigl(\bigl[\bigl(E,\overline{\partial}_E,\theta,h,g\bigr)\bigr]\bigr)$, we have that $\MC(\xi)\bigl(\bigl[\bigl(E,\overline{\partial}_E,\theta,h,g\bigr)\bigr]\bigr)= \MC^{\textnormal{ov}}(\xi)(z,\theta_e,\theta_m)$ for all~${\xi \in \mathbb{C}^* - l_{\pm}(-2{\rm i}m)}$.
\end{Theorem}

\begin{proof} We will abbreviate the notation and just write $\MC(\xi)$ and $\MC^{\text{ov}}(\xi)$.
We consider the quotient $\MC(\xi)/\MC^{\text{ov}}(\xi)$ as a function of $\xi \in \mathbb{C}^*-l_{\pm}(-2{\rm i}m)$. Because both functions have the same jumping behavior along $l_{\pm}(-2{\rm i}m)$, we have that $\MC(\xi)/\MC^{\text{ov}}(\xi)$ extends to a continuous function on $\mathbb{C}^*$ that is holomorphic in $\xi \in \mathbb{C}^*-l_{\pm}(-2{\rm i}m)$. Furthermore, because of our matching of parameters and the asymptotics \eqref{asy1} and \eqref{asy2} of both functions, we have that
\begin{gather}
 \lim_{\xi \to 0}\frac{\MC(\xi)}{\MCOV(\xi)}= r,\label{defr}
\end{gather}
where
\begin{gather*}
 r=\exp \biggl( \operatorname{Re}\left(\int_{-\infty}^{\infty}\gamma^*A_{11}\right) -\frac{1}{2\pi {\rm i}}\sum_{s\neq 0, s \in \mathbb{Z}}\frac{1}{s}{\rm e}^{{\rm i}s\theta_e}K_0(2\pi |sz|) \biggr)\in \mathbb{R}.
\end{gather*}
The fact that $r$ is real can be easily checked by noticing that $K_0(2\pi|sz|)\in \mathbb{R}$ and that the sum is over $s\in \mathbb{Z}\setminus\{0\}$, which is invariant under $s\to -s$.

Using the reality condition of both coordinates, we also get that
\begin{equation*}
 \lim_{\xi \to \infty}\frac{\MC(\xi)}{\MCOV(\xi)}= \lim_{\xi \to \infty}\frac{\overline{\MC^{-1}\bigl(-1/\overline{\xi}\bigr)}}{\overline{\MCOVI\bigl(-1/\overline{\xi}\bigr)}}=r^{-1}.
\end{equation*}

Hence, we can extend $\MC(\xi)/\MC^{\text{ov}}(\xi)$ to a continuous function on $\CP$ which is holomorphic on $\mathbb{C}^*-l_{\pm}(-2{\rm i}m)$. An application of Morera's theorem then shows that $\MC(\xi)/\MC^{\text{ov}}(\xi)$ is a~holomorphic function on $\CP$, and hence constant. In particular, we conclude that
\smash{$
 \frac{\MC(\xi)}{\MCOV(\xi)}=r=r^{-1}$},
so that $r=\pm 1$ and $\MC(\xi)=\pm \MCOV(\xi)$. To fix the sign, notice that $r$ is an exponential of a real number (recall equation \eqref{defr}), so that we conclude that $r=1$ and then $\MC(\xi)=\MCOV(\xi)$ for all~${\xi \in \mathbb{C}^* - l_{\pm}(-2{\rm i}m)}$.
\end{proof}

\subsection[The Hyperk\"ahler structure on X\^{}\{fr\}]{The Hyperk\"ahler structure on $\boldsymbol{\mathfrak{X}^{{\rm fr}}}$}
Consider the Ooguri--Vafa space $\mathcal{M}^{\text{ov}}$ together with its torus fibration $\mathcal{M}^{\text{ov}}\to \mathcal{B}$. Throughout this section, we fix the cut-off parameter specifying $\mathcal{M}^{\text{ov}}$ to be $\Lambda=4{\rm i}$.

\begin{Definition} We denote
\begin{gather*}
\mathfrak{X}^{\textnormal{fr}}(\mathcal{B}):=\bigl\{\bigl[\bigl(E,\overline{\partial}_E,\theta,h,g\bigr)\bigr] \!\in \mathfrak{X}^{\textnormal{fr}} \mid \operatorname{Det}(\theta)=-\bigl(z^2+2m\bigr){\rm d}z^2\ \text{with $-2{\rm i}m\in \mathcal{B}$}\bigr\},\\
 \mathfrak{X}^{\textnormal{fr}}_*(\mathcal{B}):=\bigl\{\bigl[\bigl(E,\overline{\partial}_E,\theta,h,g\bigr)\bigr]\! \in \mathfrak{X}^{\textnormal{fr}} \mid \operatorname{Det}(\theta)=-\bigl(z^2+2m\bigr){\rm d}z^2\ \text{with $m\neq 0$ and $-2{\rm i}m\in \mathcal{B}$}\bigr\}.
\end{gather*}
Here $z\in \mathbb{C}\subset \CP$ denotes the fixed holomorphic coordinate from Section \ref{secdeffwhb}. Furthermore, we denote $\mathcal{B}_*=\mathcal{B}\setminus\{0\}$ and $\mathcal{M}_*^{\textnormal{ov}} \subset \mathcal{M}^{\textnormal{ov}}$ the points of the Ooguri--Vafa space over $\mathcal{B}_*$.
\end{Definition}

The results of Sections \ref{OVP} and \ref{OVT} then allow us to show the following.

\begin{Theorem}\label{HKM} There is a one-to-one correspondence between $\mathfrak{X}^{\textnormal{fr}}_*(\mathcal{B})$ and $\mathcal{M}^{\textnormal{ov}}_*$ such that ${\mathcal{X}_e=\mathcal{X}_e^{\text{ov}}}$ and $\mathcal{X}_m=\mathcal{X}_m^{\text{ov}}$. Under this correspondence $\mathfrak{X}^{\textnormal{fr}}_*(\mathcal{B})$ acquires the structure of a hyperk\"ahler manifold. For $\xi \in \mathbb{C}^*$, the twistor family $\Omega(\xi)$ of holomorphic symplectic forms on~$\mathfrak{X}^{\textnormal{fr}}_*(\mathcal{B})$ is given by
\begin{equation}\label{familyholsym2}
 \Omega(\xi)=\frac{{\rm d}\mathcal{X}_e(\xi)}{\mathcal{X}_e(\xi)}\wedge \frac{{\rm d}\MC(\xi)}{\MC(\xi)}.
\end{equation}
\end{Theorem}

\begin{proof}
Consider the open cover of $\mathcal{B}_*\times {\rm U}(1)$ given by (compare with the end of Section \ref{constructionOV})
\begin{gather*}
 U_1:=\bigl\{\bigl(z,{\rm e}^{2\pi{\rm i}x^{3}}\bigr) \in \mathcal{B}_*\times {\rm U}(1) \mid {\rm e}^{z/4{\rm i}}\not\in \mathbb{R}_{<0} \bigr\}, \\
 U_2:=\bigl\{\bigl(z,{\rm e}^{2\pi{\rm i}x^{3}}\bigr) \in \mathcal{B}_*\times {\rm U}(1) \mid {\rm e}^{z/4{\rm i}}\not\in \mathbb{R}_{>0} \bigr\}.
\end{gather*}
We then can write $U_1\cap U_2= V_1\cup V_2$, where
\begin{gather*}
 V_1=\bigl\{\bigl(z,{\rm e}^{2\pi {\rm i}x^3}\bigr) \in \mathcal{B}_*\times {\rm U}(1) \mid \operatorname{Im}\bigl({\rm e}^{z/4{\rm i}}\bigr)>0\bigr\}, \\
 U_{-}=\bigl\{\bigl(z,{\rm e}^{2\pi {\rm i}x^3}\bigr) \in \mathcal{B}_*\times {\rm U}(1) \mid \operatorname{Im}\bigl({\rm e}^{z/4{\rm i}}\bigr)<0\bigr\}.
\end{gather*}
Furthermore, let $\pi\colon \mathfrak{X}_*^{{\rm fr}}(\mathcal{B}) \to \mathcal{B}_*\times S^1$ be defined by \smash{$\pi\bigl(\bigl[\bigl(E,\overline{\partial}_E,\theta,h,g\bigr)\bigr]\bigr)=\bigl(-2{\rm i}m,{\rm e}^{2\pi {\rm i} m^{(3)}}\bigr)$}, where $m$ and \smash{$m^{(3)}$} are the associated parameters of $\bigl[\bigl(E,\overline{\partial}_E,\theta,h,g\bigr)\bigr]$ specifying the singularity. We then have that \smash{$\bigl(-2{\rm i}m, {\rm e}^{2\pi {\rm i} m^{(3)}},{\rm e}^{{\rm i}\theta_m}\bigr)$}, where $\theta_m$ is the magnetic angle from Definition \ref{magangledef}, gives coordinates for $\pi^{-1}\bigl(U_1\times S^1\bigr)\subset \mathfrak{X}_*^{{\rm fr}}(\mathcal{B})$. Similarly, Proposition \ref{magmon} implies that we can find another magnetic angle coordinate \smash{$\widetilde{\theta}_m$} over \smash{$U_2\times S^1$} such that \smash{$\bigl(-2{\rm i}m, {\rm e}^{2\pi {\rm i} m^{(3)}},{\rm e}^{{\rm i}\widetilde{\theta}_m}\bigr)$} gives coordinates for $\pi^{-1}\bigl(U_2\times S^1\bigr)$, and such that the change of coordinates map $\phi\colon (V_1\cup V_2)\times {\rm U}(1) \to (V_1\cup V_2)\times {\rm U}(1)$ is given by
\begin{equation*}
 \phi\bigl(-2{\rm i}m,{\rm e}^{2\pi {\rm i} m^{(3)}},{\rm e}^{{\rm i}\theta_m}\bigr)=\begin{cases}
 \bigl(-2{\rm i}m,{\rm e}^{2\pi {\rm i} m^{(3)}},{\rm e}^{{\rm i}\theta_m}\bigr)& \text{if} \ \bigl(-2{\rm i}m,{\rm e}^{2\pi {\rm i} m^{(3)}}\bigr) \in V_1,\\
 \bigl(-2{\rm i}m,{\rm e}^{2\pi {\rm i} m^{(3)}},{\rm e}^{{\rm i}\theta_m +2\pi {\rm i} m^{(3)}-{\rm i}\pi}\bigr) &
 \text{if}\ \bigl(-2{\rm i}m,{\rm e}^{2\pi {\rm i} m^{(3)}}\bigr)\in V_2.
 \end{cases}
\end{equation*}
This matches the description of $\mathcal{M}^{\text{ov}}_*$ given at the end of Section \ref{constructionOV}, and hence gives a one-to-one correspondence between $\mathfrak{X}^{\textnormal{fr}}_*(\mathcal{B})$ and $\mathcal{M}^{\textnormal{ov}}_*$ such that $\mathcal{X}_e=\mathcal{X}_e^{\text{ov}}$ and $\mathcal{X}_m=\mathcal{X}_m^{\text{ov}}$. The identification induces a hyperk\"ahler structure on $\mathfrak{X}_*^{{\rm fr}}(\mathcal{B})$, and by Theorem \ref{matchingTC} and \eqref{familyholsym}, the twistor family of holomorphic symplectic forms is given by the formula~\eqref{familyholsym2}.
\end{proof}

However, this is not the end of the story, since the hyperk\"ahler structure of $\mathcal{M}^{\text{ov}}_*$ actually extends to the points over $0\in \mathcal{B}$ (see, for example, \cite{GW} or \cite[Section 4.1]{GMN1}). The fiber over $0\in \mathcal{B}$ of $\mathcal{M}^{\text{ov}}\to \mathcal{B}$ is a torus with a node (recall Figure~\ref{fig1}). In the following, we show that under our identification of parameters, we get the same picture for the elements of $\mathfrak{X}^{{\rm fr}}$ over $m=0$. We remark that while $\mathcal{X}_e(\xi)$ clearly extends to the elements of $\mathfrak{X}^{{\rm fr}}(\mathcal{B})$ with $m=0$ (recall \eqref{elecmatch}), it is not clear whether $\mathcal{X}_{m}(\xi)$ does, since the condition $m\neq 0$ is heavily used in the definition and its properties. This issue already appears for $\mathcal{X}_{m}^{\text{ov}}(\xi)$, since \eqref{Xmsf} and \eqref{Xminst} are ill-defined over~${0\in \mathcal{B}}$.

\begin{subsubsection}[The central fiber of X\^{}\{fr\}]{The central fiber of $\boldsymbol{\mathfrak{X}^{{\rm fr}}}$}

Recall that given $m \in \mathbb{C}$ and $m^{(3)} \in \bigl(-\frac12,\frac12\bigr]$, we denote by $\mathfrak{X}^{{\rm fr}}\bigl(m,m^{(3)}\bigr)$ the subset of equivalence classes of $\mathfrak{X}^{{\rm fr}}$ with corresponding parameters $m$ and \smash{$m^{(3)}$} describing the singularities. We have only shown that $\mathfrak{X}^{{\rm fr}}\bigl(m,m^{(3)}\bigr)$ is not empty when $m\neq 0$ (see Lemma \ref{C1} and Proposition~\ref{constructionframe}). In Appendix \ref{ae}, we show that \smash{$\mathfrak{X}^{{\rm fr}}\bigl(0,m^{(3)}\bigr)$} is also not empty for every~\smash{$m^{(3)} \in \bigl(-\frac12,\frac12\bigr]$}.

We now prove the following Lemma, which is analogous to Proposition \ref{U1A}.

\begin{Lemma} \label{U1A2} For $m^{(3)}\neq 0$, we have that $\mathfrak{X}^{\textnormal{fr}}\bigl(0,m^{(3)}\bigr)$ is a ${\rm U}(1)$-torsor under the ${\rm U}(1)$-action defined on Proposition {\rm\ref{U1A}}.
\end{Lemma}

\begin{proof}
Let $\bigl[\bigl(E,\overline{\partial}_E,\theta,h,g\bigr)\bigr] \in \mathfrak{X}^{{\rm fr}}\bigl(0,m^{(3)}\bigr)$ for $m^{(3)}\neq 0$. Let us now show that the ${\rm U}(1)$-action is freely transitive:
\begin{itemize}\itemsep=0pt
 \item The action is free: by the taking the associated Stokes data, we see that the formal monodromy of the flat connection $\nabla^{\xi}$ turns out to be ${\rm e}^{2\pi {\rm i} m^{(3)}}\neq 1$ (notice that it does not depend on $\xi \in \mathbb{C}^*$). Using the notation from Section \ref{GMD}, we then see that the relation \smash{$1+a(\xi)b(\xi)=\mu^{-1}(\xi)={\rm e}^{2\pi {\rm i} m^{(3)}}$} implies that $a(\xi)\neq 0$ and $b(\xi)\neq 0$ for all $\xi \in \mathbb{C}^*$. On the other hand, let ${\rm e}^{{\rm i}\theta}\neq 1$. It is then easy to check that
 \begin{equation*}
 a\bigl(\xi, {\rm e}^{{\rm i}\theta}\cdot \bigl[\bigl(E,\overline{\partial}_E,\theta,h,g\bigr)\bigr]\bigr)={\rm e}^{{\rm i}\theta}a\bigl(\xi,\bigl[\bigl(E,\overline{\partial}_E,\theta,h,g\bigr)\bigr]\bigr),
 \end{equation*}
 and since $a(\xi)\neq 0$, we must have that
 \begin{equation*}
 a\bigl(\xi, {\rm e}^{{\rm i}\theta}\cdot \bigl[\bigl(E,\overline{\partial}_E,\theta,h,g\bigr)\bigr]\bigr)\neq a\bigl(\xi,\bigl[\bigl(E,\overline{\partial}_E,\theta,h,g\bigr)\bigr]\bigr).
 \end{equation*}
 Since Stokes data is an isomorphism invariant, we conclude that
\[
 \bigl[\bigl(E,\overline{\partial}_E,\theta,h,g\bigr)\bigr]\neq{\rm e}^{{\rm i}\theta}\cdot \bigl[\bigl(E,\overline{\partial}_E,\theta,h,g\bigr)\bigr],
\]
 so the action is free.
 \item The action is transitive: the same proof as Proposition \ref{U1A}.\hfill $\qed$
\end{itemize}\renewcommand{\qed}{}
\end{proof}
\end{subsubsection}

\begin{Lemma}\label{U1A3}$\mathfrak{X}^{\textnormal{fr}}(0,0)$ is just a point.
\end{Lemma}
\begin{proof}
Let \smash{$\bigl[\bigl(E_0,\overline{\partial}_{E_0},\theta_0,h_0,g_0\bigr)\bigr]\in \mathfrak{X}^{{\rm fr}}(0,0)$} be the framed wild harmonic bundle from Example~\ref{trivWHB}, and let \smash{$\bigl[\bigl(E,\overline{\partial}_E,\theta,h,g\bigr)\bigr] \in \mathfrak{X}^{{\rm fr}}(0,0)$} be any other framed wild harmonic bundle. By the same proof of the transitivity of the ${\rm U}(1)$-action of Proposition \ref{U1A}, we have that
\begin{equation*}
 \bigl[\bigl(E,\overline{\partial}_E,\theta,h,g\bigr)\bigr]={\rm e}^{{\rm i}\theta}\cdot \bigl[\bigl(E_0,\overline{\partial}_{E_0},\theta_0,h_0,g_0\bigr)\bigr]
\end{equation*}
for some ${\rm e}^{{\rm i}\theta} \in {\rm U}(1)$.

But from the description of \smash{$\bigl[\bigl(E_0,\overline{\partial}_{E_0},\theta_0,h_0,g_0\bigr)\bigr]$} in Example \ref{trivWHB} it is easy to see that \smash{$\bigl[\bigl(E_0,\overline{\partial}_{E_0},\theta_0,h_0,g_0\bigr)\bigr]={\rm e}^{{\rm i}\theta}\cdot\bigl[\bigl(E_0,\overline{\partial}_{E_0},\theta_0,h_0,g_0\bigr)\bigr]$} for every ${\rm e}^{{\rm i}\theta}\in {\rm U}(1)$, the isomorphism being the map $T\colon E_0\to E_0$ described in the global canonical frame $(e_1,e_2)$ by
\begin{equation*}
 T=\begin{bmatrix}{\rm e}^{{\rm i}\frac{\theta}{2}} & 0\\
 0 & {\rm e}^{-{\rm i}\frac{\theta}{2}}
 \end{bmatrix}.
\end{equation*}
Hence, \smash{$\bigl[\bigl(E_0,\overline{\partial}_{E_0},\theta_0,h_0,g_0\bigr)\bigr]=\bigl[\bigl(E,\overline{\partial}_E,\theta,h,g\bigr)\bigr]$} and $\mathfrak{X}^{{\rm fr}}(0,0)$ is just a point.
\end{proof}

\begin{Corollary}\label{HKMext} The hyperk\"ahler structure of $\mathfrak{X}^{\textnormal{fr}}_*(\mathcal{B})$ extends to $\mathfrak{X}^{\textnormal{fr}}(\mathcal{B})$.
\end{Corollary}
\begin{proof}
 In Theorem \ref{HKM}, we established a one-to-one correspondence between $\mathfrak{X}^{{\rm fr}}_*(\mathcal{B})$ and $\mathcal{M}_*^{\text{ov}}$ such that $\mathcal{X}_e(\xi)=\mathcal{X}_e^{\text{ov}}(\xi)$ and $\mathcal{X}_m(\xi)=\mathcal{X}_m^{\text{ov}}(\xi)$, and $z=-2{\rm i}m$, $\theta_e=2\pi m^{(3)}$. On the other hand, we know that the fiber of $\mathcal{M}^{\text{ov}}\to \mathcal{B}$ over $0\in \mathcal{B}$ is a singular torus with a node (see Figure~\ref{fig1}) and that the node corresponds to $z=\theta_e=0$. Lemmas \ref{U1A2} and \ref{U1A3} then imply that our one-to-one correspondence $\mathfrak{X}_{*}^{{\rm fr}}(\mathcal{B})\cong \mathcal{M}^{\text{ov}}_*$ naturally extends to $\mathfrak{X}^{{\rm fr}}(\mathcal{B})\cong \mathcal{M}^{\text{ov}}$. Since the hyperk\"ahler structure of $\mathcal{M}_*^{\text{ov}}$ extends over the singular fiber \cite{GW} to $\mathcal{M}^{\text{ov}}$, and the hyperk\"ahler structure on $\mathfrak{X}_*^{{\rm fr}}(\mathcal{B})$ is induced from the one-to-one correspondence of Theorem \ref{HKM}, we conclude that the hyperk\"ahler structure of $\mathfrak{X}_*^{{\rm fr}}(\mathcal{B})$ extends to $\mathfrak{X}^{{\rm fr}}(\mathcal{B})$.
\end{proof}

Overall, joining Theorems \ref{matchingTC} and \ref{HKM}, and Corollary \ref{HKMext}, we obtain our main result stated at the beginning in Theorem~\ref{maintheorem}.

\appendix
\section{Estimates for the connection form}\label{AA}

We use the setting and notation of the beginning of the proof of Proposition \ref{constructionframe}. We want to prove the following.

\begin{Lemma} In the orthonormal frame $(e_1,e_2)$, the Chern connection is expressed as
\begin{equation*}
 D={\rm d} - \frac{1}{2}\begin{bmatrix} a(v_1) & 0 \\
 0 & a(v_2) \\
 \end{bmatrix}\left(\frac{{\rm d}w}{w}- \frac{{\rm d}\overline{w}}{\overline{w}}\right) + \text{regular terms at $w=0$}.
\end{equation*}

\end{Lemma}
\begin{proof}
The gauge transformation $g$ that satisfies $(v_1,v_2)\cdot g=(e_1,e_2)$ is given by
\begin{equation*}
 g=\begin{bmatrix}
 \frac{1}{|v_1|_h}&-\frac{h(v_1,v_2)}{|v_1|^2_h|v_2-h(v_1,v_2)|v_1|_h^{-2}v_1|_h}\\
 0 & \frac{1}{|v_2-h(v_1,v_2)|v_1|_h^{-2}v_1|_h}\\
 \end{bmatrix}.
\end{equation*}
Furthermore, we know that $|v_i|_h^2=|w|^{-2a(v_i)}f_i(w)$, where $f_i(w)$ is a positive real function that is bounded near $w=0$ (this is a consequence of Theorem \ref{FT2}). If $A$ denotes the connection matrix of $D$ in the frame $(e_1,e_2)$, then we have that
\begin{equation*}
 A=g^{-1}{\rm d}g + g^{-1}\left(\begin{bmatrix} -a(v_1) & 0 \\
 0 & -a(v_2) \\
 \end{bmatrix}\frac{{\rm d}w}{w} + \text{regular}\right)g.
\end{equation*}
By using the fact that the off-diagonal terms of $g$ and the off-diagonal terms of the regular terms of the connection matrix of $D$ in the frame $(v_1,v_2)$ go to 0 exponentially as $w\to 0$, it is easy to check that
\begin{equation*}
 g^{-1}\left(\begin{bmatrix} -a(v_1) & 0 \\
 0 & -a(v_2) \\
 \end{bmatrix}\frac{{\rm d}w}{w} + \text{regular}\right)g=\begin{bmatrix} -a(v_1) & 0 \\
 0 & -a(v_2) \\
 \end{bmatrix}\frac{{\rm d}w}{w}+\text{regular}.
\end{equation*}
On the other hand, the fact that in the frame $(v_1,v_2)$ we have that
\begin{equation*}
 D=D_0 +\text{regular}= {\rm d} + H^{-1}\partial H,
\end{equation*}
where $H$ is the matrix of the hermitian metric in the frame $(v_1,v_2)$,
implies that the functions~${f_i^{-1}\partial_{z} f_i}$ are regular at $w=0$. Since the $f_i$ are real, we get that $f_i^{-1}\partial_{\overline{z}}f_i$ is also regular, so that $f_i^{-1}{\rm d}f_i=-f_i{\rm d}f_i^{-1}$ is regular at $w=0$. From this fact, we conclude that
\begin{equation*}
 g^{-1}dg= \frac{1}{2}\begin{bmatrix} a(v_1) & 0 \\
 0 & a(v_2) \\
 \end{bmatrix}\left(\frac{{\rm d}w}{w}+ \frac{{\rm d}\overline{w}}{\overline{w}}\right) + \text{regular}.
\end{equation*}

Hence, in the frame $(e_1,e_2)$, we get that
\begin{equation*}
 D={\rm d}+A={\rm d} - \frac{1}{2}\begin{bmatrix} a(v_1) & 0 \\
 0 & a(v_2) \\
 \end{bmatrix}\left(\frac{{\rm d}w}{w}- \frac{{\rm d}\overline{w}}{\overline{w}}\right) + \text{regular}.\tag*{\qed}
\end{equation*} \renewcommand{\qed}{}
\end{proof}

\section{Proof of Lemma \ref{Lemma}}\label{AC}

The goal of this appendix is to show that the asymptotics of Lemma \ref{twistasymp} hold uniformly in~${\xi \in U(\xi_0)}$, where $\xi_0 \in \mathbb{C}^*$ and $U(\xi_0)$ is some small bounded neighborhood of $\xi_0$. Because of the expression of the flat frames $\Phi_i(\xi)$ in terms of extensions of the compatible frames $g=(e_1,e_2)$ of the harmonic bundles, it is easy to check that it is enough to show that
\begin{itemize}\itemsep=0pt
 \item $\Sigma_i\bigl(\widehat{F}(\xi)\bigr) \to 1$ uniformly in $\xi \in U(\xi_0)$ as $w\to 0$ with $w\in S_i \subset \widehat{\operatorname{Sect}}_i(\xi)$. Here $S_i$ is the sector defined at the beginning of Section \ref{proofhol}.
 \item $g_{\xi}(w)\to 1$ uniformly in $\xi \in U(\xi_0)$ as $w\to 0$.
\end{itemize}

We will only prove the second statement, since the first one follows from the proof of \cite[Theorem 7]{B2}. For the proof of the second statement, we will follow similar arguments and notations as those found in \cite{BB04}, where they construct $g_{\xi}(w)$ for $\xi=1$.

We will denote by $D \subset \CP$ the unit disc centered at $w=\frac{1}{z}=0$ with radial coordinate~$r$. We also denote by $E$ the vector bundle corresponding to an element \smash{$\bigl(E,\overline{\partial}_E,\theta,h,g\bigr)\in \mathcal{H}^{{\rm fr}}$}, trivialized over $D$ by an extension of the compatible framing at $w=0$ to a local ${\rm SU}(2)$ framing. Finally, we denote by $\Gamma(D,\operatorname{End}(E))$ the set of sections (with no regularity assumed) of~${\operatorname{End}(E)\to D}$. For $\delta >0$, we then define the weighted Sobolev spaces, as in \cite{BB04}
\begin{gather*}
 L_{\delta}^p=\left\{ f\in \Gamma(D,\operatorname{End}(E)) \,\Big|\, \frac{f}{r^{\delta +2/p}} \in L^p(D,\operatorname{End}(E))\right\},\\
 L_{\delta}^{p,k}=\left\{ f\in \Gamma(D,\operatorname{End}(E)) \,\Big|\, \frac{D\bigl(\overline{\partial}_E,h\bigr)^{j}f_i}{r^{i(k-j)}} \in L^p_{\delta} \ \text{for}\ i=0,3,\, 0\leq j \leq k\right\},
\end{gather*}
where derivatives are considered in the weak sense, and where $f_0$ and $f_3$ denote the diagonal and off-diagonal components of $f\in \operatorname{End}(E)$. The reason for the strange indexing notation for the diagonal and off-diagonal part is so that our notation agrees with $\cite{BB04}$. The highest order pole of the singularity (cubic order in our case) acts non-trivially via the adjoint action on the off-diagonal part of a section $f\in \Gamma(D,\operatorname{End}(E))$ (hence the ``$3$'' subscript), while the whole singular part acts trivially on the diagonal part of $f$ (hence the ``$0$'' subscript).

Similarly, we have the Banach spaces $C^k_{\delta}$ defined by
\begin{equation*}
 C^k_{\delta}=\left\{ f \in \Gamma\bigl(\overline{D},\text{End(E)}\bigr) \,\Big|\, \frac{f}{r^{\delta}} \in C^k\bigl(\overline{D},\operatorname{End}(E)\bigr)\right\}.
\end{equation*}

The gauge transformation $g_\xi$ is built as a solution to the following problem: with respect to an extension of the compatible framing, we have the expression
\begin{equation*}
 \bigl(\nabla^{\xi}\bigr)^{0,1}=\overline{\partial} -\xi H\frac{{\rm d}\overline{w}}{\overline{w}^3} - \left(\xi\overline{m}+ \frac{m^{(3)}}{2}\right)H\frac{{\rm d}\overline{w}}{\overline{w}} + a_{\text{reg}}^{0,1}+\xi \theta_{\text{reg}}^{\dagger_h},
\end{equation*}
where $a_{\text{reg}}^{0,1}$ and $\theta_{\text{reg}}^{\dagger_h}$ denote the regular parts of \smash{$D\bigl(\overline{\partial}_E,h\bigr)^{(0,1)}$} and $\theta^{\dagger_h}$, respectively. The gauge transformation $g_{\xi}$ then satisfies
\begin{gather}
 g_{\xi}\cdot \bigl(\nabla^{\xi}\bigr)^{0,1}=\overline{\partial} -\xi H\frac{{\rm d}\overline{w}}{\overline{w}^3} - \left(\xi\overline{m}+ \frac{m^{(3)}}{2}\right)H\frac{{\rm d}\overline{w}}{\overline{w}},\qquad
 g_{\xi}(w=0)=1,\nonumber\\
 g_{\xi,0}-1\in C^0_{\delta} g_{\xi,3} \in C^0_{2+\delta} \qquad \text{for some}\quad \delta>0.\label{problem}
\end{gather}
In order to show the existence of such a $g_{\xi}$, we slightly extend some of the results in \cite[Section~7]{BB04} in order to get statements for families in~$\xi$.

\begin{Lemma} \label{BBext} Take $\delta \in \mathbb{R}-\mathbb{Z}$ and $p>2$. On the unit disk, the problem
\begin{equation}\label{1CR}
 \frac{\partial f}{\partial \overline{w}} = g
\end{equation}
has a solution $f=T(g)$ such that \smash{$|f|_{C^0_{-1+\delta}} \leq c|g|_{L^p_{-2+\delta}}$}. Furthermore, if $\lambda(\xi)$ is a continuous~function of $\xi \in U(\xi_0)$, then the same is true if $\delta - \operatorname{Re}(\lambda(\xi)) \in \mathbb{R}-\mathbb{Z}$ for the problem%
\begin{equation}\label{2CR}
 \frac{\partial f}{\partial \overline{w}} -\frac{\lambda(\xi)}{2\overline{w}}f = g.
\end{equation}
By picking $U(\xi_0)$ small enough, we have \smash{$|T_{\xi}(g)|_{C^0_{-1+\delta}} \leq c|g|_{L^p_{-2+\delta}}$} for all $\xi \in U(\xi_0)$, for a uniform constant $c$.
\end{Lemma}

\begin{proof}
We will only show the last two statements of the lemma, since the first is the same as~\cite[Lemma 7.2]{BB04}.

By the same argument given in \cite[Lemma 7.2]{BB04}, we can assume for our problem that $\delta -\operatorname{Re}(\lambda(\xi))\in (0,1)$ for $\xi \in U(\xi_0)$. By shrinking $U(\xi_0)$ further if necessary, we can assume that~${\delta -\operatorname{Re}(\lambda(\xi))\in (\delta_0,\delta_1)}$ for $0<\delta_0<\delta_1<1$ and $\xi \in U(\xi_0)$. If $T$ denotes the solution operator for the first inhomogeneous Cauchy--Riemann problem \eqref{1CR}, then \smash{$T_{\xi}(g):=r^{\lambda(\xi)}T\bigl(r^{-\lambda(\xi)}g\bigr)$} is the solution operator for \eqref{2CR}.

Now notice that by applying H\"{o}lder's inequality, we get
\begin{align}
 |T_{\xi}(g)(w)|&=\Biggl|r^{\lambda(\xi)}\int_{D}\frac{|u|^{-\lambda(\xi)}g(u)}{w-u}\big|{\rm d}^2u\big|\Biggr|\nonumber\\
 &\leq r^{\operatorname{Re}(\lambda(\xi))}|g|_{L^p_{-2+\delta}}\left(\int_D \frac{\big|{\rm d}^2u\big|}{|u|^{2-(\delta -\operatorname{Re}(\lambda(\xi)))\frac{p}{p-1}}|w-u|^{\frac{p}{p-1}}}\right)^{\frac{p-1}{p}}.\label{B7}
\end{align}
We will denote $\delta(\xi)=\delta - \operatorname{Re}(\lambda(\xi))$. By our conditions on $\delta(\xi)$ and $p$, we have
\begin{equation*}
\int_{\mathbb{C}} \frac{\big|{\rm d}^2u\big|}{|u|^{2-\delta(\xi)\frac{p}{p-1}}|w-u|^{\frac{p}{p-1}}}<\infty,
\end{equation*}
so if $D_{1/|w|}$ denotes the disk centered at the origin of radius $1/|w|$, we can write
\begin{align*}
\int_D \frac{\big|{\rm d}^2u\big|}{|u|^{2-\delta(\xi)\frac{p}{p-1}}|w-u|^{\frac{p}{p-1}}}&=\frac{1}{|w|^{(1-\delta(\xi))\frac{p}{p-1}}}\int_{D_{1/|w|}}\frac{\big|{\rm d}^2u\big|}{|u|^{2-\delta(\xi)\frac{p}{p-1}}|1-u|^{\frac{p}{p-1}}}\\
&\leq \frac{1}{|w|^{(1-\delta(\xi))\frac{p}{p-1}}}\int_{\mathbb{C}}\frac{\big|{\rm d}^2u\big|}{|u|^{2-\delta(\xi)\frac{p}{p-1}}|1-u|^{\frac{p}{p-1}}}
=\frac{1}{|w|^{(1-\delta(\xi))\frac{p}{p-1}}} c(\xi),
\end{align*}
where $c(\xi)\colon U(\xi_0) \to \mathbb{R}$ is defined by
\begin{equation*}
 c(\xi):=\int_{\mathbb{C}}\frac{\big|{\rm d}^2u\big|}{|u|^{2-\delta(\xi)\frac{p}{p-1}}|1-u|^{\frac{p}{p-1}}}.
\end{equation*}

Now it is easy to check that $c(\xi)$ depends continuously on $\xi$. To show this, consider the function $h(u)$ defined in the following way
\begin{gather*}
 h(u)=\frac{1}{|u|^{2-\delta_0\frac{p}{p-1}}|1-u|^{\frac{p}{p-1}}} \qquad \text{if} \quad |u|< 1,\\
 h(u)=\frac{1}{|u|^{2-\delta_1\frac{p}{p-1}}|1-u|^{\frac{p}{p-1}}} \qquad \text{if} \quad |u|> 1.
\end{gather*}
Then we get that $h \in L^1(\mathbb{C})$ and furthermore
\begin{equation*}
 \frac{1}{|u|^{2-\delta(\xi)\frac{p}{p-1}}|1-u|^{\frac{p}{p-1}}}\leq h(u) \qquad \text{for every}\quad \xi \in U(\xi_0)\quad \text{and for almost every} \quad u \in \mathbb{C}.
\end{equation*}

From this fact, it follows that $c(\xi)$ must be a continuous function of $\xi \in U(\xi_0)$. In particular, by further restricting $U(\xi_0)$ if necessary, we can find a constant $C$ such that $c(\xi)\leq C$ for~${\xi \in U(\xi_0)}$. Hence, by going back \eqref{B7}, we conclude that
\begin{equation*}
 |T_{\xi}(g)(w)|\leq |g|_{L^p_{-2+\delta}} \frac{|w|^{\operatorname{Re}(\lambda(\xi))}C^{\frac{p-1}{p}}}{|w|^{1-\delta(\xi)}}=|g|_{L^p_{-2+\delta}} \frac{C^{\frac{p-1}{p}}}{|w|^{1-\delta}},
\end{equation*}
so we finally get that
\smash{$
 |T_{\xi}(g)|_{C^0_{-1+\delta}}\leq |g|_{L^p_{-2+\delta}} C^{\frac{p-1}{p}} $} for all $ \xi \in U(\xi_0)$.
\end{proof}

\begin{Lemma} Let
\begin{gather*}
U_{\delta}:=\bigl\{u \in \Gamma(D,\operatorname{End}(E)) \mid u_3 \in C^0_{2+\delta},\, u_0 \in C^0_{\delta}\bigr\},\\
A_{\delta}:= \bigl\{a \in \Gamma(D,\operatorname{End}(E)) \mid a_3 \in L^{p}_{\delta+1},\, a_0 \in L^{p}_{\delta-1}\bigr\}.
\end{gather*}
 Furthermore, let \smash{$\overline{\partial}_0=
 \overline{\partial} - \bigl(\xi\overline{m}+ \frac{m^{(3)}}{2}\bigr)H\frac{{\rm d}\overline{w}}{\overline{w}}$}. Then for some $\delta'<\delta$ there is a~continuous map $T_{\xi}\colon A_{\delta'}\to U_{\delta'}$ such that \smash{$\overline{\partial}_0(T_{\xi}(c))=c$}. If we pick $U(\xi_0)$ sufficiently small, the family of solution maps $T_{\xi}$ has a uniform bound in $\xi$.
\end{Lemma}

\begin{proof}
This follows from the beginning of the proof of \cite[Lemma 7.1]{BB04} and our previous Lem\-ma~\ref{BBext}.
\end{proof}

\begin{Theorem} For a sufficiently small disk $D_{\lambda}$ centered at $w=0$ and $U(\xi_0)$ bounded and sufficiently small, we have a solution $g_{\xi}(w)$ of the problem described in \eqref{problem} that is defined on ${D_{\lambda}\times U(\xi_0)}$ and depends continuously on $\xi$. Furthermore, $g_{\xi}(w) \to 1$ uniformly in $\xi$ as~${w\to 0}$.
\end{Theorem}

\begin{proof} We will follow mostly the same argument as in the proof of \cite[Lemma 7.1]{BB04}. We put it here just to emphasize the behavior in families parametrized by $\xi$, which is not done in the aforementioned paper.

For $\lambda >0$, let $h_{\lambda}$ be the homothety $h_{\lambda}(w)=\lambda w$, and let \smash{$\varphi(w,\xi)=\exp \bigl(\bigl(\frac{\overline{\xi}}{2w^2}-\frac{\xi}{2\overline{w}^2}\bigr)H\bigr)$}. Furthermore, we denote \smash{$c(w,\xi)=-a_{\text{reg}}^{0,1}(w)-\xi \theta_{\text{reg}}^{\dagger_h}(w)$}.

If we write $g_{\xi}(w)=1+u(w,\xi)$, the problem \eqref{problem} that we are trying to solve can be rephrased as the problem of finding $u(w,\xi)$ such that
\begin{gather*}
 \overline{\partial}_0(u(w,\xi))=\left[\xi H \frac{{\rm d}\overline{w}}{\overline{w}^3},u(w,\xi)\right] +c(w,\xi)(1+u(w,\xi)),\qquad
 u(0,\xi)=0,\\
 u_0 \in C^0_{\delta'},\qquad u_3 \in C^0_{2+\delta'} \qquad \text{for} \quad \delta'>0 \quad \text{as before}.
\end{gather*}

The last two equations can be satisfied if $u\in U_{\delta'}$. On the other hand, as in the proof of \cite[Lemma 7.1]{BB04}, to solve the first equation is enough to find a fixed point of the map $\widetilde{T}_{\xi}\colon U_{\delta'}\to U_{\delta'}$ given by
\begin{equation*}
 \widetilde{T}_{\xi}(u)=h^{*}_{\lambda}(\varphi)\cdot T_{\xi}\big(\big(h^{*}_{\lambda}(\varphi)\big)^{-1}\cdot[ h^{*}_{\lambda}(c(\xi))(1+u)]\big),
\end{equation*}
where ``$\cdot$'' denotes the action by conjugation. Indeed, we have that if $v$ is a fixed point, then
\begin{equation*}
 \begin{split}
 \overline{\partial}_0(v)=\overline{\partial}_0\bigl(\widetilde{T}_{\xi}(v)\bigr)=\left[\xi H \frac{{\rm d}\overline{w}}{\lambda^2\overline{w}^3},v\right]+h^*_{\lambda}(c)(1+v).
 \end{split}
\end{equation*}
Now since $\overline{\partial}_0$ is invariant under rescaling, if we put $\widetilde{w}=\lambda w$, then $u(\widetilde{w},\xi)=v(\widetilde{w}/\lambda,\xi)$ satisfies
\begin{equation*}
 \overline{\partial}_0(u(\widetilde{w},\xi))=\left[\xi H \frac{{\rm d}\overline{\widetilde{w}}}{\overline{\widetilde{w}}^3},u(\widetilde{w},\xi)\right]+c(\widetilde{w},\xi)(1+u(\widetilde{w},\xi)),
\end{equation*}
and hence $1+u$ is the solution we seek. Now to find the fixed point, we need to show that $\widetilde{T}_{\xi}$ is a contraction. We have
\begin{align*}
 \big|\widetilde{T}_{\xi}(u)-\widetilde{T}_{\xi}(v)\big|_{U_{\delta}}&\leq C|u-v|_{U_{\delta}}|h^*_{\lambda}(c(\xi))|_{A_{\delta}}\\
 &\leq C|u-v|_{U_{\delta}}\bigl(\big|h^*_{\lambda}(a_{\text{reg}})^{0,1}\big|_{A_{\delta}}+D\big|h^*_{\lambda}\bigl(\theta_{\text{reg}}^{\dagger_h}\bigr)\big|_{A_{\delta}}\bigr),
\end{align*}
where $C$ does not depend on $\xi$; and since $U(\xi_0)$ is bounded, we have the last bound with~$D$ independent of $\xi \in U(\xi_0)$. Furthermore, by the proof of \cite[Lemma 7.1]{BB04}, we have that $\smash{\big|h^*_{\lambda}(a_{\text{reg}}^{0,1})\big|_{A_{\delta}}}=\smash{\lambda^{\delta}\big|a_{\text{reg}}^{0,1}\big|_{A_{\delta}}}$ and \smash{$\big|h^*_{\lambda}\bigl(\theta^{\dagger_h}_{\text{reg}}\bigr)\big|_{A_{\delta}}=\lambda^{\delta}\big|\theta^{\dagger_h}_{\text{reg}}\big|_{A_{\delta}}$}, so that
\begin{equation*}
 \big|\widetilde{T}_{\xi}(u)-\widetilde{T}_{\xi}(v)\big|_{U_{\delta}}\leq \lambda^{\delta}C|u-v|_{U_{\delta}}\bigl(\big|a_{\text{reg}}^{0,1}\big|_{A_{\delta}}+D\big|\theta_{\text{reg}}^{\dagger_h}\big|_{A_{\delta}}\bigr).
\end{equation*}
Hence, for $\lambda$ small enough $\widetilde{T}_{\xi}$ becomes a contraction, so we can find a fixed point.

Since \smash{$\widetilde{T}_{\xi}(u)$} is a continuous function of $w$ and $\xi$, and the rightmost term in the last inequality does not depend on $\xi$, we actually get that the fixed point must be a continuous function of both variables defined on $D_{\lambda}\times U(\xi_0)$.

Finally, by shrinking $U(\xi_0)$ and $D_{\lambda}$ if necessary, we have that our solution $g_{\xi}(w)=1+u(w,\xi)$ to problem \eqref{problem} is uniformly continuous on $D_{\lambda}\times U(\xi_0)$. Hence, given $\epsilon>0$, we can find $\delta_1,\delta_2>0$ such that $|g_{\xi}(w)-g_{\xi'}(w')|<\epsilon$ as long as $|w-w'|<\delta_1$ and $|\xi-\xi'|<\delta_2$. In particular, if $|w|<\delta_1$, we have that
$
 |g_{\xi}(w)-1|=|g_{\xi}(w)-g_{\xi}(0)|<\epsilon $ for all $\xi \in U(\xi_0)$,
so $g_{\xi}(w)\to 1$ as $w\to 0$ uniformly in $\xi \in U(\xi_0)$.
\end{proof}

\section{Proof of Proposition \ref{prelimasymp}}\label{AD}

Here we prove the first case of Proposition \ref{prelimasymp}, where $\operatorname{Re}(\lambda_1(t,\xi))>\operatorname{Re}(\lambda_2(t,\xi))$, since the other case is similar. We will use the following notation
\begin{equation*}
 \lambda_i(t,\xi)= -\xi^{-1}\gamma^*\theta_{ii} - \xi \gamma^* \theta^{\dagger_h}_{ii},\qquad
 \lambda_{ij}(t,\xi)=\lambda_i(t,\xi)-\lambda_j(t,\xi),\qquad
 R= -\gamma^*A -\xi \gamma^* \theta^{\dagger_h}_{\text{od}},
\end{equation*}
where $\theta^{\dagger_h}_{\text{od}}$ denotes the off-diagonal part of $\theta^{\dagger_h}$. If $\gamma^*A_d$ denotes the diagonal part of $\gamma^*A$, then after doing a diagonal gauge transformation of the form \smash{$(\eta_1,\eta_2)\to (\eta_1,\eta_2)\cdot \exp \bigl(-\int_{a}^{t}\gamma^*A_{d}\bigr)$}, we can gauge away the diagonal part of $\gamma^*A$, so we will assume from the beginning that we are in this gauge. We then have that $R$ only has off-diagonal elements, while $\lambda_i(t,\xi)$ is still the same as before.

The first thing we want to show is that there is a continuous solution to the following integral equation
\begin{gather}
 y_i(t,\xi) = \delta_{i2}\exp \left(\int_{a}^{t}\lambda_2(\tau,\xi){\rm d}\tau \right)- \int_{t}^{\infty}\exp \left(\int_{\tau}^{t}\lambda_i(s,\xi){\rm d}s\right)R_{ij}(\tau,\xi)y_j(\tau,\xi){\rm d}\tau\label{IntEq}
\end{gather}
for $i,j=1,2$, $i\neq j$, on any interval $(a_0,\infty)$, as long as $\xi$ is restricted to lie in a small enough neighborhood of~$0$. It is easy to check that a solution of \eqref{IntEq}, gives a solution to the original flatness equation in our chosen gauge. More explicitly, if we denote \smash{$M(t)=\exp \bigl(-\int_a^t \gamma^*A_d\bigr)$}, then in the frame~$(\widetilde{\eta}_1,\widetilde{\eta}_2)=(\eta_1,\eta_2)\cdot M$ we have a solution of the form $s(t,\xi)=y_1(t,\xi)\widetilde{\eta}_1+y_2(t,\xi)\widetilde{\eta}_2$.

If we perform the change
\begin{equation*}
 y_i(t,\xi)=z_i(t,\xi)\exp \left(\int_{a}^{t}\lambda_2(\tau,\xi){\rm d}\tau \right),
\end{equation*}
we obtain the following integral equation for the $z_i(t,\xi)$:
\begin{gather}
 z_i(t,\xi)=\delta_{i2} - \int_{t}^{\infty}\exp \left(\int_{\tau}^{t}\lambda_{i2}(s,\xi){\rm d}s\right)R_{ij}(\tau,\xi)z_j(\tau,\xi){\rm d}\tau
\label{IE2}
\end{gather}
for $i,j=1,2$, $i\neq j$.

If we use the integral relation twice and change the order of integration, we get the following integral equations for the $z_i(t,\xi)$:
\begin{align}
 z_1(t,\xi)={}&-\int_{t}^{\infty}\exp \left(\int_{t}^{\tau}\lambda_{21}\right)R_{12}(\tau){\rm d}\tau\nonumber\\
 &+ \int_{t}^{\infty}\left(\int_{t}^{s}\exp \left(\int_{t}^{\tau}\lambda_{21}\right)R_{12}(\tau){\rm d}\tau \right) R_{21}(s)z_1(s){\rm d}s\nonumber\\
:={}& \delta(t,\xi) + \int_{t}^{\infty}\epsilon_1(t,s,\xi)z_1(s,\xi){\rm d}s\nonumber\\
 z_2(t,\xi)={}& 1 + \int_{t}^{\infty}\left(\int_{t}^{s}R_{21}(\tau)\exp \left(\int_{\tau}^{s}\lambda_{21} \right){\rm d}\tau \right)R_{12}(s)z_2(s){\rm d}s\nonumber\\
:={}& 1 + \int_{t}^{\infty}\epsilon_2(t,s,\xi)z_2(s,\xi){\rm d}s,\label{IE3}
\end{align}
where we have defined
\begin{gather*}
 \delta(t,\xi):= -\int_{t}^{\infty}\exp \left(\int_{t}^{\tau}\lambda_{21}\right)R_{12}(\tau){\rm d}\tau,\\
 \epsilon_1(t,s,\xi):= \left(\int_{t}^{s}\exp \left(\int_{t}^{\tau}\lambda_{21}\right)R_{12}(\tau){\rm d}\tau \right) R_{21}(s),\\
 \epsilon_2(t,s,\xi):= \left(\int_{t}^{s}\exp \left(\int_{\tau}^{s}\lambda_{21} \right)R_{21}(\tau){\rm d}\tau \right)R_{12}(s).
\end{gather*}

Before showing that there is a solution $z_i(t,\xi)$ for the integral equations obtained above, we will say a few things about the functions $\delta(t,\xi)$ and $\epsilon_i(t,s,\xi)$. We will use the following notation: $\lambda_{21}(t,\xi)=\xi^{-1}\lambda_{21,\xi^{-1}}+\xi\lambda_{21,\xi}$, where $\lambda_{21,\xi^{-1}}=\gamma^*(-\theta_{22}+\theta_{11})$ and \smash{$\lambda_{21,\xi}=\gamma^*\bigl(-\theta^{\dagger_h}_{22}+\theta^{\dagger_h}_{11}\bigr)$}.

Notice that by our choice of gauge and WKB path, the term $\lambda_{21,\xi^{-1}}$ is constant. Hence, after integration by parts, we find the following expression for $\delta(t,\xi)$:
\begin{align}
 \delta(t,\xi)={}&\frac{\xi}{{\lambda}_{21,\xi^{-1}}}R_{12}(t,\xi)+ \frac{\xi}{\lambda_{21,\xi^{-1}}} \int_{t}^{\infty}\exp \left(\int_{t}^{\tau}\xi^{-1}\lambda_{21,\xi^{-1}}{\rm d}s\right)\frac{{\rm d}}{{\rm d}\tau} \nonumber\\
 &\times\left(\exp \left(\int_{t}^{\tau}\xi\lambda_{21,\xi}(s){\rm d}s\right)R_{12}(\tau,\xi)\right){\rm d}\tau,\label{est1}
\end{align}
where we used the fact that the elements of $R$ go to $0$ (exponentially fast) as $t\to \infty$.
Hence, using this exponential decrease of the terms of $R$, we conclude that $|\delta(t,\xi)|\leq |\xi|f(t,\xi)$, where $f(t,\xi) \to 0$ as $t \to \infty$ uniformly in $\xi$ for small enough $\xi$ (restricted to the corresponding half-plane~$\mathbb{H}_m$). Furthermore, for fixed $t$, we have $\delta(t,\xi)\to 0$ as $\xi \to 0$ with $\xi \in \mathbb{H}_m$.

Similarly, we have the following expressions for $\epsilon_i$ after integration by parts:
\begin{align}
 \epsilon_1(t,s,\xi)={}& \xi \exp \left(\int_{t}^{\tau}\lambda_{21}\right)\frac{R_{12}(\tau,\xi)}{\lambda_{21,\xi^{-1}}}\Bigg|_{\tau=t}^{\tau=s}R_{21}(s,\xi)\nonumber\\
 &- \xi \frac{R_{21}(s,\xi)}{\lambda_{21,\xi^{-1}}} \int_{t}^{s}\exp \left(\int_{t}^{\tau}\xi^{-1}\lambda_{21,\xi^{-1}}{\rm d}s\right)\frac{{\rm d}}{{\rm d}\tau}\nonumber\\
 &\phantom{- }{}\times\left(\exp \left(\int_{t}^{\tau}\xi\lambda_{21,\xi}(s){\rm d}s\right)R_{12}(\tau,\xi)\right){\rm d}\tau\nonumber\\
 \epsilon_2(t,s,\xi)={}& -\xi \exp \left(\int_{\tau}^{s}\lambda_{21}\right)\frac{R_{21}(\tau,\xi)}{\lambda_{21,\xi^{-1}}}\Bigg|_{\tau=t}^{\tau=s}R_{12}(s,\xi)\nonumber\\
 &+ \xi \frac{R_{12}(s,\xi) }{\lambda_{21,\xi^{-1}}}\int_{t}^{s}\exp \left(\int_{\tau}^{s}\xi^{-1}\lambda_{21,\xi^{-1}}{\rm d}s\right)\frac{{\rm d}}{{\rm d}\tau}\nonumber\\
 &\phantom{+ }{}\times\left(\exp \left(\int_{\tau}^{s}\xi\lambda_{21,\xi}(s){\rm d}s\right)R_{21}(\tau,\xi)\right){\rm d}\tau.\label{est2}
\end{align}
In the expressions above, we always assume $t\leq s$, since this is the range of interest for the problem. We can conclude that $|\epsilon_i(t,s,\xi)|\leq |\xi|g_i(t,\xi)$, where $g_i(t,\xi)\to 0$ as $t \to \infty$ uniformly in $\xi$ for small enough $\xi$. For this statement we use again the fact that the components of $R$ go to $0$ exponentially fast as $t \to \infty$. Furthermore, for fixed $t$, we have $\epsilon_i(t,\xi)\to 0$ as $\xi \to 0$ with~${\xi \in \mathbb{H}_m}$.

Now we are ready to prove that the integral equation $\eqref{IE3}$ has a continuous solution on the interval $[a_0,\infty)$. Because of the expressions given above for $\epsilon_i$, by restricting the $\xi \in \mathbb{H}_m$ to lie in a small neighborhood $U_{0}$ of $\xi=0$, we can ensure that for some $c\in (0,1)$
\begin{equation*}
 \int_{t}^{\infty}|\epsilon_i(t,s,\xi)|{\rm d}s<c<1 \qquad \text{for all} \quad t\in [a_0,\infty) \quad \text{and} \quad \xi \in U_{0}.
\end{equation*}

Now let
\[
z^{(0)}_2(t,\xi):=1 , \qquad z^{(0)}_1(t,\xi):=\delta(t,\xi)= -\int_{t}^{\infty}\exp \bigl(\int_{t}^{\tau}\lambda_{21}\bigr)R_{12}(\tau,\xi){\rm d}\tau.
\]
 We define \smash{$z^{(m)}_i$} recursively by plugging \smash{$z^{(m-1)}_i$} into the right side of the integral equation \eqref{IE3}; at each step we get a continuous and bounded function for $(t,\xi) \in [a_0,\infty)\times U_{0}$.

For some $M>0$, we clearly have \smash{$\big|z^{(0)}_i\big|_{C^0}<M$}, where $|\ |_{C^0}$ denotes the uniform norm on~${C^0([t_0,\infty)\times U_{0})}$. Furthermore, assume inductively that we have shown that $\smash{\big|z^{(m)}_i \!-\! z^{(m-1)}_i\big|_{C^0}}\allowbreak<c^m M$. We then have that
\begin{equation*}
 \big|z^{(m+1)}_i - z^{(m)}_i\big|_{C^0}< \bigg|\int_{t}^{\infty}|\epsilon_i(t,s,\xi)|{\rm d}s\bigg|_{C^0}\big|z^{(m)}_i - z^{(m-1)}_i\big|_{C^0}<c^{m+1}M.
\end{equation*}

From this, we see that $z_i^{(m)}$ converges uniformly in $(t,\xi)\in [a_0,\infty)\times U_{0}$ to a function $z_i \in C^0([a_0,\infty)\times U_{0})\cap L^{\infty}([a_0,\infty)\times U_{0})$. The $z_i$ clearly satisfy the integral equation \eqref{IE3}. We claim that it is also a solution of the integral equation \eqref{IE2}.

We are trying to find a solution of an integral equation of the form
$
 z= e + T(z)
$
where $e_i=\delta_{i2}$ and $T$ is the linear integral operator part of the integral equation \eqref{IE2}. We have found a solution to the integral equation $z=e+T(e)+T^2(z)$. This last solution must be unique, since if $z$ and $z'$ are solutions, then $|z-z'|_{C^0}=|T^2(z-z')|_{C^0}<\delta |z-z'|_{C^0}$, which implies that~${z=z'}$. Now notice that
\begin{equation*}
 e+T(e)+T^2(e+T(z))=e+T\bigl(e+T(e)+T^2(z)\bigr)=e+T(z),
\end{equation*}
so $e+T(z)$ is also a solution of \eqref{IE3}. By the uniqueness, we conclude that $z=e+T(z)$, so that $z$ also solves \eqref{IE2}.

With the same notation from above, notice that \eqref{IE3}, \eqref{est1} and \eqref{est2} allow us to conclude that for fixed $t\in [a_0,\infty)$ we have $T(z) \to 0$ as $\xi \to 0$; and furthermore, we have that $T(z) \to 0$ as $t\to \infty$ uniformly in $\xi \in U_{0}$.

Finally, recalling the statement and notations of Proposition \ref{prelimasymp}. Setting
\[
E_1(t,\xi):= M_{11}(t) z_1(t,\xi) \qquad \text{and} \qquad E_2(t,\xi):=M_{22}(t)(z_2(t,\xi)-1)
\]
 and using that $M(t)$ is bounded in $t\in [a_0,\infty)$ (see~\eqref{finitenessint}), we conclude the result of Proposition \ref{prelimasymp}.

\section{Proof of Lemma \ref{UHB}} \label{af}

In this appendix, we prove the following lemma.

\begin{Lemma}\label{mainlemma}
Let $m\in \mathbb{C}$ and $m^{(3)}\in \bigl(-\frac12,\frac12\bigr] \subset \mathbb{R}$. Then up to equivalence, there is a unique polystable filtered Higgs bundle $(\mathcal{E}_*,\theta)\to \bigl(\CP,\infty\bigr)$ with $\textnormal{Tr}(\theta)=0$, $\textnormal{Det}(\theta)=-\bigl(z^2+2m\bigr){\rm d}z^2$, $\textnormal{pdeg}(\mathcal{E}_*)=0$, and parabolic weights determined by $m^{(3)}$ as follows:
\begin{itemize}\itemsep=0pt
 \item if $m^{(3)} \in \bigl(-\frac12,\frac12\bigr)$, then for the eigenline decomposition near $\infty$ of the induced $\frac12$-parabolic Higgs bundle $(\mathcal{E}_{1/2},\theta)$, we have that $\pm m^{(3)}$ is the weight associated to the line corresponding to the eigenvalue $\pm\bigl(z+\frac{m}{z}+\cdots \bigr)$.
 \item if $m^{(3)}=\frac12$ then the parabolic structure of the induced $\frac12$-parabolic structure $(\mathcal{E}_{1/2},\theta)$ is the trivial filtration with weight~$\frac12$.
\end{itemize}
\end{Lemma}

Before we give the proof, we will need some notation and another lemma.

Recall that a compatibly framed connection $(E,\nabla,
\tau)\to \bigl(\CP,\infty\bigr)$ determines uniquely a~formal type $(Q,\Lambda)$ (recall Definition \ref{defFT} and Lemma \ref{FGT}). If we denote by $U_{\pm}\subset {\rm GL}(2,\mathbb{C})$ the upper (resp.\ lower) unipotent matrices, and we fix a formal type $(Q,\Lambda)$, we will denote
 \begin{align*}
 \mathcal{S}(Q,\Lambda):= \bigl\{&(S_1,S_2,S_3,S_4)\in (U_{-}\times U_{+})^2 \mid \\
 & \text{Stokes matrices of} \ (E,\nabla,
 \tau) \ \text{with formal type} \ (Q,\Lambda)\bigr\}.
 \end{align*}

Now if $(E,\nabla,
\tau)\to \bigl(\CP,\infty\bigr)$ has formal type $(Q,\Lambda)$, notice that if $T\subset {\rm GL}(2,\mathbb{C})$ denotes the set of diagonal matrices, then the set of possible compatible frames $\tau'$ for $(E,\nabla,\tau')\to \bigl(\CP,\infty\bigr)$ with formal type $(Q,\Lambda)$ is a $T$-torsor, where $t\in T$ acts on the framing in the obvious way. Furthermore, if $(E,\nabla,\tau)\to \bigl(\CP,\infty\bigr)$ has formal type $(Q,\Lambda)$ and we act on $\tau$ by $t\in T$, then the corresponding Stokes matrices in $\mathcal{S}(Q,\Lambda)$ get conjugated by $t$. Hence, we get a $T$-action on~${\mathcal{S}(Q,\Lambda)}$, and we will denote the orbits by $\mathcal{S}(Q,\Lambda)/T$. We then have the following lemma.

\begin{Lemma} \label{secondlemma} For any $m\in \mathbb{C}$ and $m^{(3)}\in \bigl(-\frac12,\frac12\bigr]$, consider $\smash{\bigl[\bigl(E_i,\overline{\partial}_{E_i},\theta_i,h_i,g_i\bigr)\bigr]}\in \mathfrak{X}^{\textnormal{fr}}\bigl(m,\allowbreak\smash{ m^{(3)}}\bigr)$ for $i=1,2$. Furthermore, let \smash{$\bigl(\mathcal{P}^h_{c}\mathcal{E}^{\xi}_i,\nabla^{\xi}_i, \tau_{c,i}^{\xi}\bigr)\to \bigl(\CP,\infty\bigr)$} be the associated compatibly framed $c$-parabolic bundles for some fixed $\xi \in \mathbb{C}^*$ and some fixed $c\in \mathbb{R}$. If we denote \smash{$Q(\xi):=\frac{1+|\xi|^2}{\xi}\frac{H}{2w^2}$} and let $\Lambda(\xi,c)$ be as in \eqref{expformon}, then
\begin{itemize}\itemsep=0pt
 \item The \smash{$\bigl(\mathcal{P}^h_{c}\mathcal{E}^{\xi}_i,\nabla^{\xi}_i, \tau_{c,i}^{\xi}\bigr)\to \bigl(\CP,\infty\bigr)$} for $i=1,2$ have Stokes data in $\mathcal{S}(Q(\xi),\Lambda(\xi,a))$.
 \item The \smash{$\bigl(\mathcal{P}^h_{c}\mathcal{E}^{\xi}_i,\nabla^{\xi}_i, \tau_{c,i}^{\xi}\bigr)\to \bigl(\CP,\infty\bigr)$} for $i=1,2$ are isomorphic as compatibly framed $c$-parabolic flat bundles if and only they have the same Stokes matrices.
 \item The \smash{$\bigl(\mathcal{P}^h_{c}\mathcal{E}^{\xi}_i,\nabla^{\xi}_i\bigr)\to \bigl(\CP,\infty\bigr)$} for $i=1,2$ are isomorphic as $c$-parabolic flat bundles if and only if, after taking some compatible framings specifying the formal type $(Q(\xi),\Lambda(\xi,a))$, their Stokes data lies in the same $T$-orbit of $\mathcal{S}(Q(\xi),\Lambda(\xi,c))$.
\end{itemize}
\end{Lemma}

\begin{proof}
The fact that the Stokes matrices of \smash{$\bigl(\mathcal{P}^h_{c}\mathcal{E}^{\xi}_i,\nabla^{\xi}_i, \tau_{c,i}^{\xi}\bigr)$} lie in $\mathcal{S}(Q(\xi),\Lambda(\xi,c))$ follows from~\eqref{singform}.

Now suppose that \smash{$\bigl(\mathcal{P}^h_{c}\mathcal{E}^{\xi}_i,\nabla^{\xi}_i, \tau_{c,i}^{\xi}\bigr)$} have the same Stokes matrices, which we denote by $S_j$ for $j=1,2,3,4$, following the conventions of Section \ref{deftwistcoords}. Let \smash{$\Phi_{j,i}$} denote the corresponding sectorial frames of flat sections for $\bigl(\mathcal{P}^h_{c}\mathcal{E}^{\xi}_i,\nabla^{\xi}_i, \tau_{c,i}^{\xi}\bigr)\to \bigl(\CP,\infty\bigr)$ defined on $\widehat{\operatorname{Sect}}_j(\xi)$, and define~\smash{$T_j\colon \mathcal{P}^h_{c}\mathcal{E}^{\xi}_1|_{\widehat{\operatorname{Sect}}_j(\xi)}\to \mathcal{P}^h_{c}\mathcal{E}^{\xi}_2|_{\widehat{\operatorname{Sect}}_j(\xi)}$} by $T_j(\Phi_{j,1})=\Phi_{j,2}$. The fact that the Stokes matrices and the formal monodromy are the same implies that the $T_j$ glue into a covariantly constant morphism $T$ over a punctured neighborhood $U_{\infty}^*$ of $\infty$, and by parallel transport, we get a~covariantly constant morphism \smash{$T\colon\mathcal{P}^h_{c}\mathcal{E}^{\xi}_1|_{\CP \setminus \{\infty\}}\to \mathcal{P}^h_{c}\mathcal{E}^{\xi}_2|_{\CP \setminus \{\infty\}}$}. Furthermore, by the expressions~\smash{${\Phi_{j,i}=\tau_{a,i}^{\xi}\cdot \Sigma_j\bigl(\widehat{F}_i\bigr)w^{-\Lambda(\xi)}e^{-Q(\xi)}}$} \big(where we abuse notation and denote by \smash{$\tau_{a,i}^{\xi}$} any local extension of the framing at $\infty$ given by \smash{$\tau_{a,i}^{\xi}$}\big), we see that for any $j=1,2,3,4$
\begin{equation*}
 T(\tau_{c,1}^{\xi})=\tau_{c,2}^{\xi}\cdot \Sigma_j\bigl(\widehat{F}_2\bigr)\Sigma_j\bigl(\widehat{F}_1\bigr)^{-1}\qquad \text{on}\quad \widehat{\operatorname{Sect}}_j(\xi),
\end{equation*}
so the $\Sigma_j\bigl(\widehat{F}_2\bigr)\Sigma_j\bigl(\widehat{F}_1\bigr)^{-1}$ glue together in $U_{\infty}^*$. Since $\Sigma_j\bigl(\widehat{F}_i\bigr)\to 1$ as $w\to 0$ for all $j=1,2, 3,4$ and~${i\!=\!1,2}$, we see that $T$ extends over the puncture to a morphism satisfying \smash{$T\bigl(\tau_{c,1}^{\xi}|_{\infty}\bigr)=\tau_{c,2}^{\xi}|_{\infty}$}. Finally, since the parabolic structures of $\bigl(\mathcal{P}^h_{c}\mathcal{E}^{\xi}_i,\nabla^{\xi}_i, \tau_{c,i}^{\xi}\bigr)$ are compatible with their irregular decompositions (see equation \eqref{irrdec}), it is easy to check that $T$ preserves the parabolic structures. Hence, $T$ gives an isomorphism between \smash{$\bigl(\mathcal{P}^h_{c}\mathcal{E}^{\xi}_i,\nabla^{\xi}_i, \tau_{c,i}^{\xi}\bigr)$} for $i=1,2$ as compatibly framed $c$-parabolic flat bundles. The other implication is trivial.

Going now to the last statement, assume that after picking compatible frames, \smash{$\bigl(\mathcal{P}^h_{c}\mathcal{E}^{\xi}_i,\nabla^{\xi}_i, \tau_{c,i}^{\xi}\bigr)$} have Stokes matrices in the same $T$-orbit of $\mathcal{S}(Q(\xi),\Lambda(\xi,c))$. Then by the previous result, for some $t\in T$, we have that \smash{$\bigl(\mathcal{P}^h_{c}\mathcal{E}^{\xi}_1,\nabla^{\xi}_1, \tau_{c,1}^{\xi}\bigr)$} is isomorphic to \smash{$\bigl(\mathcal{P}^h_{c}\mathcal{E}^{\xi}_2,\nabla^{\xi}_2, \tau_{c,2}^{\xi}\cdot t\bigr)$}. Hence, by the previous argument, we get that they are isomorphic as $c$-parabolic flat bundles (after forgetting about the framing). The other implication also follows trivially.
\end{proof}

Now we use the previous lemma to prove Lemma \ref{mainlemma}.
\begin{proof}[Proof of Lemma \ref{mainlemma}]
We divide the proof in three cases.

$m\neq 0$: We start by picking $\xi \in \mathbb{H}_m$. Notice that in the case $m\neq0$, we have that $(\mathcal{E}_*,\theta)$ must be stable (by the same argument given in Lemma \ref{C1}), so by the wild non-abelian Hodge correspondence from \cite{BB04}, we get a bijective correspondence between equivalence classes of the filtered Higgs bundles $(\mathcal{E}_*,\theta)$ we wish to count, and equivalence classes of the associated filtered flat bundles $\bigl(\mathcal{P}^h_{*}\mathcal{E}^{\xi},\nabla^{\xi}\bigr)$ for some fixed $\xi \in \mathbb{C}^*$. By taking the associated $c$-parabolic flat bundles and applying Lemma \ref{secondlemma}, we can then obtain an injection of the set we wish to count into the orbits $\mathcal{S}(Q(\xi),\Lambda(\xi,c))/T$. Now let $a(\xi)$ and~$b(\xi)$ be the non-trivial off-diagonal elements of the unipotent matrices $S_1$, $S_2$, respectively. For~${m\neq 0}$ and $\xi \in \mathbb{H}_m$, we know that after taking the compatible framings \smash{$\tau_{c}^{\xi}$}, $a(\xi)\neq 0$ for any of our~${\bigl(\mathcal{P}^h_{c}\mathcal{E}^{\xi},\nabla^{\xi}\bigr)}$ (see Proposition \ref{nonvan}). Since all the points of $\mathcal{S}(Q(\xi),\Lambda(\xi,c))$ with~${a(\xi)\neq 0}$ lie in the same $T$-orbit (recall the Stokes relations in \eqref{Stokes relations}), we conclude that there is only one $\bigl(\mathcal{P}^h_{c}\mathcal{E}^{\xi},\nabla^{\xi}\bigr)$ up to equivalence, and hence only one of the original $(\mathcal{E}_*,\theta)$ we started with, up to equivalence.

$m=0$, $m^{(3)}\neq 0$: Let first check that we cannot have a strictly polystable filtered Higgs bundle in this case. If $(\mathcal{E}_*,\theta)$ is strictly polystable, then $(\mathcal{E}_*,\theta)=(\mathcal{E}_{*,1},\theta_1)\oplus(\mathcal{E}_{*,2},\theta_2)$, with $\operatorname{pdeg}(\mathcal{E}_{*,1})=\operatorname{pdeg}(\mathcal{E}_{*,2})=0$, and it is easy to check that this cannot occur unless $m^{(3)}=0$. Hence, if $m^{(3)}\neq 0$, all our corresponding $(\mathcal{E}_*,\theta)$ are stable, and by $\cite{BB04}$, they are in bijective correspondence with equivalence classes of the associated flat filtered bundles~${\bigl(\mathcal{P}^h_{*}\mathcal{E}^{\xi},\nabla^{\xi}\bigr)}$ (for some fixed $\xi \in \mathbb{C}^*$). Notice that in this case, the formal monodromy of all of our elements turns out to be given by \smash{$M_0=\exp (-2\pi {\rm i} m^{(3)}H)\neq 1$}, so by the Stokes relations~\eqref{Stokes relations} we conclude that $\mathcal{S}(Q(\xi),\Lambda(\xi,c))/T$ is just a point, so we are done by Lemma \ref{secondlemma}.

$m=m^{(3)}=0$: in this case, we have $M_0=1$, so by the relations \ref{Stokes relations} we have that $\mathcal{S}(Q(\xi),\Lambda(\xi,c))/T$ consists of $3$ points, depending on whether $a\neq 0$ and $b=0$, $b\neq 0$ and~${a=0}$, or $a=b=0$. The case $a=b=0$ corresponds to trivial Stokes data, and it is easy to check that the filtered Higgs bundle induced from Example \ref{trivWHB} gives rise to this case. Furthermore, it is also easy to check that this is the only possible strictly polystable filtered Higgs bundle with $m=m^{(3)}=0$. Hence, the problem reduces to showing that there are no stable filtered Higgs bundles with $m=m^{(3)}=0$, giving rise to either the case with $a\neq 0$ and $b=0$, or the case $a=0$ and $b\neq 0$. Again by the Stokes relations \eqref{Stokes relations}, the case~${a=0}$ and $b\neq 0$ (resp.\ $a\neq 0$ and $b=0$) gives rise to purely upper-triangular (resp.\ lower triangular) Stokes matrices, and hence to non-stable ``Stokes representations'' (see~\cite{B3}). By the remarks in \cite[p.\ 50]{B3}, we conclude that these upper-triangular (resp.\ lower-triangular) cases cannot correspond to stable filtered Higgs bundles with $m=m^{(3)}=0$. Hence, there is only one polystable filtered Higgs bundle with $m=m^{(3)}=0$.
\end{proof}

\section[Constructing polystable parabolic Higgs bundles in the case m=0]{Constructing polystable parabolic Higgs bundles\\ in the case $\boldsymbol{m=0}$}\label{ae}

Here we explain how to construct polystable $0$-parabolic Higgs bundles whose Higgs field $\theta$ satisfies $\operatorname{Tr}(\theta)=0$ and $\operatorname{Det}(\theta)=-z^2{\rm d}z^2$.

The case where the parabolic structure consists of the trivial filtration with parabolic weights equal to $0$ is already explained in Example \ref{trivWHB}. In this case, the parabolic Higgs bundle that we find is polystable.

The next proposition deals with the rest of the possible parabolic structures.

\begin{Proposition}\label{E1} Let $m^{(3)} \in (-1,0)$. There is a stable $0$-parabolic Higgs bundle $\smash{\bigl(E^{m^{(3)}},\theta\bigr)} \to \bigl(\CP,\infty\bigr)$ with \smash{$\textnormal{pdeg}\bigl(E^{m^{(3)}}\bigr)=0$}, parabolic weights specified by \smash{$m^{(3)}$} and \smash{$-1-m^{(3)}$}, $\textnormal{Tr}(\theta)=0$, and $\textnormal{Det}(\theta)=-z^2{\rm d}z^2$.
\end{Proposition}
\begin{proof}

We start by considering $E=\mathcal{O}\oplus \mathcal{O}(-1) \to \CP$. We denote by $e_1$ and $e_2$ the usual frames over $\mathbb{C}\subset \CP$ of $\mathcal{O}$ and $\mathcal{O}(-1)$, respectively. We have that $e_1$ gives a global trivilization of $\mathcal{O}$, while the usual frame of $\mathcal{O}(-1)$ over $\CP\setminus\{0\}$ will be denoted by $f_2$. Hence, if $z$ denotes the coordinate of $\mathbb{C} \subset \CP$, then $z^{-1}e_2=f_2$.

In the frame $(e_1,e_2)$ over $\mathbb{C}\subset \CP$, we define
\begin{equation*}
 \theta=\begin{bmatrix} z & 0 \\
 2 & -z \\ \end{bmatrix}{\rm d}z.
\end{equation*}
Over $\mathbb{C}$, we can write eigenvectors corresponding to the eigenvalues $z$ and $-z$, respectively by~${v_{z}= ze_1+e_2}$ and $v_{-z}=e_2$. On the other hand, over $\CP \setminus \{z=0\}$ we can write $\widetilde{v}_z=e_1+f_2$ and~${\widetilde{v}_{-z}=f_2}$. These are eigenvectors of $\theta$ for $z\in \mathbb{C}^*\subset \CP$, with eigenvalues~$z$ and~$-z$, respectively.

We put a parabolic structure at $\infty$ by putting the weight $m^{(3)} \in (-1,0)$ to the line generated by $e_1+f_2|_{\infty}$ and the weight $-1-m^{(3)} \in (-1,0)$ to the line generated by $f_2|_{\infty}$. Denote $E$ with this parabolic structure by \smash{$E^{m^{(3)}}$}. Notice that with this parabolic structure, we have that
\begin{equation*}
 \operatorname{pdeg}\bigl(E^{m^{(3)}}\bigr)=\operatorname{deg}(E) -m^{(3)}-\bigl(-1-m^{(3)}\bigr)=0.
\end{equation*}
We claim that \smash{$\bigl(E^{m^{(3)}},\theta\bigr)\to \bigl(\CP,\infty\bigr)$} is a stable $0$-parabolic Higgs bundle. To check this, notice that on $\mathbb{C}^*\subset \CP$, we have the following relations:
\begin{gather*}
 z^{-1}v_z=\widetilde{v}_z,\qquad
 z^{-1}v_{-z}=\widetilde{v}_{-z}.
\end{gather*}
Hence, $v_z$ and $\widetilde{v}_z$ define a holomorphic line $L_{z}\cong \mathcal{O}(-1)$ while $v_{-z}$ and $\widetilde{v}_{-z}$ define $L_{-z}\cong \mathcal{O}(-1)$ ($L_{-z}$ is the same $\mathcal{O}(-1)$ summand from the splitting in the definition of $E$). These line bundles are of course preserved by $\theta$, and they are the only line bundles that can be preserved by $\theta$. We denote by $L^{m^{(3)}}_{\pm z}$ the line bundles $L_{\pm z}$ with the induced parabolic structure from $E^{m^{(3)}}$. We then have that
\begin{gather*}
 \operatorname{pdeg}\bigl(L^{m^{(3)}}_{z}\bigr)=\operatorname{deg}(L_{z}) -m^{(3)}=-1-m^{(3)}<0,\\
 \operatorname{pdeg}\bigl(L^{m^{(3)}}_{-z}\bigr)=\operatorname{deg}(L_{-z}) -\bigl(-1-m^{(3)}\bigr)=m^{(3)}<0,
\end{gather*}
so we conclude that $\bigl(E^{m^{(3)}},\theta\bigr)$ is a stable 0-parabolic Higgs bundle with parabolic degree~0.
\end{proof}

Finally, recalling that the role of $m^{(3)}$ in determining the filtered structure of the Higgs bundles is periodic mod $1$, we conclude:

\begin{Corollary} For every $m^{(3)} \in \bigl(-\frac12,\frac12\bigr]$, the set $\mathfrak{X}^{{\rm fr}}\bigl(0,m^{(3)}\bigr)$ is not empty.
\end{Corollary}

\begin{proof}
The case with $m^{(3)}=0$ is true by Example \ref{trivWHB}. For the rest of the cases, consider the filtered Higgs bundles associated to the $0$-parabolic Higgs bundles constructed in Proposition~\ref{E1}. By applying Theorem \ref{HMHB}, we obtain an adapted harmonic metric for them. Then after applying the construction of Proposition \ref{constructionframe} we obtain elements of $\mathcal{H}^{{\rm fr}}$ that define equivalence classes on $\mathfrak{X}^{{\rm fr}}\bigl(0,m^{(3)}\bigr)$ for $m^{(3)}\in \bigl(-\frac12,\frac12\bigr]\setminus\{0\}$.
\end{proof}

\subsection*{Acknowledgements}

I thank Andrew Neitzke for very helpful discussions, support, and for reading the preliminary versions of this paper. I would also like to thank the anonymous referees, who carefully read the paper and gave very useful suggestions and improvements.

\pdfbookmark[1]{References}{ref}
\LastPageEnding

\end{document}